\documentclass{article}
\usepackage{amsfonts,amssymb} 
\usepackage{latexsym}
\begin{document}
\def\dot{\!\cdot\!}
\def\mod{\mathop{\rm mod}\nolimits}
\def\Coker{\mathop{\rm Coker}\nolimits}
\def\Ker{\mathop{\rm Ker}\nolimits}
\def\Re{\mathop{\rm Re}\nolimits}
\def\Im{\mathop{\rm Im}\nolimits}
\def\sh{\mathop{\rm sh}\nolimits}
\def\ch{\mathop{\rm ch}\nolimits}
\def\rk{\mathop{\rm rk}\nolimits}
\def\Cl{\mathop{\rm Cl}\nolimits} 
\def\Gal{\mathop{\rm Gal}\nolimits}
\def\kappa{\varkappa}
\def\indlim{\mathop{{\rm ind}\,\lim}}
\def\const{\mathop{\rm const}\nolimits}
\def\Norm{\mathop{\rm Norm}\nolimits}
\def\deg{\mathop{\rm deg}\nolimits}
\def\Hilb{\mathop{\rm Hilb}\nolimits}
\def\Supp{\mathop{\rm Supp}\nolimits}
\def\BS{\mathop{\rm BS}\nolimits}
\def\upper{\mathop{\rm upper}\nolimits}
\def\lower{\mathop{\rm lower}\nolimits}
\def\unc{\mathop{\rm unc}\nolimits}
\def\GL{\mathop{\rm GL}\nolimits}
\def\ab{\mathop{\rm ab}\nolimits}

\font\mathg=eufb7 at 10pt

\title{{\bf Infinite Global  Fields and the Generalized
Brauer--Siegel Theorem} }
 
\author{{\bf M. A. Tsfasman and  S. G. Vl\u{a}du\c{t}}}
  
\date{}
 
\maketitle

\centerline{\it {To our teacher Yu.I.Manin on the occasion of his
65th birthday}}

\vskip 0.3 cm
\footnotetext {Received June 10, 2001; in revised form April 22, 2002.}
\footnotetext {Supported in part by the RFBR Grants
96-01-01378, 99-01-01204. }
\noindent {\bf Abstract.}  The paper has two purposes. First,
we start to develop a theory of infinite global fields, i.e.,
of infinite algebraic extensions either of ${\mathbb{Q}}$ or of
${\mathbb{F}}_r(t)$. We produce a series of invariants of such
fields, and we introduce and study a kind of zeta-function for
them. Second, for sequences of number fields with growing
discriminant we prove generalizations of the Odlyzko--Serre
bounds and of the Brauer--Siegel theorem, taking into account
non-archimedean places. This leads to asymptotic bounds on the
ratio ${{\log hR}/\log\sqrt{\vert D\vert}}$ valid without the
standard assumption ${n/\log\sqrt{\vert D\vert}}\rightarrow 0,$
thus including, in particular, the case of unramified towers.
Then we produce examples of class field towers, showing that this
assumption is indeed necessary for the Brauer--Siegel theorem to
hold. As an easy consequence we ameliorate on existing bounds for 
regulators. 

\vskip 0.3 cm
2000 Math. Subj. Class. 11G20, 11R37, 11R42, 14G05, 14G15, 14H05

 
\section{Introduction} 

A global field $K$ is a finite algebraic extension either of the
field ${\mathbb{Q}}$ of rational numbers, or of the field
${\mathbf{Q}}_r={\mathbb{F}}_r(t)$ of rational functions in one
variable over a finite field of constants. An {\it infinite
global field} ${\mathcal{K}}$ is either an infinite algebraic
extension of ${\mathbb{Q}}$, or such an infinite algebraic
extension of ${\mathbf{Q}}_r$ that ${\mathcal{K}}\cap\bar
{\mathbb{F}}_r={\mathbb{F}}_r$. In the first case we call
${\mathcal{K}}$ an {\it infinite number field}, in the second an
{\it infinite function field over} ${\mathbb{F}}_r$.
 
The first {\it raison d'\^etre} of our paper is an attempt to convince
ourselves and the reader that there exists a (not yet
constructed) non-trivial theory of such fields. In particular, we
produce a series of invariants, and introduce and study a
kind of zeta-function of such a field.

The second one is much more down to earth. For sequences of
number fields with growing discriminant we prove generalizations
of the Odlyzko--Serre bounds and of the Brauer--Siegel theorem,
taking into account non-archimedean places. This leads to
asymptotic bounds on the ratio ${{\log hR}/\log\sqrt{\vert
D\vert}}$ valid without the standard assumption
${n/\log\sqrt{\vert D\vert}}\rightarrow 0$, thus including, in
particular, the case of unramified towers. Then we produce
examples of class field towers, showing that this assumption is
indeed necessary for the Brauer--Siegel theorem to hold.

Wanting to study infinite global fields, we should think about
examples. For "large" fields like $\bar{\mathbb{Q}}$ or
${\mathbb{Q}}^{\ab}$ the invariants we find are trivial, but there
are numerous "smaller" ones, like the limit (i.e., the union) of
fields of a given unramified (or "not too much" ramified) tower of
fields. It is for these smaller ones that the theory we start to
develop below is interesting.

To start with, an infinite global field is always the limit of a tower
of finite ones:
$${\mathcal{K}}=\indlim_{i\to\infty}K_i=\bigcup_{i=1}^\infty K_i,
\hbox { where } K_1\subset K_2\subset K_3 \subset\dots.$$ This
tower is, of course, not unique. We are looking for invariants of
${\mathcal{K}}$, i.e., for parameters of $K_1\subset K_2\subset
K_3 \subset\dots$ that do not depend on the tower, but only on its
limit.
 
We use the following notation: Let $\{K_i\}$ for $i=1,2,\ldots$ be a
sequence of pairwise non-isomorphic global fields, either number or function;
we set $$g_i=\hbox{\rm {genus}}(K_i)$$ in the function field case, and $$g_i=\log
\sqrt{\vert D_i\vert}$$ in the number field case; we call it the {\em genus}
of a number field. 

{\it Attention:} Here and below we use the following agreement. In the 
number field case notation
$\log$ means  the natural logarithm $\log_e$. In the function
field case over ${\mathbb{F}}_r$ the same notation $\log$ means $\log_r$.
As we shall see below, this is justified by the uniformity of results obtained. 

One of the reasons to think that the definition of genus for the
number field case is natural is that for any given $g_0$ there is
only a finite number of number fields $K$ whose genus does not
exceed $g_0$. The same is true for function fields with a given
constant field.

We always doubt, whether the proper definition of genus in the
number field case should be $g=\log \sqrt{\vert D\vert}$, which
we adopt in this paper, or $g=\log \sqrt{\vert D\vert}+1$. The
latter has the advantage that, for an unramified extension, $g-1$ is
multiplied by the degree of the extension (see also
\cite{Ge/Sch}). The former one has the advantage that ${\mathbb
Q}$ is of genus $0$ and has no unramified extensions, just as a
curve of genus $0$ should. However, for infinite number fields
and other asymptotic considerations this is irrelevant, both
definitions giving same results.

We call a sequence $\{ K_i \}$ of global fields a
{\em family} if $K_i$ is non-isomorphic to $K_j$ for $i\neq j$. 
A family is called a
{\em tower} if also $K_i\subset K_{i+1}$ for any $i$. In any family
$g_i\rightarrow\infty$ for $i\rightarrow\infty$. In the function
field case we always assume that the constant field of all $K_i$
is one and the same field ${\mathbb{F}}_r$.

In the number field case let
$$n_i=[K_i:{\mathbb{Q}}]=r_1(K_i)+2r_2(K_i),$$
where $r_1$ and $r_2$ stand for the numbers of real and (pairs of) complex
embeddings.  We suppose also
that $g_i > 0$ for any $i$ , i.e., $K_i\neq {\mathbb{Q}}$ in the
number field case, and $K_i$  is non-isomorphic to  ${\mathbb{F}}_r(T)$ in the function
field case, this assumption does not restrict the generality of
our considerations.

We consider the set $A=\{{\mathbb R}, {\mathbb C}; 2, 3, 4, 5, 7,
8, 9, \dots \} $  of all prime powers plus two auxiliary symbols
${\mathbb R}$ and ${\mathbb C}$ as the set of indices. The
parameters we are going to present will be indexed by elements
$\alpha\in A$. In the function field case over ${\mathbb F_r}$
the set $A$ is reduced to $A_r=\{r, r^2, r^3, \dots \}$, meaning
that for $\alpha\in A\setminus A_r$ the parameters, we are looking
at, vanish.

For a prime power $q$ we set
$$ N_q(K_i):=  \vert\{ v \in P(K_i): \Norm(v)=q\}\vert\;,$$
where $P(K_i)$ is the set of non-archimedean places of $K_i$. We
also put $N_{\mathbb R}(K_i)=r_1(K_i)$ and
$N_{\mathbb C}(K_i)=r_2(K_i)$.

By $h_i$ we denote the class-number of $K_i$ (which equals the
number of ${\mathbb{F}}_r$-rational points on the Jacobian of
$K_i$ in the function field case); $R_i$ denotes the regulator of
$K_i$ in the  number field case and equals 1 in the function
field case.

For an infinite global field ${\mathcal K}=\bigcup K_i$ and
$\alpha\in A$ let us introduce the following quantities:

$$ \phi_\alpha=\phi_\alpha ({\mathcal K}):=\lim_{i\rightarrow \infty}{ N_\alpha(K_i) \over g_i}.$$
Of course, we need to prove that these limits exist and do not
depend on the tower.

Note that  ${\phi}_{\mathbb{R}}$ and  ${\phi}_{\mathbb{C}}$ are
finite, since the ratio $n/g_i=r_1(K_i)/g_i+2r_2(K_i)/g_i$ is
bounded on the set of all number fields $\neq {\mathbb{Q}}$ by
the Minkowski bound.

More generally, we call a  family ${\mathcal K}=\{K_i\}$,
$i=1,2,\ldots$, of global fields  {\em asymptotically exact } if
and only if for any $\alpha\in A$ there exist the limit
$$ \phi_\alpha=\phi_\alpha({\mathcal K}):=\lim_{i\rightarrow \infty}{ N_\alpha(K_i) \over g_i}.$$

We call the family ${\mathcal K}$ {\em asymptotically good} (respectively, {\em bad}) if
there exists $\alpha\in A$ with $ \phi_\alpha>0$ (respectively, $ \phi_\alpha=0$ for any
$\alpha\in A$). 

It is important to point out that, by abuse of notation, 
${\mathcal K}$ is used both for an infinite global field and
for an asymptotically exact family. This is reasonable since below
we prove that for an  infinite global field all our definitions and
results do not depend on the choice of the tower.

The notion of an asymptotically exact family is much more general
than that of a tower. In particular, a simple diagonal argument
shows that any family contains an asymptotically exact subfamily.

The quantities $\phi=\{\phi_\alpha\}$ give rise to the following
definition. The {\em limit zeta-function } of an asymptotically
exact family is defined by the product
$$\zeta_{\mathcal K}(s)=\zeta_{\phi}(s)
=\prod_{q}(1-q^{-s})^{-\phi_q}\;,$$ $q$ running over all prime
powers. 
Here and below, by raising to a complex power a function in $s$ 
defined for $\Re s > a \ge 0$ and such that its values are real positive for real $s>a$, 
we mean unique analytic continuation of what is real positive for real $s$.
The {\it "completed" zeta-function} is defined in the number field case as
$$\tilde\zeta_{\mathcal K}(s)=\tilde\zeta_{\phi}(s)
= e^s 2^{-\phi_{\mathbb{R}}}\pi^{{-s\phi_{\mathbb{R}}}/2}
(2\pi)^{-s\phi_{\mathbb{C}}}\Gamma({s\over 2})^{\phi_{\mathbb{R}}}
\Gamma(s)^{\phi_{\mathbb{C}}} \prod_{q}(1-q^{-s})^{-\phi_q} .$$
In the function field case, consistent with our convention, we set
$$\tilde\zeta_{\mathcal K}(s)=\tilde\zeta_{\phi}(s)
= r^s \prod\limits_{m=1}^\infty (1-r^{-ms})^{-\phi_{r^m}} .$$

The product defining zeta-functions $\zeta_{\phi}(s)$ and
$\tilde\zeta_{\phi}(s)$ absolutely converges for $\Re(s)\geq{1}.$ These
functions depend only on $\phi=\{\phi_\alpha\}$ and
do not depend on the particular sequence of global fields. 
Therefore, we have defined
$\zeta_{\mathcal K}(s)$ and 
$\tilde\zeta_{\mathcal K}(s)$ for any
infinite global field ${\mathcal K}$.

The zeta-function of a family is thus the limit of $g$-th roots
of usual zeta-functions of its fields $K_i$. Moreover, for $\Re
s\ge 1+\varepsilon$ the convergence is uniform. 

It is also true that a family is asymptotically exact if and only if
the limit
$\lim {\zeta_{K_i}(s)^{1/g_i}} $
exists.

Let now in the number field case
$$\xi_{\mathcal K} (s)=\xi_{\phi} (s)= (\log
\tilde\zeta_{\phi})'=\tilde\zeta_{\phi}'/\tilde\zeta_{\phi}=$$
$$1-{\phi_{\mathbb{R}}\over 2}\log\pi-\phi_{\mathbb{C}}\log{2\pi}+
{1\over 2}\phi_{\mathbb{R}}\psi({s\over 2})+
\phi_{\mathbb{C}}\psi(s)-\sum_q \phi_q{\log q\over q^s-1},$$
where $\psi(s)={\Gamma'\over\Gamma}(s)$, and in the function field case let
$$\xi_{\mathcal K} (s)=\xi_{\phi} (s)= (\log_r
\tilde\zeta_{\phi})'=1-\sum\limits_{m=1}^\infty {m\phi_{r^m}\over {r^{ms}-1}}.$$

Studying the number field case we often assume the generalized
Riemann hypothesis (GRH) to hold for number fields in question,
but the most part of our results also has an unconditional (weaker) formulation.
To distinguish between the two, we always write GRH in relevant
cases. Note that the function field case does not need it, GRH
being proved.

Under GRH the above products converge absolutely for
$\Re(s)\geq{1\over 2}$, and we have the following

{\bf GRH Theorem A} (GRH Basic Inequality). {\em  For an
infinite global field $\mathcal K$ (and for any asymptotically exact family of
global fields) 
$$\xi_{\mathcal K} ({1\over2})\ge 0.$$}

This theorem imposes severe restrictions on the possible values
of $\phi=\{\phi_\alpha\}$, namely

{\bf GRH Corollary A} (GRH Basic Inequality). {\em  For an
infinite global field (and for any asymptotically exact family of
global fields)
$$ \sum_q { \phi_q \log q \over \sqrt q-1} + \phi_{\mathbb{R}} (\log2 \sqrt
{2\pi}+{\pi \over 4} + {\gamma \over 2})+\phi_{\mathbb{C}} (\log
8\pi + \gamma) \le 1, $$ the sum being taken over all prime
powers $q$}.

In the number field case this result generalizes the GRH Odlyzko--Serre inequality on
discriminants of number fields. Indeed, all terms being
non-negative, we obtain
$$ \phi_{\mathbb{R}} (\log2 \sqrt
{2\pi}+{\pi \over 4} + {\gamma \over 2})+\phi_{\mathbb{C}} (\log
8\pi + \gamma) \le 1, $$ which means
$$D\geq (8\pi e^{\gamma+{\pi \over 2}})^
{r_1}(8\pi e^{\gamma })^{2r_2}e^{o(n)}.$$

In the function field
case over ${\mathbb F}_r$ the inequality simplifies. 
(Recall our convention that $\log$, meaning $\log_e$ in the number field case, means
 $\log_r$ in the function field case.)

{\bf Corollary A$'$} (Basic Inequality in the Function Field Case). {\em  For an
infinite function field (and for any asymptotically exact family of
function fields)
$$ \sum\limits_{m=1}^\infty { m\phi_{r^m}  \over r^{m/2}-1} \le 1. $$ }

This result generalizes the
Drinfeld--Vl\u{a}du\c{t} theorem, saying that for the number $N$ of points of degree one
on algebraic curves over a finite field ${\mathbb F}_r$ we have
$$N\le (\sqrt r -1)g+o(g),$$
as $g$ tends to infinity, i.e., that $\phi_r\le\sqrt{r}-1$.
Indeed, it is enough to omit in the sum all terms except the first one.

Contemplating the statements of Theorem A and Corollary A one gets
interested in the value of $\xi_{\mathcal K} ({1\over2})$ which equals
the {\it deficiency}, the difference between the right hand side and
the left hand side of the inequality of Corollary A. This
deficiency happens to be related to the limit distribution of
zeroes of zeta-functions.

We suppose again GRH to hold. Let ${{\mathcal K}}=\{ K_j\}$ be  an
asymptotically exact family of number fields. For each $K_j$ we
define the measure on ${\mathbb{R}}$
$$ \Delta_{K_j}:={\pi\over g_{K_j}}\sum_{\zeta_{ K_j}(\rho)=0}
\delta_{t(\rho)},$$ 
where $t(\rho)=(\rho-{ 1\over 2})/i$, and
$\rho$ runs over all non-trivial  zeroes
 of the zeta-function $\zeta_{ K_j}(s)$.  Because of GRH $t(\rho)$ is
real, and $ \Delta_{K_j}$ is a discrete measure on
${\mathbb{R}}$. Moreover, $ \Delta_{K_j}$ is a measure of slow
growth.

{\bf GRH Theorem B} (GRH Explicit Formula). {\em For an infinite 
number field ${\mathcal K}$ (and for any asymptotically exact
family of number fields)  there exists
 the limit
$$ \Delta_{\mathcal K}=\lim _{j\rightarrow \infty}\Delta_{K_j} $$
in the space of
measures of slow growth on ${\mathbb{R}}$.
 Moreover, the measure $\Delta_{\mathcal K} $ has a continuous density $M_{\mathcal K}$,
$$M_{\mathcal K}(t)=\Re \left(\xi_{\mathcal K}\left({1 \over 2}+i t\right)\right)=$$
$$1-\sum_{q}\phi_q h_q(t)\log q + $$ $${1\over 2}\phi_{\mathbb R}
 \Re\psi\left({1\over 4}+{it\over 2}\right)+\phi_{\mathbb C}
\Re\psi\left({1\over 2}+it\right)-{\phi_{\mathbb R}\over 2}
\log\pi- \phi_{\mathbb C}\log{2\pi},$$ where
$$h_q(t)={\sqrt q \cos( t\log q) -1\over q+1-2\sqrt q
\cos( t\log q)},\quad\psi(s)= {\Gamma'\over\Gamma}(s).$$}

{\bf GRH Corollary B} (GRH Basic Equality). {\em For an infinite 
number field ${\mathcal K}$ (and for any asymptotically exact
family of number fields) 
$$\xi_{\mathcal K} ({1\over 2})=M_{\mathcal K}(0),$$
i.e., 
$$ \sum_q { \phi_q \log q \over \sqrt q-1} + \phi_{\mathbb{R}} (\log2 \sqrt
{2\pi}+{\pi \over 4} + {\gamma \over 2})+\phi_{\mathbb{C}} (\log
8\pi + \gamma) = 1-M_{\mathcal K}(0). $$} 

This means that the difference between 1 and the left hand side of the Basic 
Inequality, called the {deficiency} of the infinite global field (or of the family),
is in fact the "relative number" of zeroes accumulating at the real critical point $1\over 2$.

In the function field case the same is true and much easier to prove (cf.{\cite {31}}). 
Zeta-functions being periodic,
we can of course consider the space of periodic measures on $\mathbb R$ to obtain Theorem B and
Corollary B in this case. We can also make the formulation simpler using measures on the
circle. We normalise the circle to be
${\mathbb R}/2\pi{\mathbb Z}$ represented by
$(-\pi, \pi]$. 
For a  zero  $\rho$ of the zeta-function $\zeta_{ K_j}(s)$ let $t(\rho)$ be
defined by
$$t(\rho)={{\rho-{ 1\over 2}}\over{ i}}\;(\mod 2\pi).$$

Let $$ \Delta_j:={{ \pi\over g_j}}\sum_{\zeta_{
K_j}(\rho)=0}\delta_{t(\rho)},$$
where $\delta_{t(\rho)}$ is, as usual, the Dirac measure supported at
$t(\rho)$.
Then $ \Delta_j$ is a measure of total mass $2\pi$    on
${{\mathbb{R}} / 2\pi{\mathbb{Z}}}$, and  $ \Delta_j$ is symmetric with
respect to $t\mapsto-t$. 

{\bf Corollary B$'$} (Explicit  Formula and Basic Equality in the Function Field Case). 
{\em In the function field case in the
weak topology on the space of measures on ${{\mathbb{R}} / 2\pi{\mathbb{Z}}}$   the limit
$$ \Delta_{\mathcal K}=\lim _{j\rightarrow \infty}\Delta_j $$
exists.  Moreover, the measure $\Delta_{\mathcal K} $ has a continuous density
$M_{\mathcal K}$,
$$M_{\mathcal K}(t)={\Re} (\xi_{\mathcal K}({1 \over 2}+{i\over\log_e r}
t))=1-\sum_{m=1}^{\infty}
m\phi_{r^m} h_m(t)$$
for
$$h_m(t)={ r^{m/2}\cos(mt) -1\over r^m+1-2r^{m/2}\cos(mt)}\;,$$
which depends only on the family of numbers $\phi=\{\phi_{r^m}\}$ and we have the
following
Basic Equality:
$$\xi_{\mathcal K}({1 \over 2})=1-\sum_{m=1}^{\infty}{m\phi_{r^m}\over
r^{m/2}-1}=M_{\mathcal K}(0).$$}

In the number field case with no GRH at hand, the results become considerably weaker.

{\bf Theorem C} (Unconditional Basic Inequality). {\em  For an
infinite number field ${\mathcal K}$ (and for any asymptotically exact family of
number fields) 
$$\xi_{\mathcal K} (1)\ge 0.$$}

{\bf Corollary C} (Unconditional Basic Inequality). {\em  For an
infinite number field (and for any asymptotically exact family of
number fields)
$$\sum_{q}{\phi_q\log q \over q-1}+(\gamma/2+\log 2\sqrt\pi)\phi_{\mathbb{R}}+
(\gamma+\log 2\pi) \phi_{\mathbb{C}} \leq 1.$$}

This time, omitting all non-archimedean terms, we get the
unconditional Stark inequality,
$$D\geq (4{\pi} e^{\gamma})^
{r_1}(2\pi e^{\gamma })^{2r_2}e^{o(n)}.$$

The unconditional Odlyzko inequality
$$D\geq (4\pi e^{\gamma+1})^{r_1}(4\pi e^{\gamma })^{2r_2}e^{o(n)}$$
can also be generalized, using non-archimedian places.

Our next result concerning zeta-functions concerns the behaviour
of class numbers and regulators. For an asymptotically exact
family ${\mathcal K}$ of global fields we would like to consider the limit
$$\BS({\mathcal K})=\lim_{i\rightarrow \infty}{ \log{h_i R_i} \over
g_i}.$$ 
Under certain conditions, as we shall explain below, this limit exists 
and depends only on the set of
numbers
$\phi=\{\phi_\alpha\}$. Therefore, $\BS({\mathcal K})$ is well defined
for an infinite global field ${\mathcal K}$, as well as for any 
asymptotically exact family ${\mathcal K}$. We  can also define
$${\mathbf\kappa}({\mathcal K})=\lim_{i\rightarrow \infty}
{ \log\kappa_i \over
g_i},$$
 $\kappa_i$ being the residue of $\zeta_{K_i}(s)$ at 1;
this invariant exists under the same conditions.

The value of $\BS({\mathcal K})$ is described by the 
 Brauer--Siegel theorem. 
In our terms the 
classical Brauer--Siegel theorem states:
 
{\em We have
$$\BS({\mathcal K})=1 \hbox { and } {\mathbf\kappa}({\mathcal K})=0,$$ 
if the family ${\mathcal K}$
satisfies the following two conditions:

(i) the family ${\mathcal K}$ is asymptotically bad;

(ii) either GRH holds, or all fields $K_i$ are normal over $\mathbb Q$.}

Indeed, the 
assumption ${n/\log\sqrt{\vert D\vert}}\rightarrow 0$, usually used in the statement,
 means $\phi_\alpha=0$ for all
$\alpha$, since it follows that
$\phi_{\mathbb{R}}=\phi_{\mathbb{C}}=0$ and for a prime $p$ one has
$$\sum_{m=1}^{\infty} m\phi_{p^m} \leq
\phi_{\mathbb{R}}+2\phi_{\mathbb{C}}. $$

We are going to generalize the Brauer--Siegel theorem disposing of the
first condition.

{\bf GRH Theorem D} (GRH Generalized Brauer--Siegel Theorem). {\em
For an infinite global field ${\mathcal K}$ (and for any asymptotically exact
family of global fields)  the limits $\BS({\mathcal K})$
and ${\mathbf\kappa}({\mathcal K})$ exist and we have
$$\BS({\mathcal K})=\log\tilde\zeta_{\mathcal K}(1),$$
$${\mathbf\kappa}({\mathcal K})=\log\zeta_{\mathcal K}(1).$$
}

{\bf GRH Corollary D} (GRH Generalized Brauer--Siegel Theorem). {\em
For an infinite global field (and for any asymptotically exact
family of global fields)
$$  \BS({\mathcal K}) =1+\sum_{q}\phi_q\log{ q\over  q-1}- \phi_{\mathbb{R}}\log 2-
\phi_{\mathbb{C}}\log 2\pi,$$
$${\mathbf\kappa}({\mathcal K})=\sum_{q}\phi_q\log{ q\over  q-1},$$ 
the sum being taken over all prime
powers $q$.}

In the function field case, of course,
 $$\BS({\mathcal K})=\lim_{i\to\infty}{\log_r h_i\over
g_i},$$ where $h_i$ is the number of ${\mathbb F}_r$-points on the Jacobian of the curve
$X_i$  corresponding to the field $K_i$. The other parameter ${\mathbf\kappa}({\mathcal
K})= \BS({\mathcal K})-1$, and thus becomes uninteresting.

{\bf Corollary D$'$} (Generalized Brauer--Siegel Theorem in the  Function 
Field Case). {\em
For an infinite function field ${\mathcal K}$ (and for any asymptotically exact
family of function fields)  the limit $\BS({\mathcal K})$
exists and we have
$$  \BS({\mathcal K}) =1+\sum\limits_{m=1}^\infty\phi_{r^m}\log_r{ {r^m}\over {r^m}-1}.$$
}

In the number field case half of Theorem D does not depend on GRH, namely we prove

{\bf Theorem E} (Generalized Brauer--Siegel Inequality). {\em For
an infinite number field (and for any asymptotically exact family
of number fields)
$$  \limsup_{i\rightarrow \infty}{ \log(h_iR_i) \over g_i} \leq 1+
\sum_{q}\phi_q\log{ q\over  q-1}- \phi_{\mathbb{R}}\log 2-
\phi_{\mathbb{C}}\log 2\pi,$$
$$\limsup_{i\rightarrow \infty}{ \log({\mathbf\kappa}_i) \over g_i}\le
\sum_{q}\phi_q\log{ q\over  q-1},$$
 the sum being taken over all
prime powers $q$.}

As yet, we are unable to prove the generalized Brauer--Siegel
theorem unconditionally. In the general case, even the classical
Brauer--Siegel theorem is not known, one needs normality of the
fields in question. However, for an infinite number field with an
auxiliary condition this becomes possible.

{\bf Theorem F} (Unconditional Generalized Brauer--Siegel Theorem
for Infinite Number Fields ). {\em For an infinite almost normal
asymptotically good 
number field $\mathcal K$ 
the limits $\BS({\mathcal K})$ and ${\mathbf\kappa}({\mathcal K})$
 exist and we have
$$\BS({\mathcal K})=\log\tilde\zeta_{\mathcal K}(1),$$
$${\mathbf\kappa}({\mathcal K})=\log\zeta_{\mathcal K}(1).$$
}

{\bf Corollary F} (Unconditional Generalized Brauer--Siegel
Theorem for Infinite Number Fields). 
{\em For an infinite almost normal asymptotically
good number field
$\mathcal K$ 
$$\BS({\mathcal K})=1+\sum_{q}\phi_q\log{ q\over  q-1}- \phi_{\mathbb{R}}\log 2-
\phi_{\mathbb{C}}\log 2\pi,$$
$${\mathbf\kappa}({\mathcal K})=\sum_{q}\phi_q\log{ q\over  q-1},$$ 
the sum being taken over all prime
powers $q$. }

Note that the unconditional classical Brauer--Siegel theorem for normal 
fields {\it does not} follow from our results.

Next question is that of the possible asymptotic behaviour of the
Brauer--Siegel ratios ${\log(hR)}/g$ and ${\log({\mathbf\kappa})}/g $. We prove

{\bf GRH Theorem G} (GRH Bounds). {\em  For any family of number
fields
$$   \BS_{\lower}\leq  \liminf_{i\rightarrow \infty}{ \log(h_iR_i) \over g_i}\leq
\limsup_{i\rightarrow \infty}{ \log(h_iR_i) \over g_i}\leq
\BS_{\upper},$$
$$0\leq  \liminf_{i\rightarrow \infty}{ \log({\mathbf\kappa}_i) \over g_i}\leq
\limsup_{i\rightarrow \infty}{ \log({\mathbf\kappa}_i) \over g_i}\leq
{\mathbf\kappa}_{\upper},$$
 where $$\BS_{\lower}= 1-{\log
2\pi\over\gamma+\log8\pi} \approx 0.5165...,$$
$$\BS_{\upper}= 1+{{\log{3\over 2}+\log{5\over 4}+\log{7\over 6}}\over{
{\gamma\over 2}+{\pi\over 4}+ \log{2\sqrt{2\pi}}+{\log 2\over
{\sqrt 2 -1}}+{\log 3\over {\sqrt 3 -1}} +{\log 5\over {\sqrt 5
-1}}+{\log 7\over {\sqrt 7 -1}}}} \approx 1.0938...,$$
$${\mathbf\kappa}_{\upper}= {\log 2+{\log{3\over 2}}\over{
{\gamma\over 2}+\log{2\sqrt{2\pi}}+{\log 2\over
{\sqrt 2 -1}}+{\log 3\over {\sqrt 3 -1}}}} \approx 0.2164\dots$$}

In what follows we also give some bounds specific for the totally real case and
for the totally complex one.

The function field case was treated in our paper \cite{31}. In our
terms we have

{\bf Theorem G$'$} (Function Field Bounds). {\em  For any family of
function fields over ${\mathbb F}_r$
$$   1\leq  \liminf_{i\rightarrow \infty}{ \log_r{h_i} \over g_i}\leq
\limsup_{i\rightarrow \infty}{ \log_r{h_i}  \over g_i}\leq 1+
(\sqrt r-1)\log_r{r \over r-1}.$$
}

In the number field case, as usual, without GRH Theorem G weakens.

{\bf Theorem H} (Unconditional Bounds). {\em  For any family of
number fields
$$ \limsup_{i\rightarrow \infty}{ \log(h_iR_i) \over g_i}\leq
\BS_{\unc, \upper} ,$$
$$ \limsup_{i\rightarrow \infty}{ \log({\mathbf\kappa}_i) \over g_i}\leq
{\mathbf\kappa}_{\unc, \upper} ,$$
where 
$$\BS_{\unc, \upper}=
1+{{\sum\limits_{{p=3}\atop{prime}}^{23} \log{p\over {p-1}}}\over{
{\gamma\over 2}+{1\over 2}+\log{2\sqrt{\pi}}+
2\sum\limits_{{p=2}\atop{prime}}^{23} \log p
\sum\limits_{m=1}^\infty {1\over {p^m+1}}}}
 \approx 1.1588\dots,$$ 
$${\mathbf\kappa}_{\unc, \upper}=
1+{{\sum\limits_{{p=3}\atop{prime}}^{5} \log{p\over {p-1}}}\over{
{\gamma\over 2}+\log{2\sqrt{\pi}}+
2\sum\limits_{{p=2}\atop{prime}}^{5} \log p
\sum\limits_{m=1}^\infty {1\over {p^m+1}}}}
 \approx 0.3151\dots
$$
 For an infinite almost normal asymptotically good number field ${\mathcal K}$ 
we also have the lower bound
 $$\liminf_{i\rightarrow \infty}{ \log(h_iR_i) \over g_i}\ge  \BS_{\unc, \lower},$$
  where 
$$\BS_{\unc, \lower}=1-{\log {2\pi}\over {\gamma+\log
{4\pi}}}\approx 0.4087\dots$$}

Knowing that GRH--possible values of the Brauer--Siegel ratio lie
in the interval
$$(0.5165... , 1.0938...),$$
and having in mind the classical value 1 of the Brauer--Siegel
theorem itself, we are curious to know whether in our more general
setting there exist examples when it differs from 1.

They do exist. The method to construct such examples of infinite
number fields is to take the limit of a class field tower given
by some splitting conditions. In particular, we get

{\bf GRH Theorem I.} {\em The field
$$K={\mathbb{Q}}(\cos{2\pi\over11},\sqrt{2},\sqrt{-23})$$
has an infinite unramified $2$-tower $\mathcal K$, for which 
$\BS({\mathcal K})\in(\BS_{\lower}({\mathcal K}),\BS_{\upper}({\mathcal K}))$,
where
$$\BS_{\lower}({\mathcal K})=1-{{10\log 2\pi}\over g},$$
$$\BS_{\upper}({\mathcal K})=\BS_{\lower}({\mathcal K})+{{(\sqrt{23}-1)\log{23\over 22}}\over
\log{23}}\left(1-{{10(\gamma+\log{8\pi})}\over{g}}\right),$$
i.e.,
approximately
$$0.5939\dots \le \BS({\mathcal K})\le 0.6025\dots.$$
}

Note that we do not need examples giving lower bounds for 
${\mathbf\kappa}({\mathcal K})$ since any asymptotically bad infinite 
number field, for example any tower of fields abelian over $\mathbb Q$,
attains the obvious lower bound ${\mathbf\kappa}({\mathcal K})=0$.

Without GRH the upper bound is less precise.

{\bf Theorem J.} {\em The field
$$K={\mathbb{Q}}(\cos{2\pi\over11},\sqrt{2},\sqrt{-23})$$
has an infinite unramified $2$-tower $\mathcal K$, for which 
$\BS({\mathcal K})\in(\BS_{\lower}({\mathcal K}), \BS_{\unc, \upper}({\mathcal K}))$,
where
$\BS_{\lower}({\mathcal K})$ is as above and
$\BS_{\unc, \upper}({\mathcal K})\approx 0.7108\dots$.}

The upper bound we have got shows that the condition $n/\log
\vert D\vert\rightarrow 0$ (or in our terms
 ${\phi}_{\alpha}=0$ for every $\alpha$)  in the classical
Brauer--Siegel theorem is indeed indispensable. In other words, the
Brauer--Siegel ratio $\BS({\mathcal K})$ can be strictly less than 1. Can it also be
strictly greater than 1? Can ${\mathbf\kappa}({\mathcal K})$ be strictly positive? 
Here is an example.

{\bf GRH Theorem K.} {\em The field
$$K = {\mathbb{Q}}(\sqrt{11\dot 13\dot 17\dot 19\dot 23\dot 29\dot 31\dot
37\dot 41\dot 43\dot 47\dot 53\dot 59\dot 61\dot 67})$$ has an
infinite unramified $2$-tower ${\mathcal K}$ in which nine prime
ideals lying over $2$, $3$, $5$, $7$ and $71$ split completely.
Then
$\BS({\mathcal K})\in(\BS_{\lower}({\mathcal
K}),\BS_{\upper}({\mathcal K}))$ and
${\mathbf\kappa}({\mathcal K})\in({\mathbf\kappa}_{\lower}({\mathcal
K}),{\mathbf\kappa}_{\upper}({\mathcal K}))$, where
$$\BS_{\lower}({\mathcal K})=1+{{2\log{3\over 2} +2\log{5\over 4}+2\log{7\over
6}+\log{5041\over 5040}}\over g},$$ 
$$\BS_{\upper}({\mathcal K})=\BS_{\lower}({\mathcal K})+{1\over g}\sum_{p=11}^{47}\log{p\over p-1}+$$
$${{{\sqrt{53}-1}\over{g\log 53}}
\left(g-\gamma-{\pi\over 2}-\log{8\pi}-2\sum_{p=2}^7 {{\log
p}\over{\sqrt p -1}}- {{\log 71^2}\over{70}}-\sum_{p=11}^{47}
{{\log p}\over{\sqrt p -1}}\right)}\log{{53}\over52},$$
$${\mathbf\kappa}_{\lower}({\mathcal K})={2\log 2+{2\log{3\over 2} +2\log{5\over 4}+2\log{7\over
6}+\log{5041\over 5040}}\over g},$$
$${\mathbf\kappa}_{\upper}({\mathcal K})=\BS_{\upper}({\mathcal K})-1+{2\log 2\over g},$$ the sums
being taken over prime $p$'s. Numerically
$$\BS({\mathcal K})\in(1.0602\dots , 1.0798\dots),$$
$${\mathbf\kappa}({\mathcal K})\in(0.1135\dots , 0.1331\dots).$$}

Here, as well, without GRH the upper bound changes.

{\bf Theorem L.} {\em The field
$$K = {\mathbb{Q}}(\sqrt{11\dot 13\dot 17\dot 19\dot 23\dot 29\dot 31\dot
37\dot 41\dot 43\dot 47\dot 53\dot 59\dot 61\dot 67})$$ has an
infinite unramified $2$-tower ${\mathcal K}$ in which nine prime
ideals lying over $2$, $3$, $5$, $7$ and $71$ split completely.
Then 
$\BS({\mathcal K})\in(\BS_{\lower}({\mathcal
K}),\BS_{\unc, \upper}({\mathcal K}))$  
and
${\mathbf\kappa}({\mathcal K})\in({\mathbf\kappa}_{\lower}({\mathcal
K}),{\mathbf\kappa}_{\unc, \upper}({\mathcal K}))$, where
$\BS_{\lower}({\mathcal K})$ and ${\mathbf\kappa}_{\lower}({\mathcal K})$ are as above,
$\BS_{\unc, \upper}({\mathcal K})\approx 1.0951\dots$, 
${\mathbf\kappa}_{\unc, \upper}({\mathcal K})\approx 0.1454\dots$.}

The values of different bounds for $\BS({\mathcal K})$ developed in this paper form the
following table.

$$
\mbox{\small\begin{tabular}{c|l|cccc}
\multicolumn{2}{c|}{}&lower&lower&upper&upper\\
\multicolumn{2}{c|}{}&bound&example&example&bound\\ \hline
&all fields&0.5165&0.5939--0.6025&1.0602--1.0798&1.0938\\ 
GRH&totally real&0.7419&0.8009--0.8648&1.0602--1.0798&1.0938\\ 
&totally complex&0.5165&0.5939--0.6025&1.0482--1.0653&1.0764\\ 
\hline  
&all fields&0.4087&0.5939--0.7108&1.0602--1.0921&1.1588\\ 
unconditional&totally real&0.6625&0.8009--0.9248&1.0602--1.0921&1.1588\\ 
&totally complex&0.4087&0.5939--0.7108&1.0482--1.0951&1.0965
\end{tabular}}
$$ 

And here is the table for ${\mathbf\kappa}({\mathcal K})$. 
Note that the lower bound ${\mathbf\kappa}({\mathcal K})=0$ is always attainable.

$$
\mbox{\small\begin{tabular}{c|l|cc}
\multicolumn{2}{c|}{}&upper&upper\\
\multicolumn{2}{c|}{}&example&bound\\ \hline
&all fields&0.1135--0.1331&0.2164\\
GRH&totally real&0.1135--0.1331&0.1874\\
&totally complex&0.1162--0.1333&0.2164\\
\hline 
&all fields&0.1135--0.1454&0.3151\\ 
unconditional&totally real&0.1135--0.1454&0.2816\\ 
&totally complex&0.1162--0.1631&0.3151
\end{tabular}}
$$

In the function field case an example of  ${\mathcal K}$ with $\BS({\mathcal K})=1$, 
${\mathbf\kappa}({\mathcal K})=0$ is
provided by any tower with $\phi_\alpha=0$ for every $\alpha$. In
particular, any tower of fields abelian over ${\mathbb F}_r(t)$
has this property. An example reaching the upper bound must have
$\phi_r=\sqrt r-1$ and $\phi_\alpha=0$ for every other $\alpha$.
The existence of such towers is known only when $r$ is a square.
For a square $r$, different modular towers enjoy this property.

As an application of the Generalized Brauer-Siegel Theorem 
one obtains a lower bound for  
regulators of number fields in asymptotically good families  
which is better than Zimmert's bound.

{\bf Theorem M} (Regulator Bound). {\em For  an 
asymptotically good 
tower of number fields ${\mathcal K}=\{K_i\}$ we have
$$\liminf_{i\rightarrow \infty} {\log R_i\over g_i} \ge 
(\log \sqrt{\pi e} +{\gamma\over 2})\phi_{\mathbb R}+(\log 2 +\gamma)\phi_{\mathbb C}.$$}

Under GRH we get the same estimate for any asymptotically 
good family of number fields.

\centerline{---------------------------}

Our work resulting in this paper was started about ten years ago.
Now we are convinced that there is a non-trivial theory of
infinite global fields, though we do not yet understand what it
should really look like.

\centerline{---------------------------}

The paper starts with generalities on infinite global fields and
their zeta-functions. In Section 2 we introduce the invariants
${\phi_\alpha}({\mathcal K})$. In Section 3 we prove the first
form of the Basic Inequality. Then we introduce zeta-functions
(Section 4) and prove the Explicit Formula (Section 5). In
Section 6 we discuss possible directions of further study of
infinite global fields.

Part 2 is consacrated to the Brauer--Siegel theorem. We prove the
generalized Brauer--Siegel theorem in Section 7, as well as regulator bounds.
In Section 8 we provide the
bounds for the Brauer--Siegel ratio. Section 9 is devoted to 
class field towers. We finish by
discussing open questions.

\centerline{------------------------------}

It is a pleasure for us to acknowledge the previous work without which this
paper would have never been written.
Any unified treatment of number and function fields makes appeal to the heritage
of A.Weil.  A great part of this work develops two classical results, 
the Odlyzko--Serre inequalities and the Brauer--Siegel theorem, both
ideologically and technically. We first understood what is going on in the
function field case \cite{31}. Y. Ihara \cite{9} obtained most part of
results of Section 3 below in the particular case of unramified towers, both in the
function field case and in the number field one. (Unfortunately, we were unaware
of \cite{9} while writing \cite{31}.) His technique helped us a lot. 
A version of a particular case (the asymptotically bad one) of 
GRH Theorem B is the main result of Lang's paper \cite{12}.
Discussions with many of our friends and collegues were extremely useful. 
We would especially like to thank for many 
valuable remarks J.-P.Serre and the anonymous author of a 
13 page long referee report on one of the previous versions of this paper. 
We thank G.Lachaud for his interest in
our work and for  attracting our attention to the  question about the minimum zeta-zero.

It is our greatest pleasure to devote this paper to our teacher 
Yuri Ivanovich Manin, who taught us to consider number fields, zeta functions
and algebraic curves as different facets of one diamond. 
Congratulating him with his 65th birthday, we wish him many happy returns of the
day.

\part {Zeta-function of an infinite global field}

Let us repeat the definition. An {\it infinite global field}
${\mathcal{K}}$ is either an infinite algebraic extension of
${\mathbb{Q}}$, or an infinite algebraic extension of
${\mathbf{Q}}_r={\mathbb{F}}_r(t)$  such that ${\mathcal{K}}\cap\bar
{\mathbb{F}}_r={\mathbb{F}}_r$. In the first case we call
${\mathcal{K}}$ an {\it infinite number field}, in the second an
{\it infinite function field over} ${\mathbb{F}}_r$.

Our main problem here is to find out parameters of infinite
global fields and to construct a zeta-function of such a field.

\section{Invariants of infinite global fields}

Here we give some basics on infinite global fields and asymptotically exact
families showing
that these notions are worth studying.

{\bf Lemma 2.1.} {\em For any given $g_0$ there is only a finite number of
number fields $K$ whose
genus does not exceed $g_0$. The same is true for function fields 
 over a given constant field (considered up to an
isomorphism).}

{\em Proof.} In the number field case this is proved by the geometry of
numbers
 (cf. \cite{13}, Theorem V.4.5).
In the function field case this follows from the existence of moduli spaces
for genus $g$ curves.
Indeed, those are varieties over the ground field which is finite,
and thus they have but a finite
number of points defined over it.
${\Box}$

For a global field $K$ and for a prime power $q$ we set
$$ N_q(K):=  \vert\{ v \in P(K): \Norm v =q\}\vert,$$
where $P(K)$ is the set of non-archimedean places of $K$. We
also write $N_{\mathbb R}(K)=r_1(K)$  for the number of real places and
$N_{\mathbb C}(K)=r_2(K)$ for that of complex ones. In the function field
case we set
$N_{\mathbb R}(K)=N_{\mathbb C}(K)=0$.
 
We call a sequence $\{ K_i \}$ of global fields a
{\em family} if $K_i$ is non-isomorphic to $K_j$ for $i\neq j$. In any family
$g_i=g(K_i)\rightarrow\infty$ for $i\rightarrow\infty$. (In the function
field case we always assume that the constant field of all $K_i$
equals one and the same ${\mathbb{F}}_r$.)

{\bf Definition 2.1.}  We call a  family $\{K_i\}$  {\em asymptotically
exact }
 if and only if for any
$\alpha\in A=\{{\mathbb R}, {\mathbb C}; 2, 3, 4, 5, 7,
8, 9, \dots \} $
there exists a limit
$$ \phi_\alpha:=\lim_{i\rightarrow \infty}{ N_\alpha(K_i) \over g_i}.$$

Note that   we can as well divide by $g_i-1$
instead of $g_i$, since for almost all $i$ we have $g_i>1$
and $g_i\to \infty$. We shall use this division  in the function field case.

{\bf Lemma 2.2}. {\em Every family of global fields $\{K_i\}$ contains an asymptotically
exact subfamily.}

{\em Proof.}
In the number field case, for any $q$ we have $N_q(K_i)\le n(K_i)$, where $n(K)=
[K:{\mathbb{Q}}]$. On the other hand, since $K_i\neq{\mathbb{Q}}$
there exists a universal constant $C$ such that
$\vert D_{K_i}\vert^{1/n}\ge C$. This follows, e.g., from the
Odlyzko--Serre inequalities
(or even from the Minkowski constant argument), and
Lemma 2.1. Thus ${{n(K_i)}/{g(K_i)}}\le {{1}/{\log C}}$, and all
 limit points of ${{N_q(K_i)}/{g(K_i)}}$ lie between 0 and
${{1}/{\log C}}$.

Therefore, for each separate $q$ there is a limit over some subsequence
of ${\bf K}=\{K_i\}$. The same is true for ${{r_\alpha(K_i)}/{g(K_i)}},
\; \alpha={\mathbb R}, {\mathbb C}$.

Choose such a subsequence ${\bf K}_0$ that $\phi_{\mathbb{R}}$ exists. Take
its subsequence ${\bf K} _1$ where $\phi_{\mathbb{C}}$ exists, then its
subsequence ${\bf K} _2$ where $\phi_2$ exists, and so on. 
Now take a
diagonal sequence, i.e., such that $K_1\in {\bf K}_1$, $K_2\in {\bf K}_2$, etc.
For this sequence $\phi_\alpha$ exists for any $\alpha$.

The function field case is treated similarly.
${\Box}$

{\bf Remark 2.1}. {As it was pointed out by J.-P. Serre, this proof only
uses
the fact that
the space of positive measures with mass 1 is a metric compact space.}

The lemma shows that the study of any family can be reduced in a certain
sense to that of
asymptotically exact families, and from now on we mostly suppose all our families
to be asymptotically exact.

{\bf Lemma 2.3}. {\em Let $L\supseteq K$. Then, in the number field case, for $g(K)\ge 1$, 
$${N_{\mathbb R}(L)\over{g(L)}}+2 {N_{\mathbb C}(L)\over{g(L)}}\le
{N_{\mathbb R}(K)\over{g(K)}}+2 {N_{\mathbb C}(K)\over{g(K)}},$$
and, for any prime $p$ and any $n\ge 1$,
$$\sum\limits^n_{m=1}{m N_{p^m}(L)\over{g(L)}}\le \sum\limits^n_{m=1}{m
N_{p^m}(K)\over{g(K)}}.$$
In the function field case, for $g(K)\ge 2$ and for any $n\ge 1$, 
$$\sum\limits^n_{m=1}{m N_{r^m}(L)\over{g(L)-1}}\le \sum\limits^n_{m=1}{m
N_{r^m}(K)\over{g(K)-1}}.$$
}

{\em Proof.} In the number field case, for an extension $L\supset K$ we have
$$\vert D_{L}\vert\ge\vert D_{K}\vert^{[L:K]},$$ and in the
function field case we have $$g(L)-1\ge [L:K] (g(K)-1).$$
On the other hand, if a place $v$ of $K$ is decomposed
into a set $\{v_1, v_2, \dots\}$ of places of $L$ then  $$\prod \Norm v_i \le
(\Norm v )^{[L:K]}.$$
Therefore, $$ \sum^n_{m=1}m N_{p^m}(L)\le
[L:K] \sum^n_{m=1}m N_{p^m}(K)$$
for any prime $p$ and any $n\ge 1$. Dividing, we get the required
monotonicity. 

The archimedean inequality and that for function fields are
similar.
${\Box}$

{\bf Lemma 2.4}. {\em Any infinite tower $K_0\subset K_1\subset K_2\subset
\dots$ is an
asymptotically exact family.}

{\em Proof.} For a given prime $p$, by the second inequality of Lemma 2.3, the
sequence $\sum\limits^n_{m=1}{m N_{p^m}(K_i)\over{g(K_i)}}, i=1,2,\ldots$
is non-increasing for any fixed $n$, and hence has a limit. (In the function field case we
divide by $g(K_{i})-1$.)
Taking $n=1$ we see that $\phi_p$ exists,
setting $n=2$ we derive the existence of $\phi_{p^2}$, etc.

Using the first inequality of Lemma 2.3 we prove the
existence of $\phi_{\mathbb R}$ and then of  $\phi_{\mathbb C}$.

The function field case is treated similarly.
${\Box}$

Let ${\mathcal{K}}$ be an infinite global field. Then
${\mathcal{K}}=\bigcup_{i=1}^\infty K_i$,
 where $K_1\subset K_2\subset K_3 \subset\dots$, and we can define
 the corresponding parameters $\phi_\alpha$, $\alpha\in A=\{{\mathbb R},
{\mathbb C}; 2, 3, 4, 5, 7,
8, 9, \dots \} $.

{\bf Lemma 2.5.} {\em For an infinite global field ${\mathcal{K}}$ the
parameters $\phi_\alpha$ do not depend on the choice of the tower
$K_1\subset K_2\subset K_3
\subset\dots$}.

{\em Proof.} Let ${\mathcal{K}}=\bigcup_{i=1}^\infty
K_i=\bigcup_{i=1}^\infty L_i$
be two representations. Since each $K_i$ is generated by a finite
number of elements of ${\mathcal{K}}$, it is contained in some
$L_j$. Using the inequalities of Lemma 2.3 we see that the limit
of the ratio in question for the tower $\{L_j\}$ is less than or equal
to the corresponding value for the tower $\{K_i\}$. As in the proof of
Lemma 2.4, we have to treat $\phi_p$ first, and then to use
induction over the degree. We get 
$\phi_\alpha ({\mathcal{L}})\le \phi_\alpha ({\mathcal{K}})$
and, by symmetry, vice versa.
${\Box}$

{\bf Lemma 2.6.} {\em Let ${\mathcal{K}}\subset {\mathcal{L}}$ be infinite
global fields. Then
$$ \phi_{\mathbb R}({\mathcal{K}})+2\phi_{\mathbb C}({\mathcal{K}})\ge
\phi_{\mathbb R}({\mathcal{L}})+2\phi_{\mathbb C}({\mathcal{L}})$$ and for
any prime $p$
and any $n\ge 1$
$$\sum_{m=1}^{n} m\phi_{p^m}({\mathcal{K}})\ge \sum_{m=1}^{n}
m\phi_{p^m}({\mathcal{L}}).
$$ In particular, $$\phi_{p}({\mathcal{K}})\ge \phi_{p}({\mathcal{L}}).$$}

{\em Proof.} We follow the same lines as above. Consider a tower $\{L_j\}$
filtering
${\mathcal{L}}$. Then $\{K_j={\mathcal{K}}\cap L_j\}$ filters
${\mathcal{K}}$. Using the monotonicity of Lemma 2.3 for the pair
$K_j\subset L_j$
we get the result.
${\Box}$

{\bf  Definition 2.2.} An infinite global field ${\mathcal{K}}$ (or an
asymptotically exact family
${\mathcal K} =\{K_i\}$)
  is called {\em asymptotically
bad} if and only if
all $\phi_\alpha=0$. If at least one of the parameters $\phi_\alpha > 0$, it
is called {\em
asymptotically good}.

For example, any sequence of global fields of bounded degree is
asymptotically bad.

{\bf Example 2.1.} Families $\{K_i\}$ satisfying the condition $$
{n/\log\sqrt{\vert D\vert}}\rightarrow 0$$of the
Brauer--Siegel theorem are asymptotically bad. Indeed, this condition
clearly implies
$\phi_{\mathbb{R}}=\phi_{\mathbb{C}}=0$
and then all $\phi_q=0$ as well, since one has

{\bf Lemma 2.7.} {\em For  any asymptotically exact family  ${\mathcal{K}}=\{K_i\}$ 
of number fields and for
any prime $p$
one has
$$\sum_{m=1}^{\infty} m\phi_{p^m} \leq
\phi_{\mathbb{R}}+2\phi_{\mathbb{C}}. $$}

{\em Proof.} Indeed, 
 all places of a number field $K$ whose norm is a power of $p$ lie over $p$
in
 $\mathbb Q$ and the product of their norms is not greater than
 $p^n$, where  $n=r_1+2r_2$ is the degree of $K$.  $ {\Box}$

{\bf Lemma 2.8.} {\em Let $L$ be a (finite) global field and $L^{\ab}$
be its maximal abelian extension.
If for an infinite global field ${\mathcal{K}}$ the field ${\mathcal{K}}\cap
L^{\ab}$
is also infinite, then ${\mathcal{K}}$ is asymptotically bad.
In particular, if ${\mathcal{K}}$ is abelian over a finite global field,
then
${\mathcal{K}}$ is asymptotically bad.}

In the function field case this lemma is proved in \cite{Fr/Pe/St}.
In the number field case this can be proved using   
Artin's dicriminant-conductor formula. Since we do not use this 
result in our paper, we do not prove it here.  

{\bf Example 2.2.} Let now $K_1$ be a number field with the infinite Hilbert
class field tower
 $ \{K_i\}$ --- recall that the existence of such fields was proved in \cite{5} --- or any
other unramified tower. Then the family
$ \{K_i\}$ is asymptotically good, since the ratio $n_i/g_i$ is constant in
such a tower, and thus at least one of $\phi_{\mathbb{R}}$ 
and $\phi_{\mathbb{C}}$ is nonzero.
In fact, in the number field case we do not know essentially different
methods for
 constructing asymptotically good families (however, one can take composita
with a fixed
number field, one can mix different class field towers, etc).

{\bf Remark 2.2.} It might happen that
$\phi_{\mathbb{R}}+\phi_{\mathbb{C}}>0$ but $\phi_q=0$
for any $q$. Moreover, we suspect this to be the case for the example of
Theorem 9.2 below. As yet we are however unable to prove it.

In the function field case there exist three essentially
different techniques to construct
 asymptotically good families: a version of the Golod-Shafarevich method
due to Serre \cite{24} which applies to any finite ground field 
${\mathbb{F}}_r$ (cf. \cite{21}, \cite{18}),
the method of modular curves of different types: 
classical, Drinfeld, Shimura
 (cf. \cite{8}, \cite{29}, \cite{34}) which applies only if $r$ is a square 
(or sometimes another power) but gives much
better results, and explicit construction of the same towers (cf. \cite{4}). For
any $r$ one knows the existence of an asymptotically exact
family ${\mathcal{K}}$ with $\phi_r({\mathcal{K}})>0$ (cf. \cite{24}).  

\section{ Basic inequality}

In this section we prove the basic GRH inequality, as well as its weaker
versions that do not require GRH. We treat the number field case first.

Let $K$ be a number field of degree $n$ and disriminant $D$ with $r_1$ real
and $r_2$ pairs of complex embeddings, and let $\zeta_K (s)$ be its Dedekind
zeta-function; by $P(K)$ we denote the set of prime
divisors of $K$ which can be identified with the set  of
non-archimedean places of $K$.

We use the following form of the Guinand-Weil explicit formula
(see \cite{17}, p.122, eq.2.3).

{ \em Let $F(x)$ be a differentiable even real positive function defined on
the whole real line
${\mathbb{R}}$ such that $F(0)=1$ and there exist  positive real constants
$c$ and $\varepsilon$
such that
$$ F(x), F'(x) \le ce^{-(1/2+\varepsilon)\vert x\vert} {\hbox { {\rm as }}}
 \vert x \vert \rightarrow \infty. $$
Define
$$\Phi (s):= \int _{-\infty}^{\infty}F(x) e^{(s-1/2) x}dx\;.$$
Then we have the following formula:
$$\log{\vert D\vert}= r_1 {\pi \over 2}+ n(\gamma +\log 8\pi)-n\int_0^{\infty}
 {1-F(x) \over 2\sh{x\over 2}} dx-r_1 \int_0^{\infty}{1-F(x) \over
2\ch{x\over 2}} dx$$
$$ -4 \int_0^{\infty}F(x)  \ch {x \over 2} dx+{\sum_{\rho}}' \Phi (\rho)
 +2\sum_{P}\sum_{m=1}^{\infty} N({P})^{-m/2}
F(m\log N(P))\log N(P),$$
where in the first sum $\rho$ runs over the zeroes of $\zeta_K(s)$ in the
critical strip,
  $\sum'$ meaning that the $\rho$ and $\bar{\rho}$ terms are to be grouped together,
  the external sum in the last term is taken over all prime divisors
${P}\in P(K)$, and in the internal sum $N({P})$
denotes the absolute norm of} $P$.

\subsection {GRH basic inequality}

Let us now apply this formula to the case of an asymptotically exact
family ${\mathcal K}=\{K_i\}$ of number fields.

{\bf GRH Theorem 3.1} (GRH Basic Inequality). {\em  For an asymptotically
exact family of
number fields one has

$$ \sum_q { \phi_q \log q \over \sqrt q-1} + \phi_{\mathbb{R}}
(\log\sqrt{8\pi}+{\pi\over 4}
+ {\gamma \over 2})+\phi_{\mathbb{C}} (\log 8\pi + \gamma) \le 1\;, $$
the sum being taken over all prime powers $q$}.

(In the special case of unramified towers this theorem was
proved by Y.Ihara \cite{9}.)

{\em Proof.} Let us apply the above explicit formula to the field $K=K_i$
from our
sequence, substituting $F(x)=e^{-yx^2}$ for a real positive $y$. We get
$$\log \vert D \vert= r_1 {\pi \over 2}+ n(\gamma +\log 8\pi)-n\int_0^{\infty}
 {1-e^{-yx^2} \over 2\sh{x \over 2}} dx$$
$$-r_1 \int_0^{\infty}{1-e^{-yx^2} \over
2\ch{x \over 2}} dx -4 \int_0^{\infty}e^{-yx^2}  \ch {x \over 2} dx$$
$$+\,{\hbox{Re}}{\sum_{t}}' \int _{-\infty}^{\infty}
 e^{itx-yx^2}dx+2\sum_{  P}\sum_{m=1}^{\infty} N({  P})^{-m/2}e^{-ym^2\log^2
N({  P})}\log N({  P}),$$
where in the first sum $t$ runs over all reals such that $\zeta_K (1/2+it)=0$,
and the exernal sum in the last term is taken over all prime divisors
${  P}$ of $K$.

Dividing the last equation by $2g=\log \vert D \vert$ and using the
relation $n=r_1+2r_2$ one
 gets:

$$1 = {r_1  \over g}({\pi \over 4}+{\gamma \over 2} +{\log 8\pi \over 2})+
{r_2  \over g}(\gamma +\log 8\pi)-{n  \over g}\int_0^{\infty}
 {1-e^{-yx^2} \over 4\sh{x \over 2}} dx$$
$$-
{r_1  \over g}\int_0^{\infty}{1-e^{-yx^2} \over
4\ch{x \over 2}} dx-{2  \over g}\int_0^{\infty}e^{-yx^2}\ch{x \over 2} dx+
{\hbox{Re}}{1\over {2g}}
\sum_t \int _{-\infty}^{\infty} e^{itx-yx^2}dx$$
$$+
\sum_{q}{N(q) \over{g}}\sum_{m=1}^{\infty} q^{-m/2}e^{-ym^2\log^2 q}\log q\;,$$
where $N(q)$ is the number of prime divisors in $K$ of norm $q$.

Thus $1=T_1+T_2-T_3-T_4-T_5+T_6+T_7$ is presented as a sum of seven terms
$T_j$, $j=1, {\ldots} ,7$. Now we set
 $y=1/ {\log g}$ and tend $g$ to infinity. Then
$y$ tends to zero, and we are going to show that $T_1+T_2+T_7$ tends to the left
hand side of the
Basic Inequality, while  $T_3+T_4+ T_5 $ tends to zero, $T_6$ being
non-negative, which proves
the theorem.

Indeed, it is clear that $T_1$ tends to $\phi_{\mathbb{R}} (\log2 \sqrt {2\pi}+
{\pi \over 4}+{\gamma \over 2})$, and $T_2$ tends to $\phi_{\mathbb{C}}
(\log 8\pi + \gamma) .$
Then, for $T_7$ we have:

$$T_7=\sum_{q}{N(q) \over g}\sum_{m=1}^{\infty} q^{-m/2}e^{-ym^2\log^2 q}\log q
\leq \sum_{q}{N(q) \over g}\sum_{m=1}^{\infty} q^{-m/2}\log q $$
since $e^{-ym^2\log^2 q} \le 1$. On the other hand, since
$e^{-ym^2\log^2 q} \ge 1-ym^2\log^2 q$,
we get
$$T_7\ge \sum_{q\le \log g}{N(q) \over g}\sum_{m=1}^{\lbrack \log
^{1/4}g\rbrack}
 q^{-m/2}(1-ym^2\log^2 q)\log q $$
$$\ge (1-{(\log \log g)^2 \over\sqrt{\log g }})\sum_{q\le \log g}{N(q)
\over g}\sum_{m=1}^{\lbrack \log ^{1/4}g\rbrack}q^{-m/2}\log q. $$
These inequalities show that $T_7$ tends to
$$ \sum_{q}\phi_q\sum_{m=1}^{\infty} q^{-m/2}\log q=\sum_q { \phi_q \log q
 \over \sqrt q-1}$$
when $g\rightarrow\infty$. 

Now let us estimate $T_3$, $T_4$, and $T_5$. Since
 $$0\le \int_0^{\infty}{1-e^{-yx^2} \over
4\ch{x \over 2}} dx \le
\int_0^{\infty}{1-e^{-yx^2} \over
4\sh{x \over 2}} dx\;,$$
if we show that $T_3$ tends to zero then $T_4$ also tends to 0.
We write
$$ \int_0^{\infty}{1-e^{-yx^2} \over \sh{x \over 2}} dx=
 \int_0^{1}{1-e^{-yx^2} \over \sh{x \over 2}} dx+
\int_1^{\delta}{1-e^{-yx^2} \over \sh{x \over 2}} dx+
 \int_{\delta}^{\infty}{1-e^{-yx^2} \over \sh{x \over 2}} dx$$
for any ${\delta}>1$, all three integrals being positive. Since
$e^{-yx^2} \ge 1-yx^2$,
$$\int_0^{1}{1-e^{-yx^2} \over \sh{x \over 2}} dx\le
 \int_0^{1}{yx^2 \over \sh{x \over 2}} dx=O(y). $$
By the same reason
$$\int_1^{\delta}{1-e^{-yx^2} \over \sh{x \over 2}} dx\le (\delta -
1)\max_{1\le x\le\delta}
{1-e^{-yx^2} \over \sh{x \over 2}}\le (\delta - 1)
{1-e^{-y\delta^2} \over \sh{1 \over 2}}\le
{\delta -1\over \sh{1 \over 2}}y{\delta}^2\;. $$ 
For the third integral we have
$$ \int_{\delta}^{\infty}{1-e^{-yx^2} \over \sh{x \over 2}} dx \le
\int_{\delta}^{\infty}{1 \over \sh{x \over 2}}dx=O(e^{-\delta})\;,$$
as $\delta$ tends to infinity.
If $\delta$ tends to infinity  in such a way that $y\delta^3$ tends to zero
(e.g., put
 $\delta=y^{-1/4}=\log^{1\over 4}{g}$), these inequalities show that $T_3$
(and thus $T_4$)
tends to zero.

For $T_5$ we have
$$0\le T_5={2  \over g}\int_0^{\infty}e^{-yx^2}  \ch {x \over 2} dx
\le {2  \over g}\int_0^{\infty}e^{-yx^2+ {x \over 2}} dx
\le {2  \over g}\int_{-\infty}^{\infty}e^{-yx^2+ {x \over 2}} dx$$
$$=
{2  \over g}\sqrt{{\pi  \over y}}e^{1\over 16y}=2 g^{-15/16}\sqrt{\pi\log g},$$
which shows that it also tends to zero.

Then it is sufficient to note that all terms in the sum $T_6$ are positive;
indeed,
$$ \int _{-\infty}^{\infty} e^{itx-yx^2}dx=
e^{-t^2/4y}\int _{-\infty}^{\infty} e^{-y(x-{it\over 2y})^2}dx
=\sqrt{\pi \over y}e^{-t^2/4y}>0\;,$$
which finishes the proof.${\Box}$

Recall that in the function field case we also  have the Basic Inequality
(cf. \cite{9}, \cite{30}, \cite{31}),
which is valid unconditionally:

{\bf Theorem 3.2. }
$$  \sum_{m=1}^{\infty}{ m\beta_m\over {r^{m/2}-1}} \leq 1$$
{\em for any asymptotically exact family of function fields over
${\mathbb{F}}_r
{\;.\;\Box}$}

{\bf Question.} {\em How large the difference $\delta$ between the right
hand side and the left hand side of Theorems $3.1$ and $3.2$
can be?}

We call $\delta$ the {\it deficiency} of the family.
Of course, $0\leq\delta\leq1,$ and for asymptotically bad families
$\delta=1.$ In the function
field case, at least when $r$ is an even prime power,
families of modular curves provide examples
with $\delta=0.$ In the number field case we do not know such a family.
In fact, all known examples are
those of unramified class field towers, sometimes with extra splitting
conditions (see Section 9).
Ihara \cite{9} produced an example of a class field tower of the field
${\mathbb{Q}}(\sqrt{-3\dot 5\dot 7\dot 11\dot 13\dot 17\dot 23\dot 31})$
 with $\delta\leq 0.248...$

Then Yamamura \cite{33} presented other examples of fields having infinite
unramified class towers and very small $\delta$'s. The best of his examples
would have $\delta\leq 0.088...$
The main tool of his paper is a theorem giving a condition for a field
to have an infinite unramified class field tower. This theorem
looks however not to be true, at least, as we are going to show in
Section 9, it contradicts GRH. Unfortunately all the examples
of \cite{33} depend heavily on this theorem, and therefore cannot be
considered
as correct. In Section 9 we discuss the problem in more detail.

Let us remark that the class field tower of the Martinet field
$${\mathbb{Q}}(\cos{2\pi\over11},\sqrt{2},\sqrt{-23})$$
has $\delta\leq 0.1601...$, which was the best
one known for many years. 
Quite recently Hajir and Maire \cite{6} produced several
better examples, the best one they get is given
as follows.  Let $K={\mathbb Q}(\xi)$, where $\xi$ is a root of 
the polynomial
$x^8-9x^6+24182x^4+60281988x^2+895172213$, then $K$ has an infinite class
field tower
ramified (tamely) only at two places over 5. This tower has $\delta\leq 0.141\dots$.

\subsection{Unconditional basic inequalities}

Let us now give some partial results which can be obtained without assuming
GRH. Unfortunately,
they are much weaker.

{\bf Proposition 3.1} (Basic Inequality$'$). {\em For any asymptotically
exact
family of number fields one has
$$2\sum_{q}\phi_q\log q {\sum_{m=1}^{\infty}{1\over{q^m+1}}}+
\phi_{\mathbb{R}}(\gamma/2+1/2 +\log
2\sqrt\pi)+\phi_{\mathbb{C}}(\gamma+\log 4\pi) \leq 1\;,$$
the first sum being taken over all prime powers $q$.}

Note that archimedean coefficients are
$$\alpha_{\mathbb{R}}':=\gamma/2+1/2 +\log 2\sqrt\pi\approx2.054...,
\quad {\hbox {\rm{and}} } \quad\alpha_{\mathbb{C}}':=\gamma+\log
4\pi\approx3.108...,$$
whereas in Theorem 3.1 they are
$$\alpha_{\mathbb{R}}=\gamma/2+\pi/4+\log 2\sqrt{2\pi}\approx2.686...,
\quad {\hbox {\rm{and}} } \quad\alpha_{\mathbb{C}}=\gamma+\log
8\pi\approx3.801...\;.$$

This result ameliorates the unconditional Odlyzko--Serre inequality.

{\em Proof.} The method of the proof is exactly the same as for the Basic
Inequality
(GRH Theorem 3.1)
with the only difference in the choice of the function $F(x)$; without GRH
we choose
$F(x)={{e^{-yx^2}}\over{\ch{x\over 2}}}$ to have $\Re\Phi(s)$ positive on
the whole critical
strip, which is checked by an elementary calculation using the maximum
principle
(cf. \cite{17}, eq. 2.4). Using the Guinand--Weil explicit formula given at the
beginning
of Section 3 and dividing the obtained
equality by $2g$ we get
$$1 = {r_1  \over g}({\pi \over 4}+{\gamma \over 2} +{\log 8\pi \over 2})+
{r_2  \over g}(\gamma +\log 8\pi)-{n  \over g}\int_0^{\infty}
 {1-(e^{-yx^2}/\ch{x \over 2}) \over 4\sh{x \over 2}} dx$$
$$-{r_1  \over g}\int_0^{\infty}{1-(e^{-yx^2}/\ch{x \over 2}) \over
4\ch{x \over 2}} dx-{2  \over g}\int_0^{\infty}e^{-yx^2}
dx+{\hbox{Re}}{1\over {2g}}
\sum_{u+it} \int_{-\infty}^{\infty} {{e^{itx+(u-{1\over
2})x-yx^2}}\over{\ch{x \over 2}}}dx$$
$$+\sum_{q}{N(q) \over{g}}\sum_{m=1}^{\infty}
{2 q^{-m/2}e^{-ym^2\log^2 q}\log q\over  q^{-m/2}+q^{m/2} }\;,$$
where in the first sum $s=u+it$ runs over all zeroes of $\zeta_K(s)$ in the
critical strip,
the rest of
the notation being that of the proof of GRH Theorem 3.1. Thus
$1=T'_1+T'_2-T'_3(y)-T'_4(y)
-T'_5(y)+T'_6(y)+T'_7(y)$
is presented as a sum of seven terms, some of which depend on the value of $y$.
Note that $T'_1+T'_2=T_1+T_2$. We are going to show that  if $y$ tends
to 0
exactly as
described  in the proof of GRH Theorem 3.1 then $T'_3(y)$ tends to
 ${(\phi_{\mathbb{R}}+ 2\phi_{\mathbb{C}})\log 2 \over 2}$, $ T_4(y)$ tends to
${\phi_{\mathbb{R}}({\pi \over 4}-{1 \over 2})}$, $T'_5(y) $ tends to 0,
$T'_6(y)$ is
non-negative, $T'_7(y)$
tending to $2\sum\limits_{q}\phi_q\log q{\sum_{m=1}^{\infty}(q^m+1)^{-1}}$,
which proves the
result. Indeed, the statement
on $T_6$ is obvious from the very choice of $F(x)$, as explained above.
The statement on $T'_5(y)$ follows from the  bound $\vert T'_5(y) \vert \le
T_5$ implied by
the inequality $ ch{x \over 2}\ge 1$. Note that here, as above, we take
$y={{1}\over{\log g}}$.
To prove the statements on
$T_i'(y)$ for $i=3,4$ and 7 it is sufficient to notice that
$T'_3(0)={n\log 2 \over 2g}$, $ T'_4(0)={r_1 \over g}({\pi \over 4}-{1
\over 2})$,
which is an elementary calculation of integrals, $T'_7(0)=
2\sum\limits_{q}\phi_q\log q{\sum_{m=1}^{\infty}(q^m+1)^{-1}}$,
and that $\vert T'_i(y)- T'_i(0)\vert \le T_i$
for $i=3,4, 7$,  which follows  from the same inequality $\ch{x \over 2}\ge
1$ as well.
${\Box}$

We shall also present the following weaker result,
proved by Y.Ihara \cite{9} for the case of unramified towers. It is sometimes
easier to
calculate with (cf. the proof of Theorem 9.7), and has a nice
interpretation in terms of the limit
zeta-function (cf. Remark 4.3).

{\bf Proposition 3.2} (Basic Inequality$''$). {\em For any asymptotically
exact family of
number fields one has

$$\sum_{q}{\phi_q\log q \over q-1}+(\gamma/2+\log 2\sqrt\pi)\phi_{\mathbb{R}}+
(\gamma+\log 2\pi)
\phi_{\mathbb{C}} \leq 1\;.$$
}
 
Note that archimedean coefficients are
$$\alpha_{\mathbb{R}}'':=\gamma/2+\log 2\sqrt\pi\approx1.554...,
\quad {\hbox {\rm{and}} } \quad\alpha_{\mathbb{C}}'':=\gamma+\log
2\pi\approx2.415...\;.$$
 
{\em Proof.} To prove the result one uses Stark's formula (cf. \cite{17}, p. 120)
$$\log{\vert
D\vert}=r_1(\log\pi-\psi(s/2))+2r_2(\log(2\pi)-\psi(s))-{2\over
s}-{2\over{s-1}}$$
$$+2{\sum_{\rho}}'{1 \over s-\rho}+2\sum_{P}\sum_{m=1}^{\infty}
N(P)^{-ms}\log N(P)\;,$$
where in the first sum $\rho$ runs over the zeroes of   $\zeta_K (s)$ in the
critical strip,
 and $\sum'$ means that the $\rho$ and $\bar \rho$ are to be grouped
together, while
 the external sum in the double sum is taken over all prime divisors
${P}$ of $K$, and $\psi(s)=\Gamma'(s)/\Gamma(s)$. 
We then  apply this formula to $K=K_j$ for $s=1+{1\over \sqrt{g_j}} $,
where as usual $g_j= \log \sqrt{\vert D_j\vert} $, and divide it by $2g_j$.
We get
$$1={r_1\over 2g_j}(\log\pi-\psi(s/2))+{r_2\over g_j}(\log 2\pi-\psi(s))-
{1\over \sqrt{g_j}(1+\sqrt{g_j}) }-{1\over{\sqrt{g_j}}}$$
$$+{{1}\over{g_j}}{\sum_{\rho}}'{1 \over s-\rho}+{{1}\over{g_j}}\sum_{P}
\sum_{m=1}^{\infty} N(P)^{-ms}\log N(P).$$

When $g_j$ grows, the first two terms tend, respectively,
to $\alpha_{\mathbb{R}}''\phi_{\mathbb{R}}$ and
$\alpha_{\mathbb{C}}''\phi_{\mathbb{C}}$ since $\psi({1\over
2})=-\gamma-\log 4$ and
$\psi(1)=-\gamma$; the third and the forth terms tend to zero, the fifth
being
positive, and the last term tends to $\sum\limits_{q}{\phi_q\log q \over
q-1}$,
which finishes the proof.
${\Box}$

Proposition 3.2 ameliorates Stark's inequality $$\alpha_{\bf
R}''\phi_{\mathbb{R}} +\alpha_{\mathbb{C}}''\phi_{\mathbb{C}} \leq 1.$$

\section{Zeta function}

We define the {\em limit zeta function } of an asymptotically exact family 
${\mathcal K}$ as
$$\zeta_{\mathcal K}(s)=\zeta_{\phi}(s):=\prod_{q}(1-q^{-s})^{-\phi_q},$$
and its {\em completed limit zeta function } as
$$\tilde\zeta_{\mathcal K}(s)=\tilde\zeta_{\phi}(s):= e^s
2^{-\phi_{\mathbb{R}}}\pi^{{-s\phi_{\mathbb{R}}}/2}
(2\pi)^{-s\phi_{\mathbb{C}}}\Gamma({s\over 2})^{\phi_{\mathbb{R}}}
\Gamma(s)^{\phi_{\mathbb{C}}}
\prod_{q}(1-q^{-s})^{-\phi_q}$$
in the number field case, and
$$\tilde\zeta_{\mathcal K}(s)=\tilde\zeta_{\phi}(s):= r^s 
\prod_{q}(1-q^{-s})^{-\phi_q}$$
in the function field case;
$q$ runs over all prime powers in the number field case 
and over powers of $r$ in the function field case.
As usual, by raising to a complex power a function in $s$ 
defined for $\Re s > a \ge 0$ and such that its values are real positive for real $s>a$, 
we mean unique analytic continuation of what is real positive for real $s$. For an expression
$(1-q^{-s})^{-\phi_q}$ this is the same as to take the value given by the binomial
series.
 
Note that our definition of $\tilde\zeta_{\phi}$ slightly differs in the function
field case from that of \cite{31}.
Note also that $\zeta_{\phi}(s)$ and $\tilde\zeta_{\phi}(s)$ depend only on $\phi=\{\phi_\alpha\}$ and do not
depend on the particular sequence of global fields.

{\bf GRH Proposition 4.1.} 
{\em   For any asymptotically exact family of global fields the product defining
$\zeta_{\phi}$ (and $\tilde\zeta_{\phi}$) converges absolutely
 in the closed half-plane   $\Re(s)\geq{1\over 2}$, and defines 
an analytic function in
$\Re(s)> {1\over 2}$.
In the function field case the result is unconditional.}

{\em Proof.} The product converges absolutely  if and only if the series
$$\sum\limits_{q}\phi_q\log {\vert {1 \over 1-q^{-s}}\vert }$$ does,
but for ${\hbox{Re}}(s)\geq 1/2$ one has
$$\sum\limits_{q}\phi_q\log {\vert {1 \over 1-q^{-s}}\vert }\le \sum_{q}\phi_q\log {1 \over
1-q^{-1/2}}\;,$$ 
which in its turn converges, since,  starting from some $q$, 
$$ \phi_q\log {1 \over
1-q^{-1/2}}\le {\phi_q \log q \over \sqrt q-1},$$
and the series 
$$\sum_q { \phi_q \log q \over \sqrt q-1}\;$$
converges by GRH Theorem 3.1.
$\Box$

{\bf Remark 4.1.} In fact since the coefficients of the Dirichlet series corresponding to
$\zeta_{\phi}$ are positive and $\zeta_{\phi}(1/2)$ is finite,  the product converges absolutely  
for ${\Re}(s)>{1\over 2}-\varepsilon(\phi)$, where $\varepsilon(\phi)$ depends on how
large $\phi_q$ are. It can even happen that $\tilde\zeta_{\phi}(s)$ is an entire function. For
asymptotically bad families we have $\tilde\zeta_{\phi}(s)=e^s$ or $r^s$ depending on the case under
consideration. In fact, we do not know a single example of an infinite global field for which the
product does not converge in the half-plane ${\Re}(s)>0$; note also that the archimedean factors are
analytic in this half-plane.

Let now $$\xi_{\phi} (s):= (\log \tilde\zeta_{\phi})'=\tilde\zeta_{\phi}'/\tilde\zeta_{\phi}=$$
$$1-{\phi_{\mathbb{R}}\over 2}\log\pi-\phi_{\mathbb{C}}\log{2\pi}+
{1\over 2}\phi_{\mathbb{R}}\psi({s\over 2})+
\phi_{\mathbb{C}}\psi(s)-\sum_q \phi_q{\log q\over q^s-1}\;$$ in the number field case 
(where, as before, $\psi(s)={\Gamma'\over\Gamma}(s)\;),$ and 
$$\xi_{\phi} (s):= (\log_r \tilde\zeta_{\phi})'={1\over \log r}\left ({\tilde\zeta_{\phi}'\over
\tilde\zeta_{\phi}}\right )=1-\sum_{m=1}^\infty {m\phi_{r^m}\over r^{ms}-1}\;$$
in the function field case.

Then one can express the   Basic Inequality (GRH Theorem 3.1 and Theorem 3.2) as
$$\xi_{\phi} (1/2)\ge 0\;,$$
and the Generalized Brauer--Siegel Theorem (GRH Theorem 7.2 and Theorem 7.3 below) either as
$$\lim_{i\rightarrow \infty}{ \log{h_i R_i} \over g_i}=\log\tilde\zeta_{\phi}(1),  $$
or as
$$\lim_{i\rightarrow \infty}{ \log\kappa_i \over
g_i}=\log\zeta_{\phi}(1),  $$
where  $\kappa_i$ is the residue of $\zeta_{K_i}(s)$ at 1.  
The function field case of the Generalized Brauer--Siegel Theorem reads either as
$$\lim_{i\rightarrow \infty}{ \log_r h_i \over g_i}=\log_r(\tilde\zeta_{\phi}(1))$$
or as
$$\lim_{i\rightarrow \infty}{ \log_r \kappa_i \over
g_i}=\log_r(\zeta_{\phi}(1))\;.$$

Unconditionally, we have  

{\bf   Proposition 4.2.}
{\em a$)$ For any asymptotically exact family of global fields the product defining
$\zeta_{\phi}$ (and
$\tilde\zeta_{\phi}$) converges absolutely
 in the closed half-plane   $\Re s\geq 1\ \;$, and defines an analytic function in
$\Re s> 1 \;$.

b$)$ For $\Re s>1$ we have the pointwise limits

$$\zeta_{\phi}^0(s)=\lim_{j\rightarrow \infty } \zeta_{K_j}(s)^{1/g_j},$$
$$\tilde\zeta_{\phi}(s)=\lim_{j\rightarrow \infty } \tilde\zeta_{K_j}(s)^{1/g_j},$$
where 
$ \tilde\zeta_{K_j}(s)$ is the completed zeta-function defined by
$$ \tilde\zeta_{K_j}(s)=\mid D_j\mid^{s/2} 2^{-r_1(K_j)}\pi^{-s r_1(K_j)/2}
(2\pi)^{-s r_2(K_j)}\Gamma(s/2)^{r_1(K_j)}
\Gamma(s)^{r_2(K_j)}
\zeta_{K_j}(s)$$
in the number field case, and by
$$ \tilde\zeta_{K_j}(s)=r^{s(g-1)}\zeta_{K_j} (s)$$
in the function field case. For any $\varepsilon >0$ the 
convergence is uniform in the half-plane
$\Re s>1+\varepsilon $
 $($and, thus on compact subsets in the half-plane
$\Re s>1 )$. 

c$)$ Let $s_0>1$, $s_0\neq 2,4$, be a real number such that  the limit 
$$ \lim_{j\rightarrow \infty } \tilde\zeta_{K_j}(s_0)^{1/g_j}   $$
exists for some family $\{K_j\}.$ Then  the family is asymptotically exact. }

{\em Proof.} The proof of a) is the same as that of GRH Proposition 4.1, but in
place of
GRH Theorem 3.1 we use Proposition 3.2.

Let us prove b). The proof is essentially the same as in the 
function fields case considered in \cite{31}. 
Let $K_j$ be a field from our family. Note that it is sufficient to consider
the case of  real $s$, since $\mid \zeta_{K_j}(x)^{1/g}/\zeta_\phi(x)-1\mid \le
\mid \zeta_{K_j}(s)^{1/g}/\zeta_\phi(s)-1\mid$ if $x=  s+it$ with $s>1$ (look at
the corresponding Dirichlet series). For a real
$s>1$ we have

$$\zeta_{K_j}(s)^{1/g}/\zeta_\phi(s)=\prod_{q} 
(1-q^{-s})^{N_q(K_j)/g_j-\phi_q}.
$$

Let $g_0$ be a positive integer such that for any $q\le M$ (here $M$ is
a positive integer to be specified below) one has
$|\phi_q-N_q(K_j)/g_j|\le\delta_1$ for $g_j\ge g_0$, where a real positive
$\delta_1$ is also to be specified.  Then we have

$$\prod_{q\le M}(1-q^{-s})^{\delta_1}\le
\prod_{q\le M}(1-q^{-s})^{N_q(K_j)/g_j-\phi_q}\le
\prod_{q\le M}(1-q^{-s})^{-\delta_1}.
$$

For any $s\ge1+\epsilon$ the product $\prod_{q\le M}(1-q^{-s}) $ satisfies
the inequalities

$$F_M( \epsilon)^{-1}\le\prod_{q\le M}(1-q^{-s})\le F_M( \epsilon)$$
for some real $F_M( \epsilon)>1$. Thus

$$F_M( \epsilon)^{-\delta_1}\le\prod_{q\le M}(1-q^{-s})^{N_q(K_j)/g_j-\phi_q}\le
F_M( \epsilon)^{\delta_1}.
$$

Let us now estimate the ``tail''
$\prod_{q\ge M+1}(1-q^{-s})^{N_q(K_j)/g_j-\phi_q}$. 

We show first that always  $\vert N_q(K_j)/g_j-\phi_q\vert \le 3q$  for sufficiently  large $q$.
Indeed,  in the function field case we have 
$$N_q(K_j)/g_j \le {q+1+2g_jq^{1/2}\over  g_j}\le 2q 
$$
and $\phi_q\le {(q^{1/2}-1)/ m}<q^{1/2}$ by the Basic Inequality, 
which proves the assertion in
this case. In the number field case $N_q(K_j)\le  n_j=
\deg K_j \le \const\cdot g_j$ and
thus
$N_q(K_j)/g_j$ is bounded by a constant; on the other 
hand the unconditional Basic Inequality
implies that
$\phi_q<q$,   which proves the assertion (with  much room to spare) 
in the number field case.

The assertion implies that
the tail lies between
$$G( \epsilon,M)=\prod_{q\ge M+1}^\infty\big(1-q^{-(1+\epsilon)}\big)^{3q }
$$
and its inverse.  

Since $(1-q^{-(1+\epsilon)})^{3q}=
\big((1-q^{-(1+\epsilon)})^{3q^{(1+\epsilon)}}\big)^{q^{-\epsilon}}$
and $(1-q^{-(1+\epsilon)})^{3q^{(1+\epsilon)}}$ tends to $e^{-3}$ for $q\to\infty$
we get the tail to be between
$C\exp\Big(\sum\limits_{q\ge M+1}^\infty q^{-\epsilon}\Big)$ and its
inverse for any $C>1,\ s\ge 1+\epsilon$, and $M$ sufficiently large.
Choosing $C$ and $\delta_1$ such that 
$F_M( \epsilon )^{\delta_1}\le\sqrt{1+\delta}$ and $M$ such that
$C\exp\Big(\sum\limits_{q\ge M+1}^\infty q^{-\epsilon}\Big)\le\sqrt{1+\delta}$, we
get the result.

Let us prove c). Let us suppose that the family is not 
asymptotically exact; it means that there
exists $\alpha$ such that the sequence $N_{\alpha}(K_j)/g_j $ 
has at least two different limit
points, which we denote $ \phi'_{\alpha}$ and $ \phi''_{\alpha}$. 
We can choose two
asymptotically exact subfamilies $ \mathcal K'$ and 
$ \mathcal K''$ of our family such that 
$ \phi_{\beta }(\mathcal K')= \phi_{\beta}(\mathcal K'')$ for all 
$\beta\neq \alpha$ and
$ \phi_{\alpha }(\mathcal K')= \phi'_{\alpha}$, 
$ \phi_{\alpha }(\mathcal K'')= \phi''_{\alpha}$.
Then using b) we see that our condition implies $(f_\alpha(s_0))^{\phi'_{\alpha}}
=(f_\alpha(s_0))^{\phi''_{\alpha}}  $ where $f_\alpha$ is the factor in the product defining
$\tilde\zeta_{\phi}$, corresponding to $\alpha$. 
Since $f_\alpha(s_0)\neq 1$ for $s_0\neq 2,4 $ (for
non-archemedean factors this is true for any $s_0$, but for the gamma factors one needs to
throw out the values $s_0= 2,4$) we get the result.
${\Box}$

{\bf Remark 4.2.} Let us sketch briefly   another way to prove Proposition 4.2, b. Each
$\zeta_{K_j}(s)^{1/g_j}$ can be expressed as an ordinary Dirichlet series
$\sum^{\infty}_{m=1}a_j(m)m^{-s}$ convergent for $ \Re (s)>1$ if the terms in the Euler product, raised to
the $1/g_j$ power are expanded using the binomial theorem. For each fixed $m\ge 1$ the sequence
$a_j(m)$, $j=1,2,\ldots$, goes to a limit as $j\rightarrow \infty$. 

{\bf Remark 4.3.} Proposition 3.2 can be rewritten as $\xi_{\phi}(1) \ge 0$.

{\bf Remark 4.4.} Part b) of Proposition 4.2 shows in which sense the sequence of
zeta-functions of an
asymptotically exact family of global fields tends to the limit zeta-function.

{\bf Remark 4.5.} The definition of $\tilde\zeta_{\phi}(s)$ in the number field case is
chosen on the one hand so as to write the
Basic Inequality and the Generalized Brauer--Siegel Theorem in the shortest
possible way, and
on the other hand so that it is the natural analogue of the function
$$\tilde{\zeta}_K(s)=\vert D\vert^{s/2} 2^{-r_1}\pi^{-s r_1/2}
(2\pi)^{-s r_2}\Gamma(s/2)^{r_1}
\Gamma(s)^{r_2}
\zeta_K (s)$$
which is invariant under $s\mapsto 1-s$.  
 Note however, that the condition to be invariant under $s\mapsto 1-s$
does not change if the function is multiplied by a constant, and our function
$\tilde{\zeta}_K(s)$ differs from
the function ${\Lambda}_K(s)$ used in \cite{13} by the factor
$2^{-r_1}$. The above formulae strongly suggest our normalization to be the
natural one.

 {\bf Remark 4.6.} Comparing the definitions and results for the number
field case and for the  function field one, we conclude that the ``ground field'' of a
number field ``is'' of cardinality $e\approx 2.718281828459045...$

 \section{Zeta zeroes and the explicit formula}

We are going to study the asymptotic distribution of zeroes of $\zeta_{ K}$
for $g_K$ tending to infinity. We start with number fields and
suppose GRH to hold. Note that in the case of asymptotically bad families this result was
essentially obtained in \cite{12}.
 
\subsection{Number field case}

Let us recall several standard notions and facts from the theory of distributions
 (cf. \cite{23}). Let
${\mathcal S} = {\mathcal S}({\mathbb{R}})$ be the space of  complex valued
infinitely
differentiable functions on ${\mathbb{R}}$ which are rapidly (i.e., faster
than any polynomial)
decreasing together with all their derivatives. This vector space is
naturally equipped with  a standard
topology, so that
the Fourier transform is a topological automorphism of ${\mathcal S}$. Its
topological dual
${\mathcal S}'$
is called the space of tempered distributions. By duality, the Fourier
transform is also defined on
${\mathcal S}'$ and it is also a topological automorphism there.
The space ${\mathcal S}'$ is contained in the space ${\mathcal D}'$
of all distributions, which is the topological dual of the space ${\mathcal
D}$  of complex valued
infinitely
differentiable functions with compact support on ${\mathbb{R}}$. The space
of measures ${\mathcal M}$ is
the topological dual of the space of complex valued continuous functions
with compact support on
${\mathbb{R}}$.
Of course, ${\mathcal M}\subset {\mathcal D}'$. The space of measures
${\mathcal M}$ contains
the cone of positive measures ${\mathcal M}_+$, i.e., of those measures
whose value at a positive
real-valued function is positive. The space of distributions ${\mathcal
D}'$ also contains
the cone of positive distributions ${\mathcal D}'_+$. It is known that
${\mathcal D}'_+ = {\mathcal M}_+$ (cf. \cite{23}, Thm.V of Ch.I). The
intersection
 ${\mathcal M}_{sl}={\mathcal M}\cap {\mathcal S}'$ is called
the space of measures of slow growth. The criterion for a measure to be of
slow growth is that
for some positive integer $k$ the integral
$$I_k= \int_{-\infty}^{\infty}(x^2+1)^{-k} d\mu$$
converges (cf. \cite{23}, Thm.VII of Ch.VII).

Let
${ \mathcal F}=\{ K_j\}$ be  an asymptotically exact family of number fields. For
each $K_j$ we define the measure
$$ \Delta_{K_j}:={\pi\over   g_{ j}}\sum_{\zeta_{ K_j}(\rho)=0}
\delta_{t(\rho)} \;,$$
where $g_{j}:=g_{K_ j},$ $t(\rho)=(\rho-{ 1\over 2})/  i$, and  $\rho$ runs over all
non-trivial  zeroes
 of the zeta-function $\zeta_{ K_j}(s)$; here $\delta_a$ denotes the atomic (Dirac) 
measure at $a$.  Because of GRH
$t(\rho)$ is real, and
$ \Delta_{K_j}$
is a discrete measure on ${\mathbb{R}}$. Moreover, $ \Delta_{K_j}$
is a measure of slow growth, which follows,
e.g., from the Weil Explicit Formula (see  the proof of GRH Theorem 5.1 below).

Now we are ready to formulate the main result of this  section,
expressing the limit
distribution of zeta zeroes in terms of the parameters $\phi=\{\phi_\alpha\}$
of the asymtotically exact family.

{\bf GRH Theorem 5.1} (GRH Explicit Formula). {\em For an asymptotically exact family $\mathcal K$, 
in the space of measures of
slow growth on
${\mathbb{R}}$
 the limit
$$ \Delta=\Delta_{\mathcal K}:=\lim _{j\rightarrow \infty}\Delta_{K_j} $$
exists.
 Moreover, the measure $\Delta $ has a continuous density $M_{\phi}\;$,
$$M_{\phi}(t)=\Re \left(\xi_{\phi}\left({1 \over 2}+i t\right)\right)=$$
$$1-\sum_{q}\phi_q h_q(t)\log q + {1\over 2}\phi_{\mathbb R}
 \Re\psi\left({1\over 4}+{it\over 2}\right)+\phi_{\mathbb C}
\Re\psi\left({1\over 2}+it\right)-{\phi_{\mathbb R}\over 2}
\log\pi- \phi_{\mathbb C}\log{2\pi}$$
where
$$h_q(t)={\sqrt q \cos( t\log q) -1\over q+1-2\sqrt q
\cos( t\log q)},\quad\psi(s)=
{\Gamma'\over\Gamma}(s)\;.$$}

{\bf Remark 5.6.} This density depends only on the parameters
$\phi=\{\phi_\alpha\}$.

{\bf GRH Corollary 5.1} (GRH Basic Equality). {\em We have
$$\sum_q { \phi_q \log q \over \sqrt q-1} + \phi_{\mathbb{R}}
(\log\sqrt{8\pi}+{\pi\over 4}
+ {\gamma \over 2})+\phi_{\mathbb{C}} (\log 8\pi + \gamma)=
1-M_{\phi}(0)\;,$$
 in other words, the deficiency
$$\delta_{\mathcal K}=\xi_\phi({1\over 2})=M_{\phi}(0)\; .$$}

{\em Proof of GRH Corollary 5.1.} Put $t=0$ in the formula for
$M_{\phi}(t)$.$\Box$

{\em Proof of GRH Theorem 5.1.} Here is the strategy of the proof: First we
are going to prove
that $\Delta$ exists as a limit in ${\mathcal S}'$, therefore, being positive,
it lies in ${\mathcal M}_{sl,+}={\mathcal S}'_+$. The next point is to show
that
this measure is absolutely continuous, i.e., of the form $F(t)dt$; to do it
we have
to prove that neither skyscraper, nor singular component occurs. Then we shall
compare this measure $\Delta=F(t)dt$ with the measure
$\Delta_0=M_{\phi}(t)dt$ of the theorem:
we first show that they coincide on the set of some specific functions
$H_{y,a}(t)$, and then that
this is enough for the measures to be equal.

We begin by proving that $ \Delta$ is well-defined as a
tempered distribution.
We are going to use  the  Weil
Explicit Formula in the form presented in \cite{19}, Section 1. We use the notation 
of Section 3; here we suppose that $F \in {\mathcal S}({\mathbb{R}})$, and that it
satisfies the condition
$$ F(x), F'(x) \le ce^{-(1/2+\varepsilon)\vert x\vert}  {\hbox { as }} 
\vert x \vert \rightarrow \infty . \quad\quad\quad\quad\quad\quad\quad\quad (*)$$

 Note that  for $s={1\over 2}+it$ 
$$\Phi (s):= \int _{-\infty}^{\infty}F(x) e^{(s-1/2) x}dx=\hat F (t)\; ,$$
where
$$ \hat F (t)= \int _{-\infty}^{\infty} F(x) e^{itx}dx \in {\mathcal
S}({\mathbb{R}})$$
is the Fourier transform of $F.$

The abovementioned Weil Explicit Formula reads as follows:

{ \em Let $K$ be a number field, then the limit
$$ {\sum}' \Phi(\rho):=\lim_{T\rightarrow \infty}
 {\sum_{\vert \rho \vert < T}} \Phi (\rho) =\lim_{T\rightarrow \infty}
 {\sum_{{{\vert {1\over 2}+it \vert < T}\atop {\zeta_K({1\over 2}+it)=0
}}}} \hat F (t) $$
exists, where in the  sum $\rho$ runs over the set of zeroes of
$\zeta_K(s)$ on  the critical line $\Re(s)=1/2$ $($which is
supposed to be the set of all critical zeroes of $\zeta_K(s)\;)$, and we have
 the following formula$:$

$${\sum}'\Phi (\rho)-\Phi(0)-\Phi(1)=
F(0)(\log{\vert D_K \vert}-r_1 \log\pi- 2r_2\log{2\pi})$$
 $$-\sum_{P}\sum_{m=1}^{\infty} N({P})^{-m/2}[F(m\log N(  P))
+F(-m\log N(  P))]\log N(  P) +$$
$${r_1\over 2\pi} \int_{-\infty}^{\infty}{\hat F(t) +\hat F(-t)\over
2} \Re\psi\left({1\over 4}+{it\over 2}\right)dt+{r_2\over \pi}
\int_{-\infty}^{\infty}{\hat F(t)
+\hat F(-t)\over 2}\Re\psi\left({1\over 2}+it\right)dt\;,$$
where  the external sum is taken over all primes
$P$ of $K$, and $N({P})$ denotes the absolute
norm of $P$.}

 This is exactly the formula of \cite{19}, Section 1. We apply the formula only
to functions from
$ {\mathcal S}({\mathbb{R}})$ which clearly satisfy the other conditions of
the theorem of
\cite{19}, Section 1.
Though there
the function $F$ is assumed to be even, we can apply the result to any
function $F$ from
$ {\mathcal S}({\mathbb{R}})$, replacing it by $ { F(t)+F(-t)\over 2}$.

One can rewrite the sum
$$\sum_{P}\sum_{m=1}^{\infty} N({P})^{-m/2}
[F(m\log N(  P))+F(-m\log N(  P))]\log N(  P) $$
as
$$\sum_{q} N_q(K)\sum_{m=1}^{\infty} q^{-m/2} [F(m\log q)+F(-m\log q)]\log q
,$$
the sum being taken over all prime powers $q$.

 Let $\hat {\mathcal D} \subset   {\mathcal S} $ be the  Fourier dual  of
${\mathcal D}$;
 it is a dense subspace of ${\mathcal S} $ (since ${\mathcal D} $ is dense
in ${\mathcal S}$).
Let $K=K_j$ be a field from our
family. Take any $f\in \hat {\mathcal D}$ and let $f=\hat F$, $F \in
{\mathcal D}$. We have $f(t)=\Phi({1\over 2}+it)$. The function $F$
 satisfies the above condition (*), and,
dividing by $2g_{ K_j}=\log{\vert D_{ K_j}\vert}$, we get:
$$\Delta_{K_j}(f)=
{\pi\over g_{ K_j}} \mathop{{\sum}'}_{\zeta_K({1\over 2}+it)=0} f(t) =$$
$$2\pi \left ({\Phi(0)+\Phi(1)\over 2g_{ K_j}}+
F(0)\left (1-{r_1\over 2g_{ K_j}}\log\pi- {r_2\over g_{
K_j}}\log{2\pi}\right)\right)$$

$$-2\pi \left (\sum_{q} {N_q(K_j) \over 2g_{ K_j}}\sum_{m=1}^{\infty} q^{-m/2} [F(m\log
q)+F(-m\log q)]\log q\right) +$$

$${r_1\over 2g_{ K_j}} \int_{-\infty}^{\infty}{f(t) +f(-t) \over 2}
 \Re\psi\left({1\over 4}+{it\over 2}\right)dt+{r_2 \over g_{ K_j}}
\int_{-\infty}^{\infty}
{f(t) +f(-t)\over 2}\Re\psi\left({1\over 2}+it\right)dt\;.$$

If we fix $f$ and tend $j$ to infinity then the right hand side tends to
$$ \Delta(f):=2\pi F(0)\left (1-{\phi_{\mathbb R}\over 2}\log\pi- \phi_{\mathbb
C}\log{2\pi}\right)$$
$$-\pi\sum_{q} {\phi_q}\sum_{m=1}^{\infty} q^{-m/2}
 [F(m\log q)+F(-m\log q)]\log q + $$
$${\phi_{\mathbb R}\over 2} \int_{-\infty}^{\infty}{f(t) +f(-t) \over 2}
 \Re\psi\left({1\over 4}+{it\over 2}\right)dt+\phi_{\mathbb C}
\int_{-\infty}^{\infty}
{f(t) +f(-t)\over 2}\Re\psi\left({1\over 2}+it\right)dt\;,$$
since the family is asymptotically exact. This shows the existence of the
 limit
 $$\Delta(f)=\lim_{j\rightarrow\infty}\Delta_{ K_j}(f)\;$$
for any $f\in \hat {\mathcal D}$. The map $f\mapsto \Delta(f)$ is obviously
 linear. Thus in order to prove that $\Delta $ is a
tempered  distribution it is sufficient to verify that $\Delta $ is
 continuous on $\hat {\mathcal D}$ in the topology of ${\mathcal S}$, since
$\hat {\mathcal D}$
 is dense in ${\mathcal S}$.
 Let now $\{f_i\}$, $i=1,2,...$, be
a sequence with  $f_i\in \hat {\mathcal D}$ tending to zero in the topology
of ${\mathcal S}$. Since the Fourier transform is a
topological automorphism of ${\mathcal S}$, we conlude that
the sequence $\{F_i\}$, where $f_i=\hat F_i$, tends to zero as well.
In particular, both sequences $\{f_i\}$ and $\{F_i\}$ tend
to zero uniformly. Let us then  show   that $\Delta$ is continuous. Indeed,

$$ \Delta(f_i)=T_1(f_i)-T_2(f_i)+T_3(f_i)+T_4(f_i)\;,$$
where
$$T_1(f_i)=2\pi F_i(0)\left( 1-{\phi_{\mathbb R}\over 2}\log\pi- { \phi_{\mathbb
C}}\log{2\pi}\right)=
\delta(\hat f_i)\left( 1-{\phi_{\mathbb R}\over 2}\log\pi- { \phi_{\mathbb
C}}\log{2\pi}\right)\;,$$
$$T_2(f_i)=\pi\sum_{q} {\phi_q}\sum_{m=1}^{\infty} q^{-m/2}
[F_i(m\log q)+
F_i(-m\log q)]\log q\;, $$

$$T_3(f_i)={\phi_{\mathbb R}\over 2} \int_{-\infty}^{\infty}{f_i(t)
+f_i(-t) \over 2}
 \Re\psi\left({1\over 4}+{it\over 2}\right)dt\;,$$
and
$$T_4(f_i)=\phi_{\mathbb C} \int_{-\infty}^{\infty}{f_i(t) +f_i(-t) \over 2}
 \Re\psi\left({1\over 2}+it\right)dt\;.$$

The term $T_1(f_i)$ is clearly continuous as well as the terms
$T_3(f_i)$ and $T_4(f_i)$, since the measures $\Re\psi\left({1\over
2}+it\right)dt$
and $\Re\psi\left({1\over 4}+{it\over 2}\right)dt$ are of slow growth (cf.
\cite{13}, XVIII.1). 
To prove the continuity of $T_2(f_i)$ one notes that from the GRH version
of the Basic
Inequality (GRH Theorem 3.1 above) it follows that
$$T_2(f_i)\le2\pi \sup_{x\in \mathbb R}\vert F_i(x)\vert  $$
which finishes the proof of the continuity of $\Delta$.

Thus we see that $\Delta$ is a tempered distribution, $\Delta\in{\mathcal
S}'$. Note 
that by its very definition this distribution is positive, i.e.,
$\Delta(f)$ is real non-negative for a real non-negative $f\in
{\mathcal S}({\mathbb{R}})$. Since any positive distribution is a (positive)
measure, one concludes that $\Delta$ is a measure of slow growth. 
Then we have to show that $\Delta =M_{\phi} dx$.
To do that we need the following lemmata.

For any $a\in {\mathbb{R}}$ and any $y>0$ let us
define the function $H_{y,a}(x) \in {\mathcal S}({\mathbb{R}}) $

$$ H_{y,a}(x):={1\over2 \sqrt{\pi  y}}\exp ({-(a-x)^2 \over 4y})$$
so that 
$$F_{y,a}(x)={\hat H_{y,a}(x)\over 2\pi}:={1\over 2\pi} \exp ({-yx^2+iax}).
$$

{\bf Lemma 5.1.}  {\em We have}
$$ \lim_{y\rightarrow +0}\Delta(H_{y,a})=M_{\phi}(a)\;.$$

Note that for $y$ tending to zero the function $H_{y,a}  $ tends
 to $\delta_a$ in the sense of distributions, where $\delta_a$ is the
Dirac measure concentrated in $a$.

{\bf Lemma 5.2.}  {\em Let $\mu$ be a positive mesure on ${\mathbb{R}}$
 such that for any $a\in {\mathbb{R}}$ one has
$$\lim_{y\rightarrow +0}\mu(H_{y,a})=M(a) \eqno (2) $$
with a  function $M$ continuous on ${\mathbb{R}}$. Then
$\mu=M(x)dx.$}

One notes that the theorem follows  from these two lemmata.
 Let us now prove them.

{\em Proof of  Lemma 5.1.} Let us apply  the explicit formula  to $f=H_{y,a}$.
Since  $H_{y,a} $ is even, we get

$$\Delta_{ K_j}(H_{y,a})={\pi(H_{y,a}(0)+H_{y,a}(1))\over g_{K_j}}+
2\pi {F}_{y,a}(0)\left (1-{r_1\over 2g_{ K_j}}\log\pi- {r_2\over g_{
K_j}}\log{2\pi}\right)-
$$
$$\sum_{q} {N_q(K_j) \over 2g_{K_j}}\sum_{m=1}^{\infty} q^{-m/2}
[\hat H_{y,a}(m\log q)+ {\hat H_{y,a}}(-m\log q)]\log q +
$$
$${r_1\over 2g_{K_j}} \int_{-\infty}^{\infty}  H_{y,a}(t)
\Re\psi\left({1\over 4}+{it\over 2}\right)dt +
{r_2 \over g_{K_j}} \int_{-\infty}^{\infty} H_{y,a}(t)
\Re\psi\left({1\over 2}+it\right)dt\; .
$$
Since
$$ \Delta( H_{y,a})=\lim _{j\rightarrow \infty}\Delta_{K_j}( H_{y,a}),
$$
$$ 2\pi {F}_{y,a}(0)=1,$$
we get
$$\Delta(H_{y,a})= 1-{\phi_{\mathbb R}\over 2}
\log\pi- \phi_{\mathbb C}\log{2\pi}-{1\over 2}\sum_{q} \phi_q\log q
\sum_{m=1}^{\infty}  q^{-m/2}
e^{-y(m\log q)^2}(q^{iam}+q^{-iam}) +
$$
$${\phi_{\mathbb R}\over 2}  \int_{-\infty}^{\infty} H_{y,a}(t)
\Re\psi\left({1\over 4}+{it\over 2}\right)dt+
\phi_{\mathbb C}\int_{-\infty}^{\infty}
H_{y,a}(t) \Re\psi\left({1\over 2}+it\right)dt\;.
$$
Let us tend $y$ to zero. Then the expression
$e^{-y(m\log q)^2}(q^{iam}+q^{-iam})$ tends
to $q^{iam}+q^{-iam}=2\cos(am\log q)$ and thus
$${1\over 2}\sum_{q} \phi_q\log q \sum_{m=1}^{\infty}  q^{-m/2}
e^{-y(m\log q)^2}(q^{iam}+q^{-iam})\rightarrow \sum_{q} \phi_q h_q(a)\log q\;,
$$
since
$h_q(a)=\sum_{m=1}^{\infty}  q^{-m/2} \cos(am\log q).$
Since  $H_{y,a}(t)$ tends to $\delta_a$ in the space ${\mathcal M} $,
and  $\Re\psi\left({1\over 4}+{it\over 2}\right)$ is a $C^{\infty}$-function,
 $H_{y,a}(t) \Re\psi\left({1\over 4}+{it\over 2}\right)$ tends to \linebreak
$\delta_a\Re\psi({1\over 4}+i{a\over 2})$,
and thus
 $${\phi_{\mathbb R}\over 2}  \int_{-\infty}^{\infty} H_{y,a}(t)
 \Re\psi\left({1\over 4}+{it\over 2}\right)dt\rightarrow {\phi_{\mathbb
R}\over 2}
 \Re\psi\left({1\over 4}+{ia\over 2}\right)$$
for $y$ tending to zero.
A similar  argument shows that the term
$$\phi_{\mathbb C} \int_{-\infty}^{\infty} H_{y,a}(t) \Re\psi\left({1\over
2}+it\right)dt\rightarrow
\phi_{\mathbb C} \Re\psi\left({1\over 2}+ia\right)$$
for $y\rightarrow 0$, which finishes the proof.$\Box$

{\em Proof of  Lemma 5.2.}  Since $\mu$ is a positive measure one can write
 $\mu=dG$ for a non-decreasing function $G$,
and  the standard decomposition $G=G_0+G_1+G_2$ with absolutely continuous
$G_0$, singular $G_1$ and a jump-function $G_2$ (cf. \cite{10} Ch. VI, section
4) shows that
$\mu=dG=dG_0+dG_1+dG_2$. Let us prove that the  property (1) implies 
$G_2=G_1=0$. Indeed, it is sufficient to show that
$$ \int H_{y,a}(t)  dG_i, \quad i=1,2$$
cannot be bounded for $y\rightarrow +0$.
For $G_2$ it is almost obvious, since
$$dG_2=\sum_n s_n\delta_{t_{n}},$$
where
$$G_2(t)=\sum_{n:t_{n}\le t} s_n,$$
and $\delta_{t_{n}}$ is the Dirac measure at $t_n$. If
$G_2\neq 0$, i.e., if the sum
 is non-empty, let $s_i>0$ and consider
$$ \int H_{y,t_i}(t) dG_2=\sum_n s_n H_{y,t_i} (t_{n})\geq s_i H_{y,t_i}(t_i)
={s_i\over 2\sqrt {\pi y}}, $$
which obviously tends to infinity for $y\rightarrow +0$. Thus, $G_2=0$. If
$G_1\neq 0$ then
 there exist $a,b\in {\mathbb R}$, $a<b$, with $G_1(a)<G_1(b)$. Recall that
since $G_1$ is
 singular, its derivative $G'_1$ is zero almost everywhere. Let us now
show that there exists
$x_0\in [a,b]$ such that
$$\limsup_{\varepsilon \rightarrow +0}{G_1(x_0+\varepsilon)
-G_1(x_0-\varepsilon)
\over \varepsilon } =\infty. \eqno (2) $$
Indeed, if
$${G_1(x_0+\varepsilon) -G_1(x_0-\varepsilon)
\over \varepsilon }\le M $$
for any $x_0\in [a,b]$ and any $\varepsilon>0$ then one can cover the set
 ${\hbox {Supp}}\; G_1\bigcap [a,b]$
by the union of a countable set of intervals  $(a_i, b_i)$ with $\sum
_i(b_i-a_i)<\varepsilon$
(which is possible since the Lebesgue measure of ${\hbox {Supp}}\;
G_1\bigcap [a,b]$ is zero) and thus
 deduce that
$$G_1(b)-G_1(a)\le M\varepsilon $$
for any $\varepsilon>0$, which would imply that $G_1(b)=G_1(a)$. Let us
fix $x_0\in [a,b]$ satifying
the condition (2) and let us consider the value of 
$$\int^{\infty}_{-\infty} H_{\varepsilon^2,x_0}(x) dG_1. $$
Since for $\mid x-x_0\mid\le {\varepsilon}$
$$ H_{\varepsilon^2,x_0}(x)\ge {e^{-1/4}\over 2\sqrt {\pi}\varepsilon}\;, $$
 we get
$$\int^{\infty}_{-\infty}H_{\varepsilon^2,x_0}(x) dG_1\ge
 \int^{x_0+ {\varepsilon}}_{x_0- {\varepsilon}}
H_{\varepsilon^2,x_0}(x) dG_1\ge {e^{-1/4}\over 2\sqrt {\pi}\varepsilon}
 (G_1(x_0+{\varepsilon}) -G_1(x_0-{\varepsilon})),$$
and thus
$$\limsup_{\varepsilon \rightarrow +0}\int^{\infty}_{-\infty}
H_{\varepsilon^2,x_0}(x) dG_1=\infty$$
which gives the desired contradiction and shows that $G_1=0$.

Therefore, the measure  $\mu $ is absolutely continuous, $\mu=D(t) dt$
for a non-negative density function $D(t)$. Then it is sufficient to show
that $D(t)=M(t)$ almost
everywhere. Indeed, it is sufficient to show that
$$ \int_{-\infty}^{\infty} H_{y,a}(t) D (t)dt $$
tends to $D(a)$ for $y$ tending to zero for any point $a$ at which $D(t)$
is continuous. If not,
one supposes that, say,
$$F(a):=\lim_{y\rightarrow 0}\int_{-\infty}^{\infty} H_{y,a}(t) D(t)dt> D(a) $$
(if $F(a)<D(a) $, the argument is the same).
Let us choose $\varepsilon >0, \;\delta >0$ such that $D(t)\le
F(a)-\varepsilon$ for any
$t\in (a-\delta,a+\delta).$ Then we have
$$\int_{-\infty}^{\infty} H_{y,a}(t) D (t)dt =\left(
\int_{a-\delta}^{a+\delta} +
\int_{-\infty}^{a-\delta} +\int_{a+\delta}^{\infty}\right) H_{y,a}(t) D
(t)dt. $$

Since
$$ {1\over 2\sqrt{\pi y}}\exp({-(x-a)^2\over 4y})={1\over 2\sqrt{\pi
y}}\exp({-(x-a)^2\over 8y})\cdot
\exp({-(x-a)^2\over 8y})$$
$$=\sqrt 2 H_{2y,a}(x)\cdot \exp({-(x-a)^2\over 8y}), $$
we see that
$$\mid \int_{-\infty}^{a-\delta} H_{y,a}(t) D
(t)dt+\int_{a+\delta}^{\infty} H_{y,a}(t)
 D (t)dt \mid  \le  $$
$$\mid \int_{-\infty}^{\infty} \sqrt 2 H_{2y,a}(t) D (t)dt\mid 
e^{{-\delta ^2\over 8y }}$$ 
and thus tends to zero for $y$ tending to zero. Therefore,
$$F(a)=\lim_{y\rightarrow 0}\int_{-\infty}^{\infty} H_{y,a}(t) D (t)dt =
\lim_{y\rightarrow 0} \int_{a-\delta}^{a+\delta} H_{y,a}(t) D (t)dt\le $$
$$(F(a)-\varepsilon)\lim_{y\rightarrow 0}\int_{a-\delta}^{a+\delta}
H_{y,a}(t) dt \le
(F(a)-\varepsilon)\lim_{y\rightarrow 0}\int_{-\infty}^{\infty} H_{y,a}(t)
dt=F(a)-\varepsilon,
$$
which gives a contradiction, and finishes the proof both of the lemma and
the theorem.
${\Box}$

{\bf Remark 5.1.} GRH Theorem 5.1 and GRH Corollary 5.1 gives a partial
answer to the following
 question of Odlyzko (Open Problem 6.2 of \cite{17}):

{\em Do the zeroes of $\zeta_K(s)$ in the critical strip
approach the
 real axis as $n\rightarrow \infty$, and if they do, how fast do they so,
 and how many of
 them are there? }

{\bf Remark 5.2.} Theorem 5.1 implies that zeta zeroes are asymptotically
uniformly distributed for any asymptotically bad family, e.g., 
for a family of number
fields of fixed absolute degree.
This is the main result of \cite{12}.

 \subsection{Function field case }
In the function field case, the analogue of Theorem 5.1 is also true. In
\cite{31} we proved the  Asymptotic Explicit Formula which gives
the asymptotic distribution law for Frobenius angles for asymptotically
exact families of function fields, or, which is the same, the limit distribution law for
zeroes of their zeta-functions. Let ${\mathcal K}=\{ K_j\}$ be such a family.
For a  zero  $\rho$ of the zeta-function $\zeta_{ K_j}(s)$ let $t(\rho)$ be
defined by
$$t(\rho):={{\rho-{ 1\over 2}}\over{ i}}.$$

Clearly, $t(\rho)$ is a real number (Weil's theorem)
 defined modulo $2\pi$, and we suppose that $t(\rho) \in (-\pi ,\pi ]$, 
which determines it
uniquely.

Let $$ \Delta_j:={{ \pi \over g_j}}\sum_{\zeta_{
K_j}(\rho)=0}\delta_{t(\rho)},$$
where $\delta_{t(\rho)}$ is, as usual, the Dirac measure supported at
$t(\rho)$.
Then $ \Delta_j$ is a measure of total mass $2\pi $    on
${{\mathbb{R}} / 2\pi {\mathbb{Z}}}$, and  $ \Delta_j$ is symmetric with
respect to $t\mapsto-t$. Points
of ${{\mathbb{R}} / 2\pi {\mathbb{Z}}}$ are given by their representatives in
$(-\pi,\pi ]$.

{\bf Theorem 5.2.} {\em In the function field case, 
for an asymptotically exact family, in
the weak topology on the space of measures on ${{\mathbb{R}} / 2\pi {\mathbb{Z}}}$   
the limit
$$ \Delta:=\lim _{j\rightarrow \infty}\Delta_j\; $$
exists.  Moreover, the measure $\Delta $ has a continuous density
$M_{\phi}$, and the following Asymptotic  Explicit  Formula holds:
$$M_{\phi}(t)={\Re} (\xi_{\phi}({1 \over 2}+{i\over\log r}
t))=1-\sum_{m=1}^{\infty}
m\phi_{r^m} h_m(t)$$
for
$$h_m(t)={ r^{m/2}\cos(mt) -1\over r^m+1-2r^{m/2}\cos(mt)}\;,$$
which depends only on the family of numbers $\phi=\{\phi_{r^m}\}$ and we have the
following Basic Equality:
$$\xi_{\phi}({1 \over 2})=1-\sum_{m=1}^{\infty}{m\phi_{r^m}\over
r^{m/2}-1}=M_{\phi}(0)\;.{\Box} $$}

{\bf Remark 5.3.} If one applies Asymptotic  Explicit  Formula to the case
of the maximal
 family ${\mathcal K}$ (i.e., with $\phi_r=\sqrt{r}-1$ and $\phi_{r^m}=0$ for $m\ge
2$)  over
${\mathbb{F}}_{r}$, $r=p^2$, given
 by the reductions of curves $X_0(n)$ to characteristic $p$, then, using
 the Eichler-Shimura relation, one obtains a particular case (modular
 forms of weight two) of Serre's   results on the asymptotic distribution of eigenvalues of
 Hecke operators (cf. \cite{25}, especially Sections 3 and 7), namely:

{\bf Proposition 5.1.} {\em Let $p$ be a fixed prime number, and let $X_n
\subset [-2,2]$
for a positive integer $n$ coprime with $p$ be the set of eigenvalues of
the   operators
 $T_p'(n)= T_p(n)/\sqrt p,$ where $T_p(n)$ is the Hecke operator acting on
the space of cusp
forms of weight 2 and level $n$. Then for $n\rightarrow \infty$ the set
$X_n$ becomes equidistributed
with respect to the measure}
$$ \mu_p={p+1 \over \pi}\cdot{ p\sqrt{1-x^2/4}\, dx \over
(p+1)^2-px^2}.{\Box} $$

\subsection{Lowest zeta-zero}

Theorem 5.1 makes it possible to prove (under GRH) that  the lowest zero of the
zeta-function of  a global field tends to 1/2. Most probably, this result is known to experts, but
we have not found it in the literature (cf. however,  \cite{Omar}).  
 Let us denote this  lowest zero by $\rho_0(K)=1/2+t_0(K)i$. 

{\bf GRH Proposition 5.2.} 
 {\em  For any family $\{K_i\}$ of global fields 
$$ \lim_{g(K_i)\rightarrow \infty} t_0(K_i)=0.$$}

{\em Proof.}   Let us suppose the
contrary: $$ \liminf_{g(K_i)\rightarrow \infty} t_0(K_i)=\varepsilon>0.$$
Passing to a subsequence if necessary, one can suppose that there exists a sequence $
K_1,K_2,\ldots
$
  with $t_0(K_j) \ge {\varepsilon  }$ for any $j$.  Passing to a subsequence once again we
suppose  that our sequence is asymptotically exact  and thus we can apply Theorems 5.1 and 5.2. The
condition $t_0(K_j) \ge {\varepsilon }$ implies that the corresponding limit
density $M_{\phi}$   identically vanishes on the interval 
$(-\varepsilon, \varepsilon)$, which is impossible since $M''_{\phi}(0)>0$ for $\phi\neq 0$, and
$M_{\phi}$  equals  1 identically for $\phi=0$  (note also that $M_{\phi}(0)>0$ if the deficiency
of the family is positive).$\Box$

\section{Further theory}
  In this section we discuss some directions of further study of infinite global fields,
and more generally, of asymptotically exact families. Most of the following problems look
very difficult, but some of  them seem to be easier than the others, and we hope to return to
them elsewhere.

\subsection{Structure of the parameter set}
For an infinite global field $\mathcal K$ we have defined the sequence $\phi_{\mathcal K}$ 
 of parameters:
$$\phi_{\mathcal K}= (\phi_{\mathbb{R}},\phi_{\mathbb{C}},\phi_2,\ldots )$$
for the number field case, and 
$$\phi_{\mathcal K}= (\phi_r,\phi_{r^2},\phi_{r^3},\ldots )$$
for the function field one.
We have shown that these sequences are sometimes nontrivial (i.e., nonzero). Then it is natural
to ask about the structure of possible parameters. We define  the sets $\Phi$ and $\Phi_r$ in 
${\mathbb{R}}^{\infty}$ by
$$\Phi=\{\phi_{\mathcal K}: \; {\mathcal K} {\hbox{ an infinite number field}}\}, $$
$$\Phi_r=\{\phi_{\mathcal K}: \; {\mathcal K} {\hbox{ an infinite function field over } {\mathbb F}_r}\}. $$
We also introduce the sets $\tilde \Phi \subset {\mathbb{R}}^{\infty}$ and $\tilde \Phi_r \subset
{\mathbb{R}}^{\infty}$ defined exactly as $\Phi$ and  $\Phi_r$ , but for all asymptotically exact families. Clearly,
$ \{0\}\in \Phi\subseteq \tilde \Phi$, $ \{0\}\in \Phi_r\subseteq \tilde \Phi_r$.

Above considerations show that neither of these sets reduce to $ \{0\}$. However, their
structure remains mysterious. Let us put some natural questions on this structure, for
brevity, only in the case of $\Phi$; exactly the same questions are equally
interesting for the other three sets.

{\bf Problem 6.1.} {\em  Is  $\Phi$ closed in some natural topology on
${\mathbb{R}}^{\infty}$? }

One can propose a natural class of topologies to consider. Let
$a=(a_{\mathbb{R}}, a_{\mathbb{C}}, a_{2}, \ldots )$ be a sequence of
positive real numbers indexed  exactly as the sequences $\phi$ are. Then one can define the
weighted spaces $l_{p, a}$ with  $ 1\le p\le \infty$ using the norm
$$ \mid\mid x\mid\mid_{p, a}:=(\sum_{\alpha} a_{\alpha}\mid x\mid^p)^{1/p}. $$ 

Our Basic Inequality says that $\Phi\subseteq l_{1, a}$ with an appropriate sequence
$a$ depending on the version of the Basic Inequality used.

{\bf Problem 6.2.} {\em  Is  $\Phi$ (relatively) compact in some natural topology on
${\mathbb{R}}^{\infty}$? }

{\bf Problem 6.3.} {\em  Does  $\Phi$ contain a non-empty open set in 
${\mathbb{R}}^{\infty}$?
Is it convex? Is it a restricted cone, i.e., does 
$  \phi_0\in \Phi $ imply $ \mu\phi_0\in
\Phi$ for $\mu \in [0,1]$? Is it true that $ \Phi \subset l_{\infty}$?} 

{\bf Problem 6.4.} {\em Show that the cardinality of $\Supp \phi$ is unbounded on $ \Phi$
where $\Supp \phi$ is the set of indices $\alpha$ with  
$\phi_\alpha\neq 0$. Does there exist
$
\phi \in \Phi$ with infinite $\Supp \phi$? }

The last question is also related with the Unramified Fontaine-Mazur conjecture, see
Subsection  6.3 below, and with the convergence abscissa for $\zeta_\phi(s)$, cf. Remark
4.1.
 
\subsection{The deficiency problem}
Since a complete description of the sets $\tilde \Phi$ and $\Phi$ is, most
probably, very difficult, one can be also interested in possible values of the deficiency,
which has his own importance:

{\bf Problem 6.5.} {\em Does there exist an infinite global field $\mathcal K$ (an asymptotically exact
family) with zero deficiency $\delta_{\mathcal K}$? }

If the answer is positive, one would like to have an explicit construction of
 such a family. Note that the positive answer is known for function global fields over ${\bf F}_r$
with a square $r$ (it is given by appropriate infinite modular function fields).

More generally, one can put the following 

{\bf Problem 6.6.} {\em Describe the set  of  possible values of the deficiency 
for infinite global fields (asymptotically exact families).  }

 One can also be interested in properties of this set: whether it is closed, convex 
(i.e., an interval), of  positive measure, etc.

 At the moment these problems seem to be inaccessible and we would like to put a
more modest question concerning amelioration of existing 
estimates for $\delta_K$.

{\bf Problem 6.7.} {\em Produce an example of an infinite global field $\mathcal K$ 
(or of an asymptotically
exact family) with the value $\delta_{\mathcal K}$ as small as possible.}

The best
example known in the number field case is that of \cite{6}, with
$\delta_{\mathcal K}\leq 0.141...$

In particular one should consider families with $\phi_q\neq 0$ for at least one 
prime power $q$. 
It looks promising to search for ramified towers  ${\mathcal K}$ of number fields with 
good  values of $\delta_{\mathcal K}$ to replace the class field towers of Section 9 below.

\subsection{Unramified Fontaine-Mazur conjecture}
Let $K$ be a number field,
 $p$ a prime, let $T$ be a finite set of primes of $K$, none above $p$, and 
let $G^{(p)}_{K,T}=\Gal({\mathcal K}^{(p)}_T)$ be the Galois group of the maximal 
algebraic pro-$p$ extension
${\mathcal K}^{(p)}_T$ of $K$ unramified outside $T$.
Then the unramified  Fontaine--Mazur conjecture reads

{\em Any continuous representation 
$$\rho: G^{(p)}_{K,T}\rightarrow \GL_n({\mathbb Z}_p) $$
has a finite image.}

One can consider a {\em just-infinite} subextension ${\mathcal L}/K$ of ${\mathcal K}^{(p)}_T$,
i.e., an extension
 which contains no proper infinite subextensions of $K$. Then Unramified Fontaine-Mazur
 Conjecture is equivalent to the finiteness of the image for any representation
 $$\rho: \Gal({\mathcal L}/K)\rightarrow \GL_n({\mathbb Z}_p) $$ 
for all just-infinite extensions ${\mathcal L}/K$.

One says that a pro-$p$ group is 
{\em torsion-riddled} if all its open subgroups have torsion.  
N. Boston \cite{1} put forth the following conjecture: 

{\em The Galois group $\Gal({\mathcal L}/K)$ is torsion-riddled for any
 just-infinite ${\mathcal L}$. }

 This conjecture would imply the unramified  Fontaine--Mazur conjecture and it 
is ultimately connected with the following 

{\bf Problem 6.8.} {\em Does there exist an infinite number field ${\mathcal K}$  for which the set
$S_{\mathcal K}$  of prime powers $q$ such that $\phi_q>0$ 
is infinite? If so, characterize such families.}

 Boston's conjecture would follow from 

{\bf Conjecture 6.1.} {\em Let  ${\mathcal K}$ be an infinite number field  which is  
 just-infinite   over $K$. Then $S_{\mathcal K}$ is infinite. }

Indeed, if it is the case, there exist in $\Gal({\mathcal K}/K)$  infinitely many Frobenius
elements  of finite order. 

Let us remark that Ihara (cf. \cite{9}, p. 695) conjectured the existence of an unramified
extension
 ${\mathcal K}$ with  $\delta_{\mathcal K}=0$ and  infinite $S_{\mathcal K}$.

\subsection{Around the asymptotic explicit formula}
One easily sees that almost all GRH results of our paper have their unconditional
counterparts, with one notable exception, namely, GRH Theorem 5.1 (and its consequences:
GRH Corollary 5.1 and GRH Theorem 5.3). It is but natural to pose 

{\bf Problem 6.9.} {\em What are   unconditional analogues of GRH Theorem 5.1., Corollary
5.1 and GRH Theorem 5.3? }

At the moment we have no approach to this problem.

\subsection{Results specific for the function field case and corresponding problems}
In the function field case we have some specific results which do not yet have their number field
counterparts. Let us discuss some  of them.
 
In this case we can get rather good estimates for the coefficients of zeta-functions. More
precisely, let 
 $$ Z_K(t)=\sum_{m=1}^{\infty} D_m t^m$$
for a function field $K$; one knows that $D_m=D_m(K)$ is the number of positive divisors 
of degree $m$. Then we have (cf. \cite{31}, Proposition 4.1 and Theorem 4.1)
 
{\bf Theorem 6.1.} {\em Let ${\mathcal K}=\{ K_j \}$ be an asymptotically exact family 
of function fields over ${\mathbb F}_r$, and let 
$$ \mu_0=\sum_{m=1}^{\infty}{m\phi_{r^m} \over r^m -1}=1-\xi_{\phi}(1). $$
Then for any real $\mu>0$ 
 $$\lim_{i\rightarrow  \infty}{D_{[\mu {g_{i}} ]}\over {g_{i}} }=\mu \log \Lambda+ 
\sum_{m=1}^{\infty}\phi_{r^m} \log{\Lambda^m\over \Lambda^m -1},$$
where $\Lambda=\Lambda(\mu)$ equals $r$ for $\mu\ge\mu_0$, and is defined from the equation 
$$\sum_{m=1}^{\infty}{m\phi_{r^m} \over \Lambda^m -1}=\mu$$
for $\mu\le\mu_0$.${\Box}$   }

Moreover, one can get even more precise result for the ratio ${D_m/h}.$
Let $h_j=h(K_j)$ be the class number. Then (cf. \cite{31}, Theorem 5.1)

{\bf Theorem 6.2.} {\em Let ${\mathcal K}=\{ K_j\}$ be  an asymptotically exact family. Then
for
any $\varepsilon >0$ and any $m$ with $m/g\geq \mu_1+\varepsilon$, we have
$$ {D_m(K_j) \over h_j}= {r^{m-g+1}\over r-1}(1+ o(1))$$
for $j\rightarrow \infty$, $o(1)$ being uniform in m. Here $\mu_1=\mu_1(\phi) $ 
is defined as the largest of the two 
roots of the equation
$${\mu\over2}+\mu\log_r{\mu\over2}+(2-\mu)\log_r(1-{\mu\over2})=
-2\log_r\zeta_{\phi}(1).{\Box}$$}

Note that the estimate of Theorem  6.2 is much more precise than that of Theorem  6.1 (there
we have an exponential $o(1)$ instead of the multiplicative one of Theorem 6.2).

{\bf Problem 6.10.} {\em What is the number field counterpart of the above results 
on the number of positive divisors?}

It does not look likely that Theorem 6.2 has a proper number field analogue. On the other hand, one
can  hope to obtain an analogue of Theorem  6.1. (Cf. Lemma 7.5 below.)

In \cite{31} we also  cosidered the asymptotic behaviour of $w_m({\mathcal K})$,
the number of classes of positive divisors of degree $m$. More precisely,
let  ${\mathcal K}=\{ K_j\}$ be  an asymptotically exact  family, and let
$$ w({\mathcal K} ,\mu)_{\inf}:=\lim \inf (w_{[ \mu g]}(K_j))^{1/g_j},$$
$$ w({\mathcal K},\mu)_{\sup}:=\lim \sup (w_{[ \mu g]}(K_j))^{1/g_j}$$
for a real number $\mu \in (0,1)$. Clearly, $$ w({\mathcal K},\mu)_{\inf}\leq
 w({\mathcal K},\mu)_{\sup}\leq d({\mathcal K},\mu)=\lim{D_{[ \mu g]}(K_j)\over
h_j}\;.$$

In \cite{31} (Proposition 6.1 and Theorem 6.1, cf. also \cite{32}) we proved that for
$\mu\in (0, 1/r)$ one has $ w({\mathcal K},\mu)_{\inf}= w({\mathcal K},\mu)_{\sup}=
d({\mathcal K},\mu)$ and that for $\mu >1/r$  the ratio
 $w({\mathcal K},\mu)_{\inf}/d({\mathcal K},\mu)$ is bounded from below by
 $ r^{-\phi_r R_r(1-\mu/\phi_r)}$ for any asymptotic upper  bound $R_r$ for $r$-ary 
linear codes
(recall that $R_r(\delta)$ is a decreasing continuous function on $[0,{r-1 \over r}]$ with
$R_r(0)=1, \:\:R_r({r-1 \over r})=0$). 

{\bf Problem 6.11.} {\em Is it true that for any $\mu\in (0, 1)$ one has $ w({\mathcal K},\mu)_{\inf}=
w({\mathcal K},\mu)_{\sup}= d({\mathcal K},\mu)$} ?

The above results use geometric arguments,  in particular,
the construction of algebraic geometry codes and have no evident number theory
counterparts.

\part {Around the Brauer--Siegel Theorem}

Part 1 was devoted to the general theory of infinite global fields and asymptotically 
exact families. In this Part we are considering a specific parameter of these fields, 
which we call the Brauer--Siegel ratio.

For an asymptotically exact
family ${\mathcal K}=\{K_i\}$ of global fields consider the limits
$$\BS({\mathcal K})=\lim_{i\rightarrow \infty}{ \log{h_i R_i} \over g_i}$$
 and $$\mathbf\kappa({\mathcal K})=\lim_{i\rightarrow \infty}{ \log{\mathbf\kappa_i}
\over g_i}.$$ Here $h_i$ is the class number, $R_i$ the regulator, and $\mathbf\kappa_i$
the zeta-residue at $s=1$.
We are going to show that these limits exist and depend only on the numbers
$\phi=\{\phi_\alpha\}$. Therefore, $\BS({\mathcal K})$ is well defined
for an infinite global field ${\mathcal K}$.

Let us start with the function field case. It was treated in our papers \cite{30} and
\cite{31},  therefore we do
not  present any proofs here. We set $R=1$ and, of course, $\phi_{\alpha}$ can be nonzero
only for
$\alpha=r^m$, $m=1, 2, 3, \dots$.
 
First of all we have (cf. \cite{30},
Corollary 2) the  following Generalized Brauer--Siegel Theorem:

{\em For an asymptotically exact family of function 
fields over ${\mathbb{F}}_r$ we have 
$$
\BS({\mathcal K})=\lim_{i\rightarrow\infty}{{\log_r {h_i}}\over{g_i}}=1 + \sum\limits_{m=1}^\infty
\phi_{r^m}\log_r{{r^m}\over{{r^m}-1}}. 
$$
}

We get the function field case analogue of the classical Brauer--Siegel theorem in the
asymptotically bad case (i.e., $\phi_\alpha=0$ for any $\alpha$). Then $\BS({\mathcal
K})=1$.

Next (cf. \cite{30}, theorem 5; \cite{31}, theorem 3.1) we have the following Bounds:

{\em For any family of function fields over ${\mathbb{F}}_r$ we have
$$
1 \le \liminf_{i\rightarrow\infty}{{\log_r{h_i}}\over{g_i}}\le
\limsup_{i\rightarrow\infty}{{\log_r{h_i}}\over{g_i}} \le 1 +
(\sqrt r -1)\log_r{r\over{r-1}}.
$$
}

We also know some partial existence results.  
Both bounds 1 and $1+ (\sqrt r-1)\log_r{r \over r-1}$ are
attainable. The lower bound 1 is attained for any asymptotically bad
family, while
$1+ (\sqrt r-1)\log_r{r \over r-1}$ is attained for any asymptotically
maximal family, i.e., such that $\phi_{r}={\sqrt r} -1$ and $\phi_{r^m}=0$ 
for all $m\not= 1$. Such families
(and even such towers) are known to exist for $r$ being a square. 

Another limit parameter $\mathbf\kappa({\mathcal K})$ gives no new information
in the function field case, since $\mathbf\kappa({\mathcal K})=\BS({\mathcal K})-1$.

In what follows we are going to present the number field analogues of 
these results which happen to be much more complicated.

\section{The Generalized Brauer--Siegel Theorem}
In this section we prove a generalization of the Brauer--Siegel theorem and present some
corollaries.
 
\subsection{Statements} 

{\bf Theorem 7.1} (Generalized Brauer--Siegel Inequality). {\em For an asymptotically
 exact family  of number fields one has
$$  \limsup_{i\rightarrow \infty}{ \log(h_iR_i) \over g_i} \leq 1+
\sum_{q}\phi_q\log{ q\over  q-1}- \phi_{\mathbb{R}}\log 2- \phi_{\mathbb{C}}\log 2\pi,$$
$$  \limsup_{i\rightarrow \infty}{ \log\mathbf\kappa_i \over g_i} \leq 
\sum_{q}\phi_q\log{ q\over  q-1},$$
the sum being taken over all prime powers $q$.}

This result is unconditional. Assuming GRH we get the equality.

{\bf GRH Theorem 7.2} (GRH Generalized Brauer--Siegel Theorem). 
{\em For an asymptotically exact family ${\mathcal K}$
of number fields the limits $\BS({\mathcal K})$ and 
$\mathbf\kappa({\mathcal K})$ exist and we have
$$\BS({\mathcal K}):=\lim_{i\rightarrow \infty}{ \log(h_iR_i) \over {g_i}} =1+
\sum_{q}\phi_q\log{ q\over  q-1}- \phi_{\mathbb{R}}\log 2- 
\phi_{\mathbb{C}}\log 2\pi,$$
$$\mathbf\kappa({\mathcal K}):=
\lim_{i\rightarrow \infty}{ \log\mathbf\kappa_i \over {g_i}} =
\sum_{q}\phi_q\log{ q\over  q-1},$$
the sum being taken over all prime powers $q$.}

If we restrict our attention to the case of {\em almost normal towers}, we can prove the
same unconditionally. To formulate the result we need one definition more.

Let $K$ be a number field. We call $K$
{\em almost normal} if there exists a finite tower of number fields ${\mathbb{Q}}=
K_0\subset K_1\subset \ldots \subset K_m=K$ such that all the extension $K_i/K_{i-1}$ 
are normal. A family is called almost normal if all its fields are. An infinite number
field is called almost normal if it is a limit of an almost normal tower.
 
 {\bf  Theorem 7.3} (Unconditional Generalized Brauer--Siegel Theorem). 
{\em For  an asymptotically good   tower 
${\mathcal K}= \{K_i\}$, $K_1\subset K_2\subset\ldots,$ of
almost normal number fields (in particular, for an infinite asymptotically good normal
number field) the limits $\BS({\mathcal K})$ and 
$\mathbf\kappa({\mathcal K})$ exist and we have
$$\BS({\mathcal K}):=\lim_{i\rightarrow \infty}{ \log(h_iR_i) \over {g_i}} =1+
\sum_{q}\phi_q\log{ q\over  q-1}- \phi_{\mathbb{R}}\log 2- 
\phi_{\mathbb{C}}\log 2\pi,$$
$$\mathbf\kappa({\mathcal K}):=
\lim_{i\rightarrow \infty}{ \log\mathbf\kappa_i \over {g_i}} =
\sum_{q}\phi_q\log{ q\over  q-1},$$
the sum being taken over all prime powers $q$.}

{\bf Remark 7.1.} The classical Brauer--Siegel theorem claims that (subject to its
conditions) $\mathbf\kappa({\mathcal K})=0$. An upper bound for
$\mathbf\kappa({\mathcal K})$  was given by Hoffstein
\cite{Hof}. We shall ameliorate on his bound below (Remark 8.2).

\subsection {Proofs}
 
{\em Proof of Theorems $7.1$, $7.2$ and $7.3$.}
We begin with the inequality of Theorem 7.1
which  does not require additional conditions. Passing to a subfamily we can suppose that
there exits the limit (may be, infinite)
$$  \lim_{i\rightarrow \infty}{ \log(h_iR_i) \over g_i}\;.$$

 For any real $s>1$ and any $K$  we have
$$\zeta _K(s)={\kappa_K\over s-1}F_K(s)\;, $$
$\kappa_K$ being the residue of $\zeta_K(s)$ at 1 and $F_K(s)$ being an analytic
function in a neighbourhood of 1 with $F_K(1)=1$. This is just a way
to write the residue.

Let us first remark that
$${\log \kappa_{K_j}\over g_j}\;\longrightarrow\;\lim_{j\rightarrow \infty}{ \log(h_jR_j) \over g_j}
-1+ \phi_{\mathbb{R}}\log 2+ \phi_{\mathbb{C}}\log 2\pi\;.$$
To see this, start with the standard formula
$$  \kappa_K={2^{r_1}(2\pi)^{r_2}(h_K R_K)\over w_K \sqrt {\vert D_K\vert}}\;. $$
For $g_j=\log \sqrt{\vert D_{K_j}\vert}\to\infty$,
i.e., for $j\rightarrow \infty$, its logarithm gives 
exactly what we want, if we note that $\log w_{K_j} / \log{\vert D_{K_j}\vert} \rightarrow 0$
since $w_{K_j}\le cn_{K_j}^2$
for an absolute constant $c$ (cf. e.g., \cite{13}, proof of Lemma 1 of XVI.1).

Let us put $s=1+\theta$ with $\theta=\theta_j >0$ being  a small positive real number, its
dependence on $j$
to be specified
later. Taking the logarithm of $$\zeta_{K_j}(s)={\kappa_{K_j}\over s-1}F_{K_j}(s)$$
and dividing by $g_j$  we get
$${\log \zeta_{K_j}(1+\theta_j) \over g_j}={\log\kappa_{K_j}\over g_j}+
{\log F_{K_j}(1+\theta_j) \over g_j}-{\log\theta_j \over g_j}\;. $$
If $j\to\infty$ then, to prove the theorem, it suffices to show
that for a proper choice of  $\theta_j$ the following three points are satisfied:

(i) $${\log\zeta_{K_j}(1+\theta_j) \over g_j}\longrightarrow\sum_{q}\phi_q\log{q\over q-1};$$

(ii) $${\log\theta_j \over g_j}\longrightarrow 0;$$

(iii) $$\lim\inf\; {\log F_{K_j}(1+\theta_j) \over g_j}\ge 0.$$

Let us first look at (i). We have
 
$$\zeta_{K_j}(1+\theta)=\prod_q (1-q^{-1-\theta})^{-N_j(q)} $$
for any $\theta>0$, where $N_j(q)$ is the number of places of $K_j$ with the norm $q$.
Let $$f_j(\theta)={\log\zeta_{K_j}(1+\theta)\over g_j}$$ and
$$f(\theta):=\sum_q \phi_q\log{1 \over 1-q^{-1-\theta}}.$$
Taking logarithms and dividing by $g_j$ we get 
$$f_j(\theta)=
\sum\limits_q {N_j(q)\over g_j}\log{1 \over 1-q^{-1-\theta}}, $$
thus $f_j(\theta)\to f(\theta)$  uniformly for $\theta \ge\theta_0>0$
by definition of $\phi_q,\;$  ${N_j(q)/ g_j}$ being bounded by an absolute constant and
$\sum\limits_q {N_j(q)\over g_j}\log{1 \over 1-q^{-1-\theta}}$ converging for $\theta>0$.
Moreover,
$$f(\theta)\;\longrightarrow\;\sum_{q}\phi_q\log{q\over q-1}\quad\hbox{ for }\quad\theta\to 0\; ,$$
the series $\sum\limits_{q}\phi_q\log{ q\over  q-1}$ being convergent (the series
$\sum\limits_{q}{\phi_q\log{q}\over q-1}$ which is at most 1 by Proposition 3.2 
provides an upper bound for it).

Then we choose a decreasing sequence $\theta(N)>0$ in such a way that  $$\vert f(\theta(N))-
\sum_{q}\phi_q\log{ q\over  q-1}\vert<1/2N\;,$$
 and we can also choose an increasing sequence $j(N)$ such that  $g_{j(N)}\ge{1\over \theta(N)}$
and also such that for any $\theta\in\lbrack\theta(N+1),\theta(N)\rbrack $ we have
$$\vert f(\theta)- f_{j(N)}(\theta)
\vert< 1/2N.$$ 
This is possible since $g_j\to\infty$  and $f_j(\theta)\to f(\theta)$
uniformly for
$\theta\ge\theta(N+1).$
Then let $N=N(j)$ be given by $j(N)\le j\le j(N+1)-1$ and put
 $ \theta_j=\theta(N(j))$; note that $N(j)\to\infty$. We see that
$$\vert f_j(\theta_j)-\sum_{q}\phi_q\log{ q\over  q-1}\vert<1/N(j)\quad\longrightarrow\quad 0$$
which proves  (i). We also get (ii) for granted since
$1/g_j\le1/g_{j(N)}\le\theta_{j}$ and hence $\log{\theta_j}/g_j\to 0.$

For (iii), keeping in mind that 
$\left({\log \zeta_{K_j}(s)\over g_j}\right)'
=\sum_{  P}\sum_{m=1}^{\infty} r^{-ms}\log r,$
we rewrite  Stark's formula 
$$\log{\vert D\vert}=r_1(\log\pi-\psi(s/2))+2r_2(\log(2\pi)-\psi(s))-{2\over s}-{2\over{s-1}}$$
$$+2{\sum_{\rho}}'{1 \over s-\rho}+2\sum_{  P}\sum_{m=1}^{\infty}
r^{-ms}\log r,$$
for 
$s=1+\theta$ as
$$\left({\log \zeta_{K_j}(1+\theta)+\log \theta \over g_j}\right)'=-1+{r_1\over 2g_j}
(\log \pi-\psi({1+\theta\over 2}))$$
$$+{r_2\over g_j}(\log 2\pi-\psi({1+\theta}))-
{1\over (1+\theta)g_j}+{\sum_{\rho}}'{1\over (1+\theta-\rho)g_j},$$
which shows that the derivative
$$\left({\log F_{K_j}(1+\theta)  \over g_j}\right)'=
\left({\log \zeta _j(1+\theta)+\log \theta \over g_j}\right)'$$
 is bounded from below by $-2$ for any small enough $\theta$ since all the terms except $-1$
and $-{1\over (1+\theta)g_j}$ are positive.
 Thus
$${\log  F_{K_j}(1+\theta_j)  \over g_j} \ge c\theta_j\longrightarrow 0\;,$$
 which proves (iii). Summing up, we get an unconditional proof of
 the desired inequality
$$\lim_{i\rightarrow \infty}{ \log(\mathbf\kappa_i) \over g_i} \le 
\sum_{q}\phi_q\log{ q\over  q-1}$$
and the corresponding one for $\lim{ \log(h_iR_i) \over g_i}$,
i.e., that of Theorem 7.1.

\vskip 0.1 cm

To prove Theorem 7.2 one supposes GRH. In fact, it is sufficient
to prove that

(iii)$'$
$$\lim\inf {\log  F_{K_j}(1+\theta_j) \over g_j}\le 0.$$

To do this we shall use the following GRH lemma.

{\bf GRH Lemma 7.1.} {\em For any asymptotically exact family
 of number fields the function
$$Z_j(s):={-(\log( \zeta_{K_j}(s))'-1/(s-1) \over g_j}$$
 tends for
$j\rightarrow\infty$ to 
$$Z_{\phi}(s):=\sum\limits_{q}\phi_q  {{\log q} \over q^{s}-1}$$
 uniformly on $\Re(s)\ge 1/2+\delta$ for any $\delta>0$.}
 
In fact, Lemma 7.1 is the key lemma of Ihara's paper (\cite{9}, p. 698), 
where it is proved in the special case
of an unramified tower; his proof stays mostly valid
in our situation as well, nevertheless we present it here.
 
{\em Proof of Lemma $7.1$.} Note first of all, that the series defining $Z_{\phi}(s)$
converges uniformly on  $\Re s\ge 1.$ Indeed, it is bounded from above by
$$\sum\limits_{q}\phi_q{{\log q}\over q-1}\le 1.$$
If one assumes GRH, the
series becomes uniformly
 convergent and hence analytic  on  $\Re s>1/2$, since it is the case for 
$$\sum\limits_{q}\phi_q{\log q   \over\sqrt q-1}\le 1.$$

Let us consider a presentation of $Z_j(s)$ and $Z_{\phi}$ as Mellin
transforms of Chebyshev step functions.
We have a well-known and easy to prove formula (cf. \cite{9}, eq. 5-5 and 5-6)
valid for $\Re s>1$:

$$ s^{-1}Z_j(s)={1 \over g_j}\int_1^{\infty}(G_j(x)-x) x^{-s-1} dx\;,$$
where
$$G_j(x):= \sum_{{{  P}, \: m\ge 1}\atop {N(P)^m\le x}} \log N({  P})=\sum_
{{q,\:m\ge1}\atop {q^m\le x}}N_q(K_j) \log q$$
is the Chebyshev step function for the field $K_j$, and the first sum is taken 
over all prime divisors
$  P$ of the field $K_j$.  

Similarly, for $Z_{\phi}$ we get for ${\hbox {Re}}(s)>1$

$$s^{-1}Z_{\phi}(s)=\int_1^{\infty}H(x) x^{-s-1} dx,  $$
where $H(x)$ is an asymptotic analogue of $G_j(x)$:
$$H(x)=\sum_{{q,\:m\ge1}\atop {q^m\le x}} \phi_q \log q. $$

Now we use the Lagarias-Odlyzko estimate for $G_j(s)$ (which uses GRH, cf. \cite{11},
Theorem 9.1):
$$ \vert G_j(x)-x \vert \le C (n_j\sqrt x (\log x)^2+2g_j\sqrt x \log x)\;,$$
where $n_j=[K_j:{\mathbb{Q}}]$ and $C$ is an absolute constant. Thus
$$ \vert G_j(x)-x \vert \le C_1g_j\sqrt x (\log x)^2$$
 with another absolute constant $C_1$.

The last formula shows that the integral in the integral representation of $s^{-1}Z_j(s)$ 
converges for ${\hbox {Re}}(s)>1/2$, and thus the representation is valid 
for ${\hbox {Re}}(s)>1/2$. The same is true for $s^{-1}Z_{\phi}(s)$ since it is analytic
for ${\hbox {Re}}(s)>1/2$ as explained above. Therefore, for ${\hbox {Re}}(s)>1/2$ we get 
$$s^{-1}Z_j(s)-s^{-1}Z_{\phi}(s)=\int_1^{\infty}\left({G_j(x)-x \over g_j}-H(x)\right) x^{-s-1} dx\;. $$

Fix $\delta>0$. Let then ${\hbox {Re}}(s)\ge 1/2+ \delta $ and let $\varepsilon>0$. 
We choose $M>1$ so that
$$C_1\int_M^{\infty}(\log x)^2 x^{-1-\delta}dx \le \varepsilon$$ and
$$ \int_M^{\infty}H(x) x^{-{3\over 2}-\delta}dx \le \varepsilon\;.$$
Let then choose $j(M)$  in such a way that for $j\ge j(M)$ we have the following two inequalities:
$$
\big\vert {G_j(x) \over g_j}-H(x)\big\vert \le \delta\varepsilon {\hbox  { for }} 1\le x\le M 
$$
and
$$
\big\vert\int_1^{M}\left({x \over g_j}\right) x^{-s-1} dx\big\vert\le
\int_1^{M}\left({x \over g_j}\right) x^{-\delta-{3\over 2}} dx=
{M^{{{1\over 2}-\delta}}-1\over g_j({1\over 2}-\delta)}\le \varepsilon\;, 
$$
which is possible since $N_q(K_j)/g_j$ tends to $\phi_q$, and since the sums in the definition
of $G_j(x)$ and $H(x)$ contain only finite (and bounded from above)  number of terms 
for $x\le M$. Here, by abuse of notation, we agree to understand
$(M^{{1\over 2}-\delta}-1)/
({1\over 2}-\delta)$ as $\log M$ if $\delta=1/2.$

We get 
$$\vert s^{-1}Z_j(s)-s^{-1}Z_{\phi}(s)\vert\;\; =\;\;
\big\vert \int_1^{\infty}\left({G_j(x)-x \over g_j}-H(x)\right) x^{-s-1} dx\big\vert\;\;= $$
$$\big\vert \int_1^{M}\left({G_j(x) \over g_j}-H(x)\right) x^{-s-1} dx-
\int_1^{M}\left({x \over g_j}\right) x^{-s-1} dx+$$
$$\int_M^{\infty}({G_j(x)-x \over g_j}-H(x))x^{-s-1} dx
\big\vert\;\; \le $$
$$ \delta \varepsilon\int_1^{M} x^{-\delta-{3\over 2}} dx +
 \int_1^{M}\left({x \over g_j}\right) x^{-\delta-{3\over 2}} dx 
+\big\vert\int_M^{\infty}\left({G_j(x)-x \over g_j}\right) x^{-s-1} dx\big\vert+$$

 $$ \int_M^{\infty}H(x) x^{-s-1} dx \le {\delta \varepsilon\over \delta+{1\over 2}}+
\big\vert{M^{{{1\over 2}-\delta}}-1\over g_j({1\over 2}-\delta)}\big\vert+
C_1\int_M^{\infty}(\log x)^2 x^{-1-\delta}dx+$$
$$\int_M^{\infty}H(x) x^{-{3\over 2}-\delta}dx
 \le 4\varepsilon  $$
for $\Re(s)\ge 1/2+ \delta $ and  $j\ge j(M)$ which proves the lemma. Note that
$j(M)$ depends on $\delta$.$\;{\Box}$

{\em End of proof of GRH Theorem} 7.2. From Lemma 7.1 it follows, in particular,  that 
$Z_j(s)$ tends to $Z_{\phi}(s)$ for
$j\rightarrow\infty$ uniformly on ${\hbox{Re}}(s)\ge 1$.
Therefore for small enough $\varepsilon$, large enough $j$, and any $\theta>0$, we have 
$$ \vert ({\log  F_{K_j}(1+\theta) / g_j})'\vert=
\vert Z_j(1+\theta)\vert \le \vert Z_{\phi}(1)\vert+\varepsilon \le 1,$$
since $\vert Z_{\phi}(1)\vert<1$ because of 
Proposition 3.2. Thus 
$${\log  F_{K_j}(1+\theta_j) / g_j} \le\theta_j.$$
This proves (iii)$'$
and the theorem.${\Box}$ 

{\bf Remark 7.2}. In fact, one notes that Corollary 7.1 remains valid under the assumption 
 that there are no zeta-zeroes with ${\hbox {Re}}(s)\ge 1-\delta$ for arbitrary fixed
$\delta>0$, so that we do not need the full strength of GRH.

To prove the opposite inequality of Theorem 7.3 it is sufficient to show that for 
an asymptotically  good almost normal tower  $\{K_j\}$ and for a 
proper choice of  $\theta_j$
the following conditions are satisfied:

 $$\quad\quad{\hbox{(i)}}'' \quad \quad \lim\inf\;
{f_j(\theta_j)}\geq\sum_{q}\phi_q\log{q\over q-1}=f(0)\; ;  \quad $$

$$\quad\quad{\hbox{(ii)}}''\quad \quad {\log\theta_j   \over g_j}\longrightarrow 0\; ;\quad
\quad \quad \quad \quad \quad \quad \quad \quad \quad \quad \quad \quad $$

$$\quad\quad{\hbox{(iii)}}''\quad \quad  \limsup {\log  F_{K_j}(1+\theta_j) / g_j}\le
0\;.\quad\quad\quad\quad\quad $$

We should stress that the choice of $\theta_j$ in the proof of Theorem  7.3 is completely
different from that in the proofs of Theorems 7.1 and 7.2. 
We need the following definition.

Let $K$ be a number field (of a finite degree $n$). A real number $\rho$
is called an {\em exceptional zero }  of 
$\zeta_K(s)$ if $\zeta_K(\rho)=0$ and 
$$ 1-(4\log \vert D_K \vert)^{-1}\le \rho<1; $$
an  exceptional zero $\rho$ of $\zeta_K(s)$  is called its {\it Siegel zero} if
$$ 1-(16\log \vert D_K \vert)^{-1}\le \rho<1. $$

It is known that  for any $K$ there exists at most one exceptional zero of
$\zeta_K(s)$.  

We are going to show that  under some conditions asymptotically exact families  
have no Siegel zero. We begin with 
the following fundamental property of 
Siegel zeroes discovered by Heilbronn \cite{7} and precised  by Stark 
(\cite{27}, Lemma 10):

{\bf Lemma 7.2.} {\em Let $K$ be an  almost normal number field,   and let $\rho$
 be a Siegel zero of $\zeta_K(s)$.
Then there is a quadratic subfield $k$ of $K$ such that $\zeta_k(\rho)=0.\;{\Box}$}

{\bf Lemma 7.3.} {\em Let ${\mathcal K}= \{K_i\}$ be an asymptotically good   
 family of almost normal number fields. Then there exists a positive integer  $I$ such that 
 $\zeta_{K_i}(s)$ has no Siegel zero for any $i\ge I$.  
In other words, in such a family almost all fields  have no Siegel  zero.}
 
{\em Proof.} In view of Lemma 7.2 it is sufficient to prove that the set 
$Q({\mathcal K})$ of quadratic 
fields $k$ contained in at least one of the fields $K_i$ is  finite. 
Indeed, if this is the case, let 
$$\beta =\max\{\rho\in{\mathbb{R}}: \zeta_k(\rho)=0 {\hbox { for }} k\in Q({\mathcal K})\}. $$
Since $Q({\mathcal K})$ is finite, the maximum exists and $\beta<1$. Now if 
$$g_{K_i}>{1\over 16(1-\beta )}$$
 then $\zeta_{K_i}(s)$ has no Siegel zeroes by Lemma 7.2. 

Let us verify the 
finiteness  of $Q({\mathcal K})$. The ratio $n/g$ is non-increasing 
in extensions, since 
$$\vert D_{K}\vert\ge\vert D_{k}\vert^{[K:k]}. $$
Moreover,  $n_{K_i}/g_{K_i}\to \phi_\infty=\phi_{\mathbb{R}}+2\phi_{\mathbb{C}}>0$, 
the family ${\mathcal K}$ being asymptotically good. Therefore,
there exists a positive real number $\varepsilon$ such that
$n_{K_i}/g_{K_i}\ge\varepsilon$ for any $i$. If $k\in Q({\mathcal K})$,  $k\subseteq K_i$ then  
$$2/g_k=n_k/g_k \ge n_{K_i}/g_{K_i}\ge \varepsilon .$$
 Therefore, $g_k\le  2/\varepsilon$ for any $k\in Q({\mathcal K})$, and $\vert D_k \vert \le 
e^{4/\varepsilon}$, which  implies  the finiteness of $ Q({\mathcal K})$. ${\Box}$

{\bf Corollary 7.1.} {\em Let  $K_1$ be a number field with infinite Hilbert class 
field tower $ \{K_i\}$. Then almost all fields $K_i$ have no Siegel zero. ${\Box}$}

Note that for any $\theta>0$
$$\left({\log F_{K_j}(1+\theta)  \over g_j}\right)'=Z_j(1+\theta)$$
where
$$Z_j(s):={-(\log( \zeta_{K_j}(s))'-1/(s-1) \over g_j}$$
which follows from the definition of $F_{K_j}$.
  
{\bf  Lemma 7.4.} {\em There exist absolute constants 
$C_0$ and $C>0$, such that for any number field $K$ which
has no  Siegel zero we have 
$$\vert Z(1+\theta)\vert\le C g^6$$
for any  real
$\theta\in (0,1)$ and for  any $g>C_0$. Here
$$Z(s):={-(\log( \zeta_{K}(s))'-1/(s-1) \over g}.$$ 
}

{\em Proof of Lemma $7.4.$} 
We use the above  presentation of $Z(s)$ 
 as the Mellin transform of the Chebyshev step function:

$$ s^{-1}Z(s)={1 \over g}\int_1^{\infty}(G(x)-x) x^{-s-1} dx\;,$$
where
$$G(x):= \sum_{{{  P}, \: m\ge 1}\atop {N(P)^m\le x}} \log N({  P}) 
=\sum_{{q,\:m\ge1}\atop {q^m\le x}}N_q(K) \log q$$

 We use then the (unconditional) Lagarias-Odlyzko estimate for $G(x)$  ( \cite{11}, Theorem 9.2):
$$\vert G(x)-x \vert \le C_1  x\exp\left(-C_2\sqrt{{\log x\over n}}\right)+ 
{x^{\rho}\over \rho}
\;$$ for $\log x\ge C_3ng^2,$ where $n=[K:{\mathbb{Q}}]$ , $C_1,C_2$ and $C_3$ being positive   
absolute constants, $\rho$ being an eventual exceptinal zero of $K$; note that since $\rho$ is
not a Siegel zero we can suppose that $ 1-(16g)^{-1}> \rho\ge 1-(4g)^{-1}$. Under that
condition one easily verifies, using that $g\ge C_4n$ for a positive   absolute constant
$C_4$, the condition
${x^{\rho}/\rho}=o\left(x\exp\left(-C_2\sqrt{{\log x/ n}}\right)\right)$, and we can 
suppose that
$$\vert G(x)-x \vert \le C_1  x\exp\left(-C_2\sqrt{{\log x\over n}}\right)\; .$$
Since
$g\ge C_4n$ for a positive   absolute constant
$C_4$ we also have 
$$\vert G(x)-x \vert \le C_1  x \exp\left(-C_5\sqrt{{\log x\over g}}\right) \;$$
for $\log x\ge C_6g^3$ and positive absolute constants $C_5$ and $C_6$.
Note that for $\log x\le C_6g^3$ we have the following  trivial 
estimate
 $$ 0\le G(x)  \le C_7 g x\log x 
$$ 
with an absolute constant $C_7$; indeed
$$G(x)=  \sum_{{q,\:m\ge1}\atop {q^m\le x}}N_q(K) \log q \le n \sum_{{q,\:m\ge1}\atop {q^m\le
x}} \log q\le C_7 g x\log x ,
$$
since $n\le C g$ and $\sum_{q^m\le x}  \log q\le C' x\log x $.

Therefore,
$$ \left \vert {Z(1+\theta)\over 1+\theta} \right \vert = 
\left \vert {1 \over g}\int_1^{\infty}(G(x)-x)
 x^{-2-\theta} dx\ \right \vert =\;$$
$${1 \over g}\left \vert\int_1^{\exp( C_6g^3)}(G(x)-x) x^{-2-\theta} dx+
\int_{\exp( C_6g^3)}^{\infty}(G(x)-x) x^{-2-\theta} dx \right  \vert  \le$$

$${(C_7 +1) }   \int_1^{\exp( C_6g^3)} x^{-1-\theta}\log x  dx    +
{C_1 \over g}  \int_{ \exp(C_6g^3)}^{\infty} \exp\left( -C_5\sqrt{{\log x
 \over g}}\right)
 x^{-1-\theta} dx.          $$
 
Then we have 
$${(C_7 +1) } \int_1^{\exp( C_6g^3)} x^{-1-\theta} \log x dx  \le C_6g^3{(C_7 +1)\over
\theta}(1-e^{-\theta C_6g^3})\le (C_7 +1)C_6^2g^6$$
and
$${C_1 \over g}  \int_{\exp( C_6g^3)}^{\infty} \exp\left( -C_5\sqrt{{\log x
 \over g}}\right) 
x^{-1-\theta} dx\   =2C_1   \int_{\exp( g\sqrt {C_6})}^{\infty} 
z^{-\theta g\log z -C_5-1}\log z  dz    $$
which can be seen by the change of variables $x=z^{g\log z}$.

Since
 $$z^{-\theta g\log z}\le z^{-\theta g^2\sqrt{C_6}}$$
and 
$$ \log z\le z^{g^{-1}C_6^{-1/2}\log(g\sqrt{C_6})}$$
for $z\ge \exp( g\sqrt {C_6})$ and $g>C_0$, for
$$ \alpha(g)= g^2\sqrt{C_6}, \quad \beta(g)=g^{-1}C_6^{-1/2}\log(g\sqrt{C_6})	$$
we get the following estimate
$$C_1  \int_{\exp( g\sqrt {C_6})}^{\infty} z^{-\theta g\log z -C_5-1}
\log z  dz 
  \le 2 C_1    \int_{\exp( g\sqrt {C_6})} ^{\infty} 
z^{-\theta\alpha(g)+\beta(g) -C_5-1} dz =$$	
$$ {2 C_1 \over (\theta\alpha(g)-\beta(g) +C_5)} \exp(- g
 (\theta\alpha(g)-\beta(g) +C_5))\sqrt {C_6}\le  C_8  \exp(- C_9 g) 
\le  C_{8} $$
with  positive absolute constants $C_8 $ and $C_9$, which implies the lemma  since
$(C_7+1)C_6^2g^6+C_{8}\le C g^6$. 
${\Box}$			
 
Now let us set $\theta_j=g_j^{-7}$, and verify condititions (i)$''$,
(ii)$''$ and (iii)$''$ for that choice, which achieves the proof of Theorem 7.3. The
conditition (ii)$''$ is obvious.

 Applying   Lemmata 7.3 and 7.4
 to any field $K=K_j$ from our tower for large enough $j$  we get
$$\left \vert \left({\log F_{K_j}(1+\theta)  \over g_j}\right)'\right \vert
=\vert Z_j(1+\theta)\vert 
\le C  g_j^6$$
for any $\theta\in (0,1)$.  

Therefore, recalling that $F_{K_j}(1) =1$, we get
$$\left \vert{\log F_{K_j}(1+\theta_j)  \over g_j}\right \vert =\left \vert
 \int_0^{\theta_j}
\left({\log F_{K_j}(1+\theta)  \over g_j}\right)' d\theta \right \vert  \le C  g_j^6
\theta_j ={C\over g_j},$$
which proves (iii)$''$.
 
Let us prove inequality (i)$''$.  
We set $$f_j(\theta)=f_j^{(1)}(\theta)+f_j^{(2)}(\theta),$$
where $$f_j^{(1)}(\theta):= 
\sum\limits_p {N_j(p)\over g_j}\log{1 \over 1-p^{-1-\theta}}$$
is the sum over prime $p$, and
$$f_j^{(2)}(\theta):= 
\sum\limits_{q=p^m, m\geq 2} {N_j(q)\over g_j}\log{1 \over 1-q^{-1-\theta}}.$$

Similarly, we set 
$$f_j^{(1)}(\theta):= 
\sum\limits_p \phi_p\log{1 \over 1-p^{-1-\theta}},$$ 
$$f_j^{(2)}(\theta):=  
\sum\limits_{q=p^m, m\geq 2} \phi_q\log{1 \over 1-q^{-1-\theta}},$$
$$f (\theta)=f_j^{(1)}(\theta)+f_j^{(2)}(\theta).$$

Since for a prime $p$ and any $j$  one has $\phi_p\le {N_j(p)\over g_j}$,
we get $f_j^{(1)}(\theta)\ge f^{(1)}(\theta)$ for any $\theta>0$. On the other hand,
$f_j^{(2)}(\theta)$ and $f^{(2)}(\theta)$ converge uniformly on $\theta \ge -\delta$ with a
positive $\delta$, and thus $f_j^{(2)}(\theta_j)$ tends to $f^{(2)}(0)$ for $\theta_j$ tending
to zero. We get 
$$\lim\inf\;
{f_j(\theta_j)}=\lim\inf (f_j^{(1)}(\theta_j)+f_j^{(2)}(\theta_j))\geq
\lim\inf (f^{(1)}(\theta_j)+f^{(2)}(\theta_j))=f(0).{\Box}$$

\subsection {Lower bounds for regulators}

As an application of the Generalized Brauer-Siegel Theorem 
one obtains a lower bound for 
regulators of number fields in asymptotically good families,  
which is better than
the general bound obtained by Zimmert \cite{Zi}.

{\bf GRH Theorem 7.4.} {\em For an asymptotically good family
of number fields
$$\liminf_{i\rightarrow \infty} {\log R_i\over g_i} \ge 
(\log \sqrt{\pi e} +\gamma/2)\phi_{\mathbb R}+(\log 2 +\gamma)\phi_{\mathbb C}.$$}

{\em Proof.}
We begin with an estimate for the class numbers of fields 
in question which could be of
independent interest. 

{\bf Proposition 7.1.} {\em   For an asymptotically exact family
 of number fields
$$\limsup_{i\rightarrow \infty}{\log h_i\over g_i}\le  1-
(\log 2\sqrt \pi+{\gamma+1\over 2})\phi_{\mathbb
R}-(\log 4\pi+\gamma)\phi_{\mathbb C}+\sum\phi_{q}
\log {q\over q -1}.$$ 
}

{\em Proof of Proposition 7.1.}
Let  ${\mathcal K}=\{K_i\}$ be an asymptotically exact family of number fields and let 
$$\zeta_{K_i}(s)=\sum_{n=1}^{\infty} D_{n}^{(i)} n^{-s}$$
be the  corresponding zeta functions.  We shall use the following result on the asymptotic
behaviour of the coefficients $D_{n}^{(i)}$: 

{\bf Lemma 7.5.} {\em  Let $n_i, \:i=1,2,\dots,$ 
be a sequence of positive integers such
that the limit 
$$\nu := \lim_{i\to\infty} {n_i\over g_i{}} $$
exists. Then
$$ \limsup_{i\rightarrow \infty} {\log D_{n_i}^{(i)}\over g_i{}}=
 \nu\limsup_{i\rightarrow \infty} {\log D_{n_i}^{(i)}\over n_i}\le
 \nu+\sum_q\phi_q\log{q\over q-1}, $$
where the sum is taken over all prime powers.}
    
{\em Proof of Lemma 7.5.}
Indeed, from the  Euler product for $\zeta_{K_i}(s)$ we see:
 
a$)$ the function $n\rightarrow D_{n}^{(i)}$ is multiplicative, i.e.,
$D_{nn'}^{(i)}=D_{n}^{(i)}D_{n'}^{(i)}$ for coprime
$n$ and $n'$. In particular,
$$D^{(i)}_n=\prod D^{(i)}_{p_j^{m_j}}$$
$n=\prod p_j^{m_j}$ being the prime factorization of $n$.

b$)$
 $$ D^{(i)}_{p^m}=\sum_{{(b_1,\dots,b_m)}\atop{b_1+2b_2+\ldots+mb_m=m}}
\prod_{s=1}^m{N_{p^s}(K_i)+b_s-1\choose b_s},$$
the sum being taken over all partitions 
of $m$, $b_i\in\mathbb Z$, $b_i\ge 0$.

Let $n_i=\prod p_j^{m_{ij}}$. Then
$$ \log D_{n_i}^{(i)}=\sum_j \log D_{p_j^{m_{ij}}}^{(i)},$$
$$ D_{p_j^{m_{ij}}}^{(i)}=\sum_{{(b_1,\dots,b_{m_{ij}})}\atop{b_1+2b_2+\ldots+{m_{ij}}
b_{m_{ij}}=
{m_{ij}}}}\prod_{s=1}^{m_{ij}}{N_{p_j^s}(K_i)+b_s-1\choose b_s}.$$
This implies
$$\log D_{p_j^{m_{ij}}}^{(i)}\le \log {\bf p}(m_{ij})+ \max_{{(b_1,\dots,b_{m_{ij}})}\atop{b_1+2b_2+\ldots+{m_{ij}}
b_{m_{ij}}=
{m_{ij}}}} 
\left ( \sum_{s=1}^{m_{ij}}\log {N_{p_j^s}(K_i)+b_s-1\choose
b_s}\right),$$ 
where ${\bf p}(x)$ is the partition function. 
Now the argument of the proof of Lemma 3.4.10 of \cite{Ts/Vl 1}   shows that 
$$\log D_{p_j^{m_{ij}}}^{(i)}\le O(\sqrt g)+ \log {p_j^{m_{ij}}} + \phi_{p_j^{m_{ij}}} \log {{p_j^{m_{ij}}}\over
{p_j^{m_{ij}}} -1}, $$ which proves the lemma.$\Box$

{\em End of proof of  Proposition 7.1.} Zimmert's theorem on twin classes
\cite{Zi} (cf. also \cite{Oes}) states that for $\lambda_{\mathbb R}=\log 2\sqrt \pi+{\gamma+1\over 2}$, 
$\lambda_{\mathbb C}=\log 4\pi+\gamma$, and any class
${\mathcal C}$ of ideals of a number field $k$
$${n_{\inf}({\mathcal
C})+n_{\inf}({\mathcal C^*})\over 2}\le 
g-\lambda_{\mathbb R}r_1-\lambda_{\mathbb C}r_2
$$ 
and thus
$${\nu_{\inf}({\mathcal
C})+\nu_{\inf}({\mathcal C^*})\over 2} \le 1-\lambda_{\mathbb R}\phi_{\mathbb
R}-\lambda_{\mathbb C}\phi_{\mathbb C},$$  where $\nu_{{\inf}}({\mathcal
C})={n_{\inf}({\mathcal C})/ g}$, $n_{\inf}({\mathcal C})$ 
being the minimum norm of an
ideal from  
${\mathcal C}$, and ${\mathcal C^*}$ being the twin class of 
the class ${\mathcal C}$
(i.e., the class 
${\mathcal D}{\mathcal C}^{-1}$, where ${\mathcal D}$ 
is the class of the different). This implies
that the class number $h_i$ of $K_i$ is  not greater 
than two times the number of "small norm" ideals, i.e., those counted
in $D_{n'_i}^{(i)}$ for $n'_i\le {\tilde n}_i=g_i(1-\lambda_{\mathbb R}\phi_{\mathbb
R}-\lambda_{\mathbb C}\phi_{\mathbb C})+o(g_i)$. By Lemma 7.5, 
$$\limsup_{i\rightarrow \infty}{\log h_i\over g_i}\le
\limsup_{i\rightarrow \infty}{1\over {g_i}}\log ({2\sum_{n'_i\le {\tilde n}_i}D_{n'_i}^{(i)}})
\le\limsup_{i\rightarrow \infty}{1\over {g_i}}\log ({\max_{n'_i\le {\tilde n}_i}D_{n'_i}^{(i)}})$$
$$\le 1-\lambda_{\mathbb R}\phi_{\mathbb
R}-\lambda_{\mathbb C}\phi_{\mathbb C}+\sum\phi_{q}
\log {q\over q -1}.\Box$$ 

{\em End of proof of  Theorem 7.4.} To finish the proof one compares the
Generalized Brauer-Siegel theorem with Proposition 7.1.
$\Box$

Passing to an asymptotically exact subfamily we easily deduce  
 
{\bf GRH Corollary 7.3.} {\em Let ${\mathcal S}$ be any family of 
number fields which does not contain an asymptotically bad subfamily. 
Then there
exists a strictly positive $A=A({\mathcal S})$ such that
$$ R(K)\ge A (\sqrt \pi e^{(\gamma+1)/2})^{r_1}(2e^{\gamma})^{r_2}$$
for any $K\in {\mathcal S}.\Box$}

{\bf Remark 7.3.} Let us recall that the resut of Zimmert 
(implicit in \cite{Zi}, but easily deduced
from the argument therein) in our notation reads
$$\liminf_{i\rightarrow \infty} {\log R_i\over g_i} \ge 
(\log 2   +\gamma)\phi_{\mathbb R}+  2 \gamma\phi_{\mathbb C}.$$

The numerical values of Zimmert's coefficients are 
$\log 2   +\gamma\approx 1.270\ldots$ and 
$2\gamma\approx 1.154\ldots$; those of ours being
$\log \sqrt{\pi e} +\gamma/2\approx 1.361\ldots$ and 
$\log 2   +\gamma\approx 1.270\ldots$, respectively.
 
Applying the same argument to the case of an asymptotically good   
tower of almost
normal number fields we get an unconditional
version of Theorem 7.4:
 
{\bf  Theorem 7.5.} {\em For an asymptotically good tower
 of  almost
normal number fields 
$$\liminf_{i\rightarrow \infty} {\log R_i\over g_i} \ge 
(\log \sqrt{\pi e} +\gamma/2)\phi_{\mathbb R}+(\log 2 +\gamma)\phi_{\mathbb C}. \Box$$
}

\section{Bounds for the Brauer--Siegel Ratios}

\subsection{Linear programming problem}

To get optimal
estimates (on both
sides) for the limit points of
$\log(h_iR_i) \over g_i$ and $\log(\mathbf\kappa_i) \over g_i$
we come to the following linear programming problem.

Let $q$ run over all prime powers and let there be given two sets of
non-negative coefficients $\lbrace a_q \rbrace$ and $\lbrace b_q \rbrace$
as well as non-negative
$a_0, b_0, a_1, b_1$ with the properties that if $a_i=0$ then $b_i=0$ for
all $i=0,1,q$.
Suppose that

(1)
$${m\over n} \ge {a_{p^m}\over a_{p^n}}\;,$$
for any $m\ge n$;
and

(2)
$${{b_{q_1}\over a_{q_1}} \ge {b_{q_2}\over a_{q_2}}}$$
for any $q_1 \le q_2$
such that $a_{q_1}\not =0,\; a_{q_2}\not =0$,

(3) if $a_0\not =0$ and $a_1\not =0$ then $a_1\ge a_0,\; b_1\ge b_0$,
and
$${b_0\over a_0} \le {b_1\over a_1}\;, $$

(4)
$$
\sum_q b_q = \sum_q a_q = \infty\;.
$$

Consider the
following optimization problem:

{\em Find the maximum and the minimum of
$$F(x)=\sum_q b_q x_q - b_0 x_0 - b_1 x_1$$
under the conditions

(i) for any $i=0, 1,\hbox{or } q$
$$x_i\ge 0;$$

(ii) for any prime $p$
$$ \sum_{m=1}^\infty m x_{p^m} \le x_0+2 x_1;$$

(iii)
$$ \sum_q a_q x_q + a_0 x_0 + a_1 x_1 \le 1;$$

(iv) if for some $i=0, 1, q$ we have $a_i=0$ then it is supposed that
the corresponding $x_i=0$.}

We consider this problem in two versions, either when for all  $i=0,1,q$
such that $a_i\not =0$
the corresponding $x_i$ are variables,
or when $x_0$ and $x_1$ are fixed and $x_q$ vary. We suppose that either
$x_0$ or $x_1$
is nonzero, since otherwise $\max F(x)=\min F(x)=0$.

{\bf Proposition 8.1.} {\em If $a_1\not =0$ then
$$\min F(x) = - {b_1 \over a_1}\;.$$
If $a_1 =0$ and $a_0\not =0$ then
$$\min F(x) = - {b_0 \over a_0}\;.$$
The same problem for fixed $x_0$ and $x_1$ has
$$\min_{x_0, x_1\atop{fixed}} F(x) =- b_0 x_0 - b_1 x_1\;.$$
}

{\em Proof.} The assertion is almost obvious. Indeed, given a vector $x$
with a nonzero
$x_q$ for some $q$, change it, putting $x_q=0$. The value of $F(x)$ then
diminishes, leaving
all the conditions satisfied. Therefore the minimum is attained when
$x_q=0$ for every $q$.
Then we see that the minimum of $- b_0 x_0 - b_1 x_1$ under the conditions
$x_i\ge 0$ and
$a_0 x_0 + a_1 x_1 \le 1$ is attained for $a_0 x_0 + a_1 x_1 = 1$ and 
one of the two
$x_i$ being 0, namely, that with the smaller ratio $b_i/a_i$.${\Box}$

{\bf Proposition 8.2.} {\em Suppose that $x_0$ and $x_1$ are fixed. Then
$$\max_{x_0, x_1\atop{fixed}}  F(x)={( x_0+2 x_1)}
{\left( \sum\limits _ {{p< {{p'}}}}
{{b_{p}}}
+ \alpha{{b_{{p'}}}}\right) - b_0 x_0 - b_1 x_1},$$
where $p'$ and $\alpha\in(0,1]$ are found from the condition
$${{\sum\limits_{{p< {{p'}}}} {{a_{p}}}
+ \alpha{{a_{{p'}}}}}} =
{{1- a_0 x_0 -a_1 x_1}\over {x_0+2 x_1}}\;.$$
}

{\em Proof.} Suppose that $x$
satisfies the requirements (i), (ii), and (iii). Let $x'$ coincide with $x$
in all coordinates
except $x_{p^m}\not =0$ and $x_{p^n}$, $n\le m$, and set
$x'_{p^m}=x_{p^m}-\varepsilon$ and
$x'_{p^n}=x_{p^n}+\varepsilon a_{p^m}/a_{p^n}$. Then, as we pass from $x$
to $x'$,
the left hand side of (iii) does not
change, that of (ii) can only get less because of (1), and $F(x') \ge F(x)$
because of
(2). This proves that we can take

(v)
$$x_{p^m}=0\;\;\hbox{ for } m>1\; , $$
and (ii) is reduced to

(ii$'$)
$$x_{p} \le x_0+2 x_1\;.$$

Now let us deal with

(iii$'$)
$$\sum_{p} a_{p} x_{p} + a_0 x_0 + a_1 x_1 \le 1\;.$$

Suppose again that $x$
satisfies the requirements (v), (i), (ii$'$), and (iii$'$).
Let $x'$ coincide with $x$ in all coordinates
but $x_{{p_1}}\not =0$ and $x_{{p_2}}$, ${{p_1}}\ge {{p_2}}$,
and set $x'_{{p_1}}=x_{{p_1}}-\varepsilon$ and
$x'_{{p_2}}=x_{{p_2}}+\varepsilon {{a_{{p_1}}}\over{a_{{p_2}}}}$.
Again, as we pass from $x$ to $x'$,
the left hand side of (iii$'$) does not
change, and $F(x') \ge F(x)$ because of (v).
Therefore, it is profitable for $F(x)$ to make for
small indices ${{p}}$
the value of $x_{{p}}$ as large as possible, i.e., to set
$$x_{p} =  x_0+2 x_1\;.$$
This we can do, until it starts to contradict (iii$'$).

Summing up, we have proved that there exists a prime $p'$ such that
the maximum of $F(x)$ is attained for some $x$ satisfying the conditions:
$$x_q=0\;\;\hbox{ for }\;\; q \not = p;$$
$$x_{p}=0\;\;\hbox{ for }\;\; {p} > {p'};$$
$$x_{p} =  x_0+2 x_1\;\;\hbox{ for }\;\;p<{p'};$$
$$x_{{p'}} = \alpha ( x_0+2 x_1)$$
for some $\alpha \in( 0 , 1 \rbrack$.
Here $p'$ and $\alpha$ are chosen in such a way that (iii$'$) becomes an
equality.
Then
$$F(x)={( x_0+2 x_1)}
{\left( \sum\limits _ {{p< {{p'}}}}
{{b_{p}}}
+ \alpha{{b_{{p'}}}}\right ) - b_0 x_0 - b_1 x_1}\;,$$
the condition being
$${( x_0+2 x_1)}{{\left(\sum\limits_{{p< {{p'}}}} {{a_{p}}}
+ \alpha{{a_{{p'}}}}\right) +a_0 x_0 +a_1 x_1}} = 1\;.\;{\Box}$$

{\bf Proposition 8.3.} {\em Suppose that $x_0,\; x_1$, and all $x_q$ vary.
Then
$$
\max F(x) =  {{\sum\limits_{{p\le {{p_0}}}} {{b_{p}}} - {{b}}} \over
{\sum\limits_{{p\le {{p_0}}}} {{a_{p}}} + {{a}}}}\;,
$$
where

$a={{a_1}\over 2}$ , $b={{b_1}\over 2}$ if
$a_0=0$ and $ a_1\not =0$,

$a=a_0,\; b=b_0$ if $a_0\not =0$ and $a_1=0$,

and if both
$a_0\not=0,\; a_1\not =0$ we have to compare two possibilities
$a={{a_1}\over 2}$ , $b={{b_1}\over 2}$ and $a=a_0,\; b=b_0$.

Here, for each choice of $a$ and $b$, we let $p'$ run over all primes such
that
$$0\le{{\sum\limits_{{p\le {{p'}}}} b_p - b} \over
{\sum\limits_{{p\le {{p'}}}} a_p + a}}\le {{b_{{p'}}
\over a_{{p'}}}}\;$$
and take
${{p_0}}$ to be the greatest of such ${{p'}}$.
}

{\em Proof.} If $x=(x_0, x_1, x_q)$ is a maximum point (one of) for
our problem, then it is also a maximum point for the problem with $x_0$
and $x_1$ fixed. Therefore, by Proposition 8.2
$$ \max_{x_0, x_1\atop{fixed}} F(x)=( x_0+2 x_1)
{\biggl( \sum\limits _ {p< {p'}}
b_p
+ \alpha b_{p'}\biggr ) - b_0 x_0 - b_1
x_1}\;.$$
Here $p'=p'(y)$ and $\alpha=\alpha(y)$ depend on and are uniquely
determined by
$$
y={{1-a_0 x_0 - a_1 x_1}\over{x_0+2x_1}}\;.
$$
Recall that $y\ge 0$ because of (iii).

Our first goal is to prove that for a fixed $y$ the maximum is attained
when either $x_0=0$, or $x_1=0$. If $a_0=0$ (or $a_1=0$) this follows from
(iv), so we can consider the case  $a_0>0$ (or $a_1>0$).

Indeed,
$$
x_0={{1-(a_1+2y)x_1}\over{y+a_0}}\ge 0
$$
hence $$0\le x_1\le {1\over a_1+2y}\;.$$ Substituting $x_0$ into the
expression for
$$\max_{x_0, x_1 \hbox{\scriptsize{fixed}}} F(x)$$
we see that
for a
fixed $y$ it is linear in $x_1$. Thus the maximum is attained at one of
the ends, i.e., either for $x_1=0$ or for $x_0=0$.

We have to maximize each of them over $y$ which is uniquely determined
by $p'=p'(y)$ and $\alpha=\alpha(y)$, so we can maximize first over
$\alpha\in(0,1]$ and then over $p'$. The expressions being linear in
$\alpha$, the maxima are attained at the end, since $\alpha=0$ and
$\alpha=1$ do not differ up to a change of $p'$.

Therefore, either $x_0=0$ and
$$
\max F(x) =  {{\sum\limits_{{p\le {{p'}}}} {{b_{p}}} - {{b_1}\over 2}} \over
{\sum\limits_{{p\le {{p'}}}} {{a_{p}}} + {{a_1}\over 2}}}
$$
or $x_1=0$ and
$$
\max F(x) =  {{\sum\limits_{{p\le {{p'}}}} {{b_{p}}} - {{b_0}}} \over
{\sum\limits_{{p\le {{p'}}}} {{a_{p}}} + {{a_0}}}}\;.
$$
The last thing to do is to maximize over $p'$. Let $(a,b)$ be either
$(a_0,b_0)$, or $({{a_1}\over 2},{{b_1}\over 2})$. When $p'$ grows, at
some point the expression
$$
\max F(x) =  {\sum\limits_{p\le {p'}} b_p- b \over
\sum\limits_{p\le {p'}} a_p+ a}\;.
$$
becomes positive because of (4). Then we just use the fact that if
$a,\;b,\;A,\;B$ are non-negative and
$$
{B\over A}\ge{b\over a}
$$
then
$$
{b\over a}\le {{B+b}\over {A+a}}\le {B\over A}\;.\;{\Box}
$$

{\bf Remark 8.1.} The same result could, of course, be obtained by
writing out the dual linear problem.

\subsection{Bounds}

{\bf GRH Theorem 8.1} (GRH Bounds). {\em  For any family of number fields
$$   \BS_{\lower}\leq  \liminf_{i\rightarrow \infty}{ \log(h_iR_i) \over g_i}\leq
\limsup_{i\rightarrow \infty}{ \log(h_iR_i) \over g_i}\leq \BS_{\upper},$$
$$0\leq  \liminf_{i\rightarrow \infty}{ \log {\mathbf\kappa}_i \over g_i}\leq
\limsup_{i\rightarrow \infty}{ \log {\mathbf\kappa}_i \over g_i}\leq
{\mathbf\kappa}_{\upper},$$
where $$\BS_{\lower}= 1-{\log 2\pi\over\gamma+\log8\pi} \approx 0.5165...,$$
$$\BS_{\upper}= 1+{{\log{3\over 2}+\log{5\over 4}+\log{7\over 6}}\over{
{\gamma\over 2}+{\pi\over 4}+
\log{2\sqrt{2\pi}}+{\log 2\over {\sqrt 2 -1}}+{\log 3\over {\sqrt 3 -1}}
+{\log 5\over {\sqrt 5 -1}}+{\log 7\over {\sqrt 7 -1}}}} \approx
1.0938\dots,$$
$${\mathbf\kappa}_{\upper}= {\log 2+{\log{3\over 2}}\over{
{\gamma\over 2}+\log{2\sqrt{2\pi}}+{\log 2\over
{\sqrt 2 -1}}+{\log 3\over {\sqrt 3 -1}}}}\approx
0.2164\dots,.$$
Moreover, if all the fields in the family are totally real then
$$
\liminf_{i\rightarrow \infty}{ \log(h_iR_i) \over g_i}\ge \BS_{{\mathbb{R}},\lower}
$$
where
$$
\BS_{{\mathbb{R}},\lower}=1-{\log 2\over{{\gamma\over 2}+{\pi\over
4}+\log{2\sqrt{2\pi}}}} \approx 0.7419\dots,$$
and 
$$\liminf_{i\rightarrow \infty}{ \log({\mathbf\kappa}_i) \over g_i}\le
{\mathbf\kappa}_{{\mathbb{R}},{\upper}}
$$
where
$$
{\mathbf\kappa}_{{\mathbb{R}},{\upper}}=
{\log 2+{\log{3\over 2}}\over{
{\gamma\over 2}+\log{2\sqrt{2\pi}}+{\pi\over 4}+{\log 2\over
{\sqrt 2 -1}}+{\log 3\over {\sqrt 3 -1}}}} \approx 0.1874\dots,$$
If all the fields are totally complex then
$$
\limsup_{i\rightarrow \infty}{ \log(h_i R_i) \over g_i}\le
\BS_{{\mathbb{C}},\upper},
$$
where
$$
\BS_{{\mathbb{C}},\upper}=1+
{{\sum\limits_{{p=2}\atop{prime}}^{13} \log{p\over
{p-1}}-{1\over 2}\log{2\pi}}\over{{\gamma\over
2}+\log{2\sqrt{2\pi}}+\sum\limits_{{p=2}\atop{prime}}^{13}
{\log p \over {\sqrt p -1}}}} \approx 1.0764\dots\;.
$$
}

{\em Proof.} Since any family contains an asymptotically exact one
(Lemma 2.2), any limit point of the ratio is a limit for some
asymptotically exact family, and it is enough to prove the theorem for
such families. The Generalized Brauer--Siegel Theorem (GRH Theorem 7.2) gives
us the limit value of the ratio in terms of $\phi=\{\phi_\alpha\}$. 
The Basic Inequality (GRH Theorem 3.1)
gives us a
restriction. Up to a constant 1 we get an optimization problem of the
type described above with
$$
\begin{array}{ll}
b_0=\log 2\approx 0.693..., &
a_0=\log 2\sqrt{2\pi}+{\pi\over 4}+{\gamma\over 2}\approx 2.686..., \cr\cr
b_1=\log {2\pi}\approx 1.837..., \;\;\;\;\; &
a_1=\log {8\pi}+\gamma\approx 3.801...,\cr\cr
b_q=\log{q\over{q-1}}, &
a_q={{\log q}\over{\sqrt q -1}}.
\end{array}
$$
We have to check the conditions (1)---(4) above (see the beginning of 
Subsection 8.1). For (1) we see that
$${m\over n} \ge{{m(p^{n/2}-1)}\over {n(p^{m/2}-1)}} = {a_{p^m}\over
a_{p^n}}\;\;\hbox{ for }\; n\le m.$$
To check (2) it is enough to prove that
$$
f(x)={{(\sqrt x -1)\log{x\over{x-1}}}\over{\log x}}
$$
is decreasing for $x\ge 2$. This is quite straightforward.

As for (3), it is obvious because of the numerical values given.

To prove (4), just note that
$$
\sum_q a_q\ge \sum_q b_q\ge -\sum_{p\atop {prime}}\log(1-{1\over
p})=\log\zeta(1)=\infty.
$$

So we come to the above optimization problem,
where $1+F(x)$ is the right hand side of the Generalized Brauer--Siegel,
(i) corresponds to non-negativity of $\phi_q$, $\phi_{\mathbb{R}}$ and
 $\phi_{\mathbb{C}}$, (ii) is the condition of Lemma 2.4, (iii) is the GRH
Basic Inequality, and (iv) is empty since all $a_i\not=0$.

We can now use Proposition 8.1. If $a_1\not=0$ we get
$$
\min F(x)=-{b_1\over a_1}=-{{\log{2\pi}}\over{\log{8\pi}+\gamma}}
\approx -0.4834...,
$$
which gives the value of $\BS_{\lower}=1+\min F(x)$. If all the fields are totally
real, i.e., $a_1=0$, Proposition 8.1 gives
$$
\min F(x)=-{b_0\over
a_0}=-{{\log{2}}\over{{\gamma\over 2}+{\pi\over
4}+\log{2\sqrt{2\pi}}}}\approx -0.2580...
$$
As for the maxima, Proposition 8.3 gives two possibilities, either
$$
\max F(x)=C_{p'}^0=
{{\sum\limits_{{p\le {{p'}}}} {{b_{p}}} - {{b_0}}} \over
{\sum\limits_{{p\le {{p'}}}} {{a_{p}}} + {{a_0}}}}\;,
$$
or
$$
\max F(x)=C_{p'}^1=
{{\sum\limits_{{p\le {{p'}}}} {{b_{p}}} - {{b_1}\over 2}} \over
{\sum\limits_{{p\le {{p'}}}} {{a_{p}}} + {{a_1}\over 2}}}\;,
$$
and in both cases we still have to find out $p'$. Note that
$b_2=b_0=\log 2$.  The values of $C_{p'}^0$ are easily computable, and
we have
$$
0=C_2^0<C_3^0<C_5^0<C_7^0\approx 0.0938...
$$
and
$$
C_7^0>{b_{11}\over a_{11}}\approx 0.092...
$$
Doing the same for $C^1_{p'}$'s we get
$$
C_2^1<C_3^1<C_5^1<C_7^1<C_{11}^1<C_{13}^1\approx 0.0764...
$$
and
$$
{b_{17}\over a_{17}}\approx 0.066...
$$
Therefore, by Proposition 8.3, $\BS_{\upper}=C_7^0$ and
$\BS_{{\mathbb{C}},\upper}=C^1_{13}$. 

The bounds for $\mathbf\kappa$ are obtained in the same way, 
and we leave details to
the reader.
${\Box}$

{\bf Theorem 8.2} (Unconditional Upper Bound). {\em  For any family of number
fields
$$\limsup_{i\rightarrow \infty}{ \log(h_iR_i) \over g_i}\leq
1+{{\sum\limits_{{p=3}\atop{prime}}^{23} \log{p\over {p-1}}}\over{
{\gamma\over 2}+{1\over 2}+\log{2\sqrt{\pi}}+
2\sum\limits_{{p=2}\atop{prime}}^{23} \log p \sum\limits_{m=1}^\infty
{1\over {p^m+1}}}}
 \approx 1.1588\dots,$$
$$ \limsup_{i\rightarrow \infty}{ \log {\mathbf\kappa}_i \over g_i}\leq
 {\mathbf\kappa}_{\unc, \upper}=
1+{{\sum\limits_{{p=3}\atop{prime}}^{5} \log{p\over {p-1}}}\over{
{\gamma\over 2}+\log{2\sqrt{\pi}}+
2\sum\limits_{{p=2}\atop{prime}}^{5} \log p
\sum\limits_{m=1}^\infty {1\over {p^m+1}}}}
 \approx 0.3151\dots
 $$
Moreover, if all the fields are totally complex then 
$$\limsup_{i\rightarrow \infty}{
\log(h_iR_i) \over g_i}\leq 1+{\sum\limits_{{p=2}\atop{prime}}^{179} \log{p\over
{p-1}}-\log\sqrt{2\pi}
\over {\gamma\over 2}+\log{2\sqrt{\pi}}+2\sum\limits_{{p=2}\atop{prime}}^{179}
 \log p \sum\limits_{m=1}^\infty {1\over p^m+1}
   } \approx 1.0965\dots,$$
and if all the fields are totally real then 
$$ \limsup_{i\rightarrow \infty}{ \log {\mathbf\kappa}_i \over g_i}\leq
1+{{\sum\limits_{{p=3}\atop{prime}}^{5} \log{p\over {p-1}}}\over{
{\gamma\over 2}+\log{2\sqrt{\pi}}+{1\over2}+
2\sum\limits_{{p=2}\atop{prime}}^{5} \log p
\sum\limits_{m=1}^\infty {1\over {p^m+1}}}}
 \approx 0.2816\dots.$$}

{\em Proof.} Along the same lines as above. This time we use  the
unconditional Basic Inequality$'$ (Proposition 3.1) and the Generalized
Brauer--Siegel Inequality (Theorem 7.1). We get a maximization problem
of the same type with
$$
\begin{array}{ll}
b_0=\log 2\approx 0.693..., &
a'_0=\log 2\sqrt{\pi}+{1\over2}+{\gamma\over 2}\approx 2.054..., \cr\cr
b_1=\log {2\pi}\approx 1.837...,\;\;\;\;\; &
a'_1=\log {4\pi}+\gamma\approx 3.108...,\cr\cr
b_q=\log{q\over{q-1}}\;, &
a'_q=2\log q\sum\limits_{m=1}^\infty(q^m+1)^{-1}\;.
\end{array}
$$
Again, we have to check the conditions (1)---(4). This is done as in the previous proof.
The only tedious point to check is (2), and again it is enough to show that
$$
f(x)={\log{x\over{x-1}}\over{2\log
x\sum\limits_{m=1}^\infty(x^m+1)^{-1}}}
$$
is decreasing, which is again straightforward.

Using, as above, Proposition 8.3 we have to maximize
$$
c^0_{p'}=
{{\sum\limits_{{p\le {{p'}}}} {{b_{p}}} - {{b_0}}} \over
{\sum\limits_{{p\le {{p'}}}} {{a'_{p}}} + {{a'_0}}}}
$$
and
$$
c^1_{p'}=
{{\sum\limits_{{p\le {{p'}}}} {{b_{p}}} - {{b_1}\over 2}} \over
{\sum\limits_{{p\le {{p'}}}} {{a'_{p}}} + {{a'_1}\over 2}}}\;.
$$
We get
$$
0=c^0_2<c^0_3<c^0_5<\ldots <c^0_{23}\approx 0.1588...,
$$

$$
{b_{29}\over a'_{29}}\approx 0.150...\;\;.
$$
We also have
$$
c^1_2<c^1_3<c^1_5<\ldots <c^1_{179}\approx 0.0965...,
$$

$$
{b_{181}\over a'_{181}}\approx 0.0964...
$$
Again, we leave $\mathbf\kappa$ to the reader.${\Box}$

{\bf Remark 8.2.} In \cite{Hof} it is proved that
$$ \mathbf\kappa({\mathcal K})\le 0.958-1.936\phi_{\mathbb{R}}-2.936\phi_{\mathbb{C}}.$$
Using our estimates one easily gets
$$ \mathbf\kappa({\mathcal K})\le 0.946-1.936\phi_{\mathbb{R}}-2.936\phi_{\mathbb{C}},$$
and also
$$ \mathbf\kappa({\mathcal K})\le 0.654-1.343\phi_{\mathbb{R}}-2.032\phi_{\mathbb{C}},$$
which is always better (for $\phi_{\mathbb{R}}$ and $\phi_{\mathbb{C}}$ allowed by the
Odlyzko bound). Still better estimates of $ \mathbf\kappa({\mathcal K})$ in terms of
$\phi_{\mathbb{R}}$ and $\phi_{\mathbb{C}}$ follow from Proposition 8.2. 
(Note, however, that in \cite{Hof} the result is not just asymptotic,
 but effective.) 

On the other side we get

{\bf Theorem 8.3} (Unconditional Lower Bound). {\em If for a tower of
almost normal number fields there exists $\alpha>0$ such that $\inf {N_{i}\over
{g_i}}\ge\alpha$, then
$$\liminf_{i\rightarrow \infty}{ \log(h_iR_i) \over g_i}\geq
1-{{\log{2\pi}}\over{\gamma+\log{4\pi}}}\approx 0.4087...$$
If, in addition, the fields are totally real, then
$$\liminf_{i\rightarrow \infty}{ \log(h_iR_i) \over g_i}\geq
1-{2{\log{2}}\over{\gamma+1+\log{4\pi}}}\approx 0.6625...$$
}

{\em Sketch of proof.} The same argument as in the proof of GRH Theorem 8.1 is
applied. We use the unconditional
Generalized Brauer--Siegel Theorem (Theorem 7.3) and Proposition 3.1. ${\Box}$

\section{Class field towers}

\subsection{Infinite Unramified Towers with Splitting Conditions}

We are going to present some examples. The main goal of this section
is to show that the Brauer--Siegel
ratio does not necessarily tend to 1. Because of the Brauer--Siegel theorem,
we need families  of fields for which $[K:{\mathbb{Q}}]/ g_K$ does not tend
to 0.
One's thought turns immediately to unramified towers, and the only 
infinite examples we know are Hilbert class field towers satisfying
some extra conditions.

Recall that for a field $K$ its Hilbert class field $K_{\Hilb}$ is defined
as the maximal unramified abelian extension, and $\Gal(K_{\Hilb}/K)= \Cl_K$.
We fix a prime $\ell$
and consider the maximal unramified abelian $\ell$-extension
$K_{\Hilb,\ell}$ with
 $\Gal(K_{\Hilb, \ell}/K)= \Cl_{\ell,K}$, where $\Cl_{\ell,K}$ is the Sylow
$\ell$-subgroup
of $\Cl_K$.
Put $K_0=K, K_1=K_{\Hilb, \ell}, K_2=(K_1)_{\Hilb, \ell}$, etc.
There are two
possibilities,
either $K_n=K_{n+1}=\ldots$ for some $n$, or all these fields are
different, i.e., the tower
$K=K_0\subset
K_1\subset K_2\subset \ldots$ is infinite. The latter is the situation we
are looking for.

Note that if $K_0$ is totally real (respectively, totally complex), such are all
fields of the tower.

For a group $A$, let $d_{\ell}(A)=\dim_{{\mathbb{F}}_\ell}(A^{\ab}/\ell)$
denote its $\ell$-rank.
For a field $F$ let $r_1(F)$ be the number of its real, and $r_2(F)$ the
number of its complex
places.
  
Consider a degree $\ell$ extension of number fields $K/k$. Set
$r_1=r_1(k)$, $r_2=r_2(k)$,
let 
$r$ be the number of prime ideals of $k$ ramified in $K/k$ and 
$\rho$ be the number of real places  ramified in $K/k$ (i.e., becoming
complex),
$\delta_\ell=d_\ell(W_K)$, $W_K$ being the group of roots of 1 lying in $K$,
i.e., $\delta_\ell=0$ if there is no $\ell$-root of 1
in $k$, and $\delta_\ell=1$ otherwise. In  \cite{15} the following statement is
proved.

{\bf Proposition 9.1} (J.Martinet). {\em In the above notation, if
$$r\ge r_1+r_2+\delta_\ell+2-\rho+2\sqrt{\ell(r_1+r_2-\rho/2)+\delta_\ell}$$
then $K$ has an infinite unramified Hilbert class field $\ell$-tower
$K=K_0\subset K_1\subset K_2\subset \ldots$.${\Box}$}

Here are the best specimens found in the hunt for small discriminants,
obtained with the help of
Proposition 9.1.

{\bf Corollary 9.1} (J.Martinet). {\em The fields
$${\mathbb{Q}}(\sqrt{3\dot 5\dot 13\dot 29\dot 61}),$$
$${\mathbb{Q}}(\sqrt{2},\sqrt{3\dot 5\dot 7\dot 23\dot 29}),$$
$${\mathbb{Q}}(\sqrt{-3\dot 5\dot 17\dot 19}),$$
$${\mathbb{Q}}(\cos {2\pi\over 11},\sqrt{2},\sqrt{-23})$$
have infinite unramified class field $2$-towers.}

{\em Proof.} This is a non-obvious corollary of Proposition 9.1 which is
proved in \cite{15}, Examples 3.2, 4.2, 5.3 and 6.2.${\Box}$

We shall also need unramified towers with some extra splitting
conditions. To study them, 
let us fix some notation. Let $C_K=J_K/K^*$ be the id\`ele class group. Let
$S_\infty$ be the set
of archimedean places of $K$.
Fix a finite set $S$ of prime ideals of $K$. In what follows if $S$ is empty
we omit $S$ in the corresponding notation. Let
$O_{K,S}$ be the ring of $S$-integers,
$U_{K,S}=\prod\limits_{v\in S\cup S_\infty} \!\!K_v^*
\prod\limits_{v\not\in S\cup S_\infty}\!\! O^*_{K_v}$ the group of
id\`ele $S$-units,
$E_{K,S}=O_{K,S}^*=U_{K,S}\cap K^*$ the group of $S$-units in $K$,
$I_{K,S}$ the group of fractional ideals nondivisible by prime ideals
of $S$,
$I_K^{S}$ the group of ideals generated by
prime ideals in $S$ (of course, $I_K=I_{K,S}\oplus I_K^{S}$),
$\;P_{K,S}$ the image of the principle ideal group $P_K$ in $I_{K,S} $
under $P_K\subset I_K
\rightarrow I_{K,S}$, $I_K
\rightarrow I_{K,S}$ being the natural projection,
$\Cl_{K,S}=I_{K,S}/P_{K,S}=\Cl_{K}/\Im(I_K^{S})$ the group of $S$-classes,
$W_K$ the group of roots of unity lying in $K$.

{\bf Theorem 9.1.} {\em Let $P=\{ p_1,\ldots,p_t\}$ and
$Q=\{ q_1,\ldots,q_r\}$
be disjoint sets of prime ideals of $k$, and let $t_0$ be the number of
principal ideals in $P$.
 Consider a number field
$K/k$ of prime degree $\ell$, ramified exactly at $Q$. Let $S$ be the set of
prime
ideals in $K$ lying over $P$, $s=\vert S\vert$.
If $$r\ge s-t_0+
r_1+r_2+\delta_\ell+2-\rho+2\sqrt{\ell(r_1+r_2-\rho/2)+\delta_\ell+ s}$$
then the field $K$ has an infinite unramified class field $\ell$-tower
$K=K_0\subset
K_1\subset K_2\subset \ldots$,
where $S$ splits completely.}

To prove this theorem we need some lemmata.

In fact, all constructions of unramified towers we know are based on
the following well-known lemma
(cf. \cite{20}, \cite{14}).

{\bf Lemma 9.1.} {\em If
$$d_\ell(\Cl_{K,S})\ge2+2\sqrt{d_\ell(E_{K,S})+1}$$
then $K$ has an infinite class field $\ell$-tower $K=K_0\subset
K_1\subset K_2\subset \ldots$,
where $S$ splits completely.}

{\em Proof.} Let $L$ be the union of all $K_i$, where $K_0=K$ and
$K_i/K_{i-1}$
is the abelian
$\ell$-extension corresponding by the class field theory to the
$\ell$-Sylow subgroup of
$\Cl_{K_{i-1},S_{i-1}}$, where $S_{i-1}$ consists of all places lying over $S$.
The places of $S$ split completely in $L$.

Suppose that $L /K$ is of finite degree, ${\cal G}=Gal(L /K)$,
let $S_L$
be the set of places of $L$ lying over $S$. As for any
finite $\ell$-group, we have (see \cite{20}, eq.6 )
$${1\over 4}(d_{\ell}({\cal G}))^2 - d_{\ell}({\cal G}) <
d_{\ell}(H_2({\cal G},{\mathbb{Z}})).$$
Then
$$H_2({\cal G},{\mathbb{Z}})=H^{-3}({\cal G},{\mathbb{Z}})=H^{-1}({\cal
G},C_{L})$$
by Tate's fundamental theorem (see \cite{28}, section 11.3). We have
$$0\;\longrightarrow\;U_{L,S_L}/E_{L,S_L}\;\longrightarrow\;
C_{L}\;\longrightarrow\;
\Cl_{L,S_L}\;\longrightarrow\;0$$
and the $\ell$-Sylow subgroup of $\Cl_{L,S_L}$ is trivial, otherwise $L$
would yet have another nontrivial $\ell$-extension splitting $S_L$. Hence
$$d_{\ell}(H^{-1}({\cal G},C_{L}))=
d_{\ell}(H^{-1}({\cal G},U_{L,S_L}/E_{L,S_L})).$$
The extension $L/K$ being unramified, $U_{L,S_L}$ is cohomologically
trivial,
since first each local component $\prod\limits_{w\not\in S_L\cup S_{L_\infty}}
O^*_{L_w}$ is
cohomologically
trivial, and next the points of $S$ split and hence $\prod\limits_{w\in S_L\cup
S_{L_\infty}}L_w^*$ is cohomologically trivial and by definition
$$U_{L,S}=\prod\limits_{w\in S_L\cup S_{L_\infty}} L_w^*
\prod\limits_{w\not\in S_L\cup S_{L_\infty}} O^*_{L_w}.$$
Hence
$$H^{-1}({\cal G},U_{L,S_L}/E_{L,S_L})=
{\hat H}^0({\cal G},E_{L,S_L})=
E_{K,S}/N_{L/K}(E_{L, S_L}).$$
The $\ell$-rank of the latter being less than or equal to
$d_{\ell}(E_{K,S})$ we see that
$${1\over 4}(d_{\ell}({\cal G}))^2 - d_{\ell}({\cal G}) < d_{\ell}(E_{K,S}).$$
On the other hand,
$d_{\ell}({\cal G})=d_{\ell}({\cal
G}^{\ab})=d_{\ell}(Gal(K_1/K))=d_{\ell}(\Cl_{K,S})$
and we get
$${1\over 4}(d_{\ell}(\Cl_{K,S}))^2 - d_{\ell}(\Cl_{K,S}) <
d_{\ell}(E_{K,S}),$$
i.e.,
$$d_{\ell}(\Cl_{K,S})<2+2\sqrt{d_{\ell}(E_{K,S})+1}.{\Box}$$

{\bf Lemma 9.2.} {\em We have
$$d_{\ell}(E_{K, S})= r_1(K)+r_2(K)+\delta_\ell(K)-1+ s$$
$$=\ell(r_1+r_2-\rho/2)+\delta_\ell-1+ s.$$
}

{\em Proof.} It is well known that $E_{K,S}=W_K\oplus
{{\mathbb{Z}}^{r_1(K)+r_2(K)+s-1}}$
(cf.\cite{13}, V.1). In our case, $r_1(K)=\ell (r_1-\rho) $, $r_2(K)=\ell
r_2+\ell\rho/2$,
and  $\delta_\ell(K)=\delta_\ell$ since no new $\ell$-root can appear in an
$\ell$-extension.
${\Box}$

{\bf Lemma 9.3} (J.Martinet). {\em We have
$$d_\ell(\Cl_{K})\ge
r-r_1-r_2+\rho-\delta_\ell .$$}

{\em Proof.} This is proved in \cite{15}, section 2.$\;{\Box}$

{\bf Lemma 9.4.} {\em We have
$$d_\ell(\Cl_{K,S})\ge d_\ell(\Cl_{K})-s+t_0 .$$}

{\em Proof.} Let $P_0$ be the set of principal ideals lying in $P$. Let
$\phi:I^S_K\rightarrow\Cl_K$ be the composition of natural maps
$I^S_K\rightarrow I_K$ and
$I_K\rightarrow \Cl_K$. By definition $\Cl_{K,S}=\Cl_{K}/\Im\phi$. We have
$I^S_K=\prod\limits_{p\in P}(\prod\limits_{w\vert
p}w^{\mathbb{Z}})\simeq{\mathbb{Z}}^s$. Look at the kernel of
$\phi$. Since for any $p\in P_0$ there is the relation
$\prod\limits_{w\vert p}w\in \Ker\phi$, we get
$\rk_{\mathbb{Z}}\Ker\phi\ge t_0$. Therefore, $\rk_{\mathbb{Z}}\Im\phi\le
s-t_0$. It remains to
take $d_\ell$. ${\Box}$

{\em Proof of Theorem} 9.1. By Lemmata 9.2, 9.3, 9.4 and the inequality of
the theorem we get
$$d_\ell(\Cl_{K,S})\ge d_\ell(\Cl_{K})-s+t_0$$
$$ \ge r-r_1-r_2+\rho-\delta_\ell-s+t_0 $$
$$\ge s-t_0+
r_1+r_2+\delta_\ell+2-\rho+2\sqrt{\ell(r_1+r_2-\rho/2)+\delta_\ell+ s}
-r_1-r_2+\rho-\delta_\ell-s+t_0$$
$$=2+2\sqrt{\ell(r_1+r_2-\rho/2)+\delta_\ell+ s}$$
$$=2+2\sqrt{d_\ell(E_{K,S})+1}.$$
By Lemma 9.1 this proves the theorem.$\;{\Box}$

{\bf Corollary 9.2.} {\em Let $P=\{ p_1,\ldots,p_t\}$ and
$Q=\{ q_1,\ldots,q_r\}$
be disjoint sets of primes. Consider a quadratic number field
$K/{\mathbb{Q}}$ ramified exactly at $Q$. Let $\sigma$ be the number of primes in
$P$ that split in $K$, and $s=t+\sigma$ the total number of prime
ideals in $K$ lying over $P$. Suppose that either $K$ is complex quadratic and
$$r\ge 3+\sigma+2\sqrt{2+s},$$
or $K$ is real quadratic and
$$r\ge 4+\sigma+2\sqrt{3+s}
.$$
Then $K$ has  an infinite unramified class field $2$-tower totally splitting
all prime ideals over $P$.}

{\em Proof.} Indeed, here $k={\mathbb{Q}}$, $\ell=2$, $\delta_2=1$,
$r_1=1$, $r_2=0$,
$\rho=1$ for the complex quadratic case and 0 for the real quadratic one,
$t_0=t$.$\;{\Box}$

Here are some numerical examples.

{\bf Corollary 9.3.}{\em The field
$${\mathbb{Q}}(\sqrt{11\dot 13\dot 17\dot 19\dot 23\dot 29\dot 31\dot 37\dot 41
\dot 43\dot 47\dot 53\dot 59\dot 61\dot 67})$$
has  an infinite unramified class field $2$-tower
in which nine prime ideals lying over $2$, $3$, $5$, $7$ and $71$ split
completely.}

{\em Proof.} A straightforward check shows that 2, 3, 5, 7 split in
$K/{\mathbb{Q}}$ and
71 is inert. Then we apply Corollary 9.2.$\;{\Box}$

{\bf Corollary 9.4.}{\em The field
$${\mathbb{Q}}(\sqrt{-13\dot 17\dot 19\dot 23\dot 29\dot 31\dot 37\dot 41
\dot 43\dot 47\dot 53\dot 59\dot 61\dot 73\dot 79})$$
has  an infinite unramified class field $2$-tower
in which ten prime ideals lying over $2$, $3$, $5$, $7$ and $11$ split
completely.}

{\em Proof.} Along the same lines.$\;{\Box}$

We shall exploit these examples in Subsection 9.3. We also get

{\bf Corollary 9.5} (Y.Ihara \cite{9}). {\em The field
$${\mathbb{Q}}(\sqrt{-3\dot 5\dot 7\dot 11\dot 13\dot 17\dot 23\dot 31})$$
has an infinite unramified class field $2$-tower
in which two prime ideals lying over $2$ split
completely.}$\;{\Box}$

\subsection{A remark on the deficiency problem}

In \cite{33} K.Yamamura writes
 
`` Combining Ihara's remark (\cite{9}, sect.14) to Golod-Shafarevich theory
(cf.\cite{20}) and Martinet's
result (\cite{15}), we easily obtain the following

Theorem. {\em Let $K/k$ be a cyclic extension of degree $p$ ($p$: a prime
number) of an algebraic
number field of finite degree. Let {\mathg S} be a given set of finite
primes of $K$. Let $r'$ be
the number of those finite primes of $k$ which are ramified in $K$ and none
of its extension to $K$
belongs to {\mathg S}. If
$$r'\ge r_1+r_2+\delta^{(p)}_k
+2-\rho+2\sqrt{H+p(r_1+r_2-\rho/2)+\delta^{(p)}_k},$$
then $K$ has an infinite {\mathg S}-decomposing $p$-class field tower. Here
$\rho$ denotes the
number of real primes of $k$ which are ramified in $K$, $r_1=r_1(k)$,
$r_2=r_2(k)$,
and $H=\vert${\mathg S}$\vert$.} ''

In our notation, $r'=r$, $p=\ell$, $\delta^{(p)}_k=\delta_\ell$, $H=s$,
{\mathg S}$=S$, and the
inequality reads
$$r\ge r_1+r_2+\delta_\ell+2-\rho+2\sqrt{\ell(r_1+r_2-\rho/2)+\delta_\ell+s}.$$
 
We should admit that we consider the Yamamura theorem to be not only
unproved, but most likely
false. To explain this point of view let us prove the following

{\bf Proposition 9.3.} {\em If the Yamamura theorem is true, the
generalized Riemann hypothesis is
false.}

{\em Proof.} Let $\ell=2$, $k={\mathbb{Q}}$. Consider the field
$$K={\mathbb{Q}}(\sqrt{-13\dot 17\dot 19\dot 23\dot 29\dot 31\dot 37\dot
41\dot 61\dot 101}),$$
We have $g_K\approx 17.16493$, $\rho=r_1=0$, $r_2=1$, $r=10$, $\delta_2=1$.
Let $S$ consist of ten
ideals lying over 2, 3, 5, 7, 11 (straightforward calculation shows that
these primes split in
$K/{\mathbb{Q}}\;$). Then the Yamamura theorem gives the infiniteness of
the unramified class field 2-tower over $K$ since $10>3+2\sqrt{2+10}$.
We have
$$  \sum_{q}{\phi_q\log q\over {\sqrt
q-1}}+\alpha_{\mathbb{R}}\phi_{\mathbb{R}}
+\alpha_{\mathbb{C}}\phi_{\mathbb{C}}=$$
$${1\over g}\left(\gamma+\log {8\pi}+{{2\log 2}\over{\sqrt 2 -1}}+
{{2\log 3}\over{\sqrt 3 -1}}+{{2\log 5}\over{\sqrt 5 -1}}+
{{2\log 7}\over{\sqrt 7 -1}}+{{2\log {11}}\over{\sqrt {11}
-1}}\right)\approx1.0013... >1\;.$$
This contradicts GRH Basic Inequality (cf. GRH Theorem 3.1  or \cite{9},
2-2).${\Box}$

Unfortunately enough, the rest of \cite{33} is derived from the above Yamamura theorem
and we have to discard all his examples.

In particular, the smallest known deficiency $\delta$ is that of Hajir and Maire \cite{6}
with $\delta\le 0.141\dots$ (cf. the end of Section 3.1). It has $S=\emptyset$. Ihara's example
$${\mathbb{Q}}(\sqrt{-3\dot 5\dot 7\dot 11\dot 13\dot 17\dot 23\dot 31}),$$
$S$ consisting of two divisors of 2, has $\delta\le 0.248...$.

\subsection{Examples}

{\bf GRH Theorem 9.2.} {\em Consider the Martinet field
$$K={\mathbb{Q}}(\cos{2\pi\over11},\sqrt{2},\sqrt{-23})$$
of degree $20$ over
${\mathbb{Q}}$. We have
$$
D_K=2^{30}11^{16}23^{10}\;,
$$
$$
g=g(K)=\log\sqrt{\vert D_K\vert}\approx 45.2578...
$$
This field has an infinite unramified $2$-tower ${\mathcal K}$, and
we have 
$$1-{{10\log(2\pi)}\over g}=\BS_{\lower}({\mathcal K})\le
\BS({\mathcal K})\le \BS_{\upper}({\mathcal K}),$$
$$0\le
{\mathbf\kappa}({\mathcal K})\le {\mathbf\kappa}_{\upper}({\mathcal K}),$$
where
$$\BS_{\upper}({\mathcal K})=\BS_{\lower}({\mathcal K})+{{(\sqrt{23}-1)\log{23\over 22}}\over
\log{23}}(1-{{10(\gamma+\log{8\pi})}\over{g}}),$$
$${\mathbf\kappa}_{\upper}({\mathcal K})=
{{(\sqrt{23}-1)\log{23\over 22}}\over
\log{23}}(1-{{10(\gamma+\log{8\pi})}\over{g}}),$$
i.e., approximately, 
$$0.5939\dots\le \BS({\mathcal K})\le 0.6025\dots,$$
$$0\le{\mathbf\kappa}({\mathcal K})\le 0.0086\dots.$$
The deficiency $\delta({\mathcal K})$ of this tower is at most
$$
1-{{10(\gamma+\log{8\pi})}\over{g}}\approx
0.1601\dots
$$
}

{\em Proof.} Recall first (Corollary 9.1) that this field
has an infinite unramified tower.

Let $K_0={\mathbb{Q}}(\cos{2\pi\over 11})$, $K_{11}={\mathbb{Q}}({\root 11 \of
1})$, $k={\mathbb{Q}}(\cos{2\pi\over 11}, \sqrt 2)$, $F_{23}={\bf
Q}(\sqrt{-23})$, $F_2={\mathbb{Q}}(\sqrt 2)$. The discriminant of a
cyclotomic field is well-known (cf., \cite{13}, IV.1), so we have
$$
D_{K_{11}}=11^9.
$$
Hence, $D_{K_{11}}/{\mathbb{Q}}$ and $D_{K_{0}}/{\mathbb{Q}}$ are unramified
outside of 11. Since 11 is totally ramified in $D_{K_{11}/{\mathbb{Q}}}$,
it is also totally
ramified in $K_0/{\mathbb{Q}}$, therefore, $D_{K_0}=11^4$. The field $K$
is the composite of $K_0$, $F_2$ and $F_{23}$. We get
$$
D_K=D_{K_0}^4 D_{F_2}^{10} D_{F_{23}}^{10}=11^{16}2^{30}23^{10},
$$
and derive the above value of $g$.

The deficiency $\delta$ for the tower is at most
$$
1-{{10(\gamma+\log{8\pi})}\over{g}}\approx
0.1601...
$$

Let us first prove that in $K/{\mathbb{Q}}$ we have the following decomposition
of small primes
$$
\mbox{
\begin{tabular} {c|ccccccccc}
$v$&2&3&5&7&11&13&17&19&23\\
\hline
$e_v$&2&1&1&1&5&1&1&1&2\\
$f_v$&5&10&20&10&4&10&10&20&1\\
$n_v$&2&2&1&2&1&2&2&1&10
\end{tabular}
}
$$
where $e_v$ is the ramification index, $f_v$ the inertia one, and $n_v$ is
the number
of places over $v$.

In $K_{11}/{\mathbb{Q}}$ only 11 is
ramified (totally) and
for a prime $p$ the inertia index $f_p$ equals the smallest $f$ such
that $p^f\equiv 1 (\mod 11)$.
We easily check that $f_3=f_5=5$,
$f_2=f_7=f_{13}=f_{17}=f_{19}=10$ and $f_{23}=1$.
Since 2 does not divide 5, $K_0$ being index 2 subfield of $K_{11}$,
we see that
in $K_0/{\mathbb{Q}}$ all the primes of our table except 11
and 23 are inert, i.e., $f_p=5$, that 11 is totally ramified and 23 is
totally split.
In $F_2/{\mathbb{Q}}$ only 2 is ramified, and $p$ is split if and only if 2
is a square
modulo $p$, i.e., if and only if $p\equiv \pm 1 (\mod 8)$; such
primes are 7, 17 and 23, and the rest (3, 5, 11, 13, 19) are inert. In
$F_{23}/{\mathbb{Q}}$ the only ramified prime is 23 since $-23\equiv 1 (\mod
4)$, and the splitting condition is for $-23$ to be a
square $\mod p$. Thus 2, 3, 13 are split, and 5, 7, 11, 17 and 19 are
inert.

Summing up this information we get the above table. Note also that in
$k/{\mathbb{Q}}$ there are 10 prime ideals over 23, and they all ramify in
$K/k$.

Looking at the decomposition table above, we see that the smallest norm
for which there exists a prime ideal of $K$ is 23 (indeed, $2^5>23$,
etc.)

The next thing to do is to apply the linear programming approach of
Subsection 8.1 to get the minimum and maximum of 
$\BS({\mathcal K})$. We set $a_0=0$, as well as $a_q=0$ for
$q<23$, recalling that if $a_i=0$ we have also $b_i=0$ and we do not optimize over $x_i$.
The right-hand side of GRH Theorem 7.2 becomes
$$
1-{10\over g}\log{2\pi}+F(x)=\BS_{\lower}({\mathcal K})+F(x),
$$
where
$$
F(x)=\sum_{q\ge 23} b_qx_q, b_q=\log{q\over{q-1}}, x_q=\phi_q
\hbox { for } q\ge 23.
$$
The restrictions are, as usual, (i) $x_q\ge 0$, (ii)
$\sum\limits_{m=1}^\infty mx_{p^m}\le{20 \over g} $ for any $p$,  and (iii)
${10\over g}(\gamma+\log{8\pi})+\sum\limits_{q\ge 23}
a_q x_q\le 1$, where $a_q={{\log q}\over{\sqrt q -1}}$.

The minimum of $F(x)$ is clearly 0, and it is easy to check that the
maximum is attained for all $x_q=0$ except for
$$
x_{23}={{\sqrt{23}-1}\over{\log{23}}}(1-{10\over
g}(\gamma+\log{8\pi})).
$$
Indeed, $x_{23}\approx 0.19...<{20\over g}$ which checks (ii),
and $x_{23}$ is chosen so that (iii) becomes an equality.

As for $\mathcal\kappa({\mathcal K})$, we have
$$\mathcal\kappa({\mathcal K})=\BS({\mathcal K})-1+\phi_{\mathbb
C}\log(2\pi).{\Box}$$

{\bf Remark 9.1} In all our examples, once we have some information
on $\BS({\mathcal K})$, we also have it on $\mathcal\kappa({\mathcal K})$,
the difference between the two being known. That is why, in many cases,
we do not say a word about $\mathcal\kappa({\mathcal K})$.

{\bf GRH Theorem 9.3.} {\em Consider the real Martinet field
$$K={\mathbb{Q}}(\sqrt{2},\sqrt{3\dot 5\dot 7\dot 23\dot 29})$$
of degree $4$ over
${\mathbb{Q}}$. We have
$$
D_K=2^{8}\dot (3\dot 5\dot 7\dot 23\dot 29)^2\;,
$$
$$
g=g(K)=\log\sqrt{\vert D_K\vert}\approx 13.9293...
$$
This field has an infinite unramified $2$-tower ${\mathcal K}$.
Then $$\BS({\mathcal K})\in (\BS_{\lower}({\mathcal K}),\BS_{\upper}({\mathcal K})),$$
where
$$\BS_{\lower}({\mathcal K})=1-{{4\log 2}\over g}$$
and
$$\BS_{\upper}({\mathcal K})=\BS_{\lower}({\mathcal K})+{\log{2}\over{g}}+
{{{\sqrt{7}-1}\over{g\log 7}}
\left(g-2\gamma-{\pi}-2\log{8\pi}-{{\log 2}\over{\sqrt{2}-1}}\right)
}\log{{7}\over 6},$$
i.e., approximately in the interval
$$(0.8009...\; , 0.8648...).$$
The deficiency of this tower is at most
$$
1-{{2\gamma+\pi+2\log{8\pi}}\over{g}}\approx
0.2286...
$$
}

{\em Proof.} We proceed as in the proof of Theorem 9.2. Our field
has  an infinite unramified tower (Corollary 9.1). Let
$K_1={\mathbb{Q}}(\sqrt 2)$,
$K_2={\mathbb{Q}}(\sqrt {3\dot 5\dot 7\dot 23\dot 29})$, $K=K_1\dot K_2$. Since
$3\dot 5\dot 7\dot 23\dot 29=70035\equiv3 (\mod 4)$, we have $D_{K_2}=4\dot
70035$
and 2, 3, 5, 7, 23 and 29 are ramified in $K_2$. Since $D_K=2^8\dot
70035^2$, we
see that the ideal lying over 2 is also ramified in $K/K_2$, i.e., 2 is
totally ramified
in $K/{\mathbb{Q}}$. In $K_1/{\mathbb{Q}}$ only 2 is ramified, $D_{K_1}=8$.
Since 2 is congruent to a
square modulo 7, and noncongruent to a square modulo 3 and 5,
we see that 7 splits, but 3 and 5 remain inert in $K_1/{\mathbb{Q}}$. Thus
in $K/{\mathbb{Q}}$
there is one ideal of norm 2, no ideals of norm 3 and 5, and two ideals of
norm 7. There are
4 real places and no complex ones.

Using, as above, the linear programming approach, with $a_1=a_3=a_5=0$ and
$x_2\le{1\over g}$, we get
$$\BS_{\lower}(K)=1-{{4\log 2}\over g}$$
and
$$\BS_{\upper}(K)=\BS_{\lower}(K)+\max F(x)\;,$$
where $F(x)=\sum\limits_{q\not=3,5}b_q x_q$.

The maximum is attained for $x_2={1\over g}$, $x_q=0$ for $q=4$ and $q>7$,
and the value of $x_7$ is chosen so that (iii) becomes an equality.$\;{\Box}$

{\bf Remark 9.2.} The other two fields of Corollary 9.1 give the following
numerical results. For
$$K={\mathbb{Q}}(\sqrt{-3\dot 5\dot 17\dot 19})$$
we GRH--have
$$g(K)\approx 4.9359...,$$
$$\delta({\mathcal K})\le 0.2298...,$$
and
$$\BS_{\lower}({\mathcal K})\approx 0.6276..., \; \BS_{\upper}({\mathcal K})\approx 0.6402...$$
For
$$K={\mathbb{Q}}(\sqrt{3\dot 5\dot 13\dot 29\dot 61})$$
we GRH--have
$$g(K)\approx 7.0687...,$$
$$\delta({\mathcal K})\le 0.2400...,$$
and
$$\BS_{\lower}({\mathcal K})\approx 0.8038..., \; \BS_{\upper}({\mathcal K})\approx 0.9020...$$

{\bf GRH Theorem 9.4.} {\em Consider the totally real quadratic field
$$K = {\mathbb{Q}}(\sqrt{11\dot 13\dot 17\dot 19\dot 23\dot 29\dot 31\dot
37\dot
41\dot 43\dot 47\dot 53\dot 59\dot 61\dot 67}).$$
The genus of this field is $g\approx 25.9882\dots$.
This field has an infinite unramified $2$-tower ${\mathcal K}$ in which nine prime
ideals lying over $2$, $3$, $5$, $7$ and $71$ split completely.
Then $\BS({\mathcal K})\in (\BS_{\lower}({\mathcal K}),\BS_{\upper}({\mathcal K}))$ 
and
${\mathbf\kappa}({\mathcal K})\in({\mathbf\kappa}_{\lower}({\mathcal
K}),{\mathbf\kappa}_{\upper}({\mathcal K}))$, where
$$\BS_{\lower}({\mathcal K})=1+{{2\log{3\over 2} +2\log{5\over 4}+2\log{7\over
6}+\log{5041\over 5040}}\over g},$$
$$\BS_{\upper}({\mathcal K})=\BS_{\lower}({\mathcal K})+{1\over g}\sum_{p=11}^{47}\log{p\over p-1}+$$
$${{{\sqrt{53}-1}\over{g\log 53}}
\left(g-\gamma-{\pi\over 2}-\log{8\pi}-2\sum_{p=2}^7 {{\log p}\over{\sqrt p
-1}}-
{{\log 71^2}\over{70}}-\sum_{p=11}^{47}
{{\log p}\over{\sqrt p -1}}\right)}\log{{53}\over52},$$
$${\mathbf\kappa}_{\lower}({\mathcal K})={2\log 2+{2\log{3\over 2} +2\log{5\over 4}+2\log{7\over
6}+\log{5041\over 5040}}\over g},$$
$${\mathbf\kappa}_{\upper}({\mathcal K})=
\BS_{\upper}({\mathcal K})-1+{2\log 2\over g},$$
the sums being taken over prime $p$'s.
Numerically
$$\BS({\mathcal K})\in(1.0602\dots , 1.0798\dots),$$
$${\mathbf\kappa}({\mathcal K})\in(0.1135\dots , 0.1331\dots).$$}

{\em Proof.} Let
$d=11\dot 13\dot 17\dot 19\dot 23\dot 29\dot 31\dot 37\dot
41\dot 43\dot 47\dot 53\dot 59\dot 61\dot 67.$
An easy, though tedious check shows that $d$ is congruent to 1 modulo 8, and
it is a square modulo  3, 5
and 7, but not modulo 71. Hence 2, 3, 5, 7 split in $K/{\mathbb{Q}}$ and 71
is inert.
Corollary 9.3 shows that there is an unramified tower splitting the nine
ideals over
2, 3, 5, 7 and 71. Then we use the same linear programming approach of
Subsection
7.1 to calculate $\BS_{\lower}(K)$, $\BS_{\upper}(K)$,
${\mathbf\kappa}_{\lower}({\mathcal K})$ and ${\mathbf\kappa}_{\upper}({\mathcal
K})$.${\Box}$

{\bf GRH Theorem 9.5.} {\em Consider the totally complex quadratic field
$$K={\mathbb{Q}}(\sqrt{-13\dot 17\dot 19\dot 23\dot 29\dot 31\dot 37\dot
41\dot 43\dot 47\dot 53\dot 59\dot 61\dot 73\dot 79}).$$
This field has an infinite unramified $2$-tower ${\mathcal K}$ in which ten prime
ideals lying over $2$, $3$, $5$, $7$ and $11$ split comletely. The genus of
this
field is $g\approx 27.0169...$
Then $\BS({\mathcal K})\in (\BS_{\lower}({\mathcal K}),\BS_{\upper}({\mathcal K}))$
and ${\mathbf\kappa}({\mathcal K})\in ({\mathbf\kappa}_{\lower}({\mathcal
K}),{\mathbf\kappa}_{\upper}({\mathcal K}))$, where
$$\BS_{\lower}(K)=
1-{1\over g}\log{2\pi}+{2\over g}\left({{\log 2+\log(3/2)+\log(5/4)+\log(7/6)
+\log(11/10)}}\right),$$
$$\BS_{\upper}({\mathcal K})=\BS_{\lower}({\mathcal K})+{1\over g}\sum_{p=13}^{61}\log{p\over p-1}$$
$$+{{{\sqrt{67}-1}\over{g\log 67}}
\left(g-\gamma-\log{8\pi}-2\sum_{p=2}^{11} {{\log p}\over{\sqrt p -1}}-
\sum_{p=13}^{61}
{{\log p}\over{\sqrt p -1}}\right)}\log{{67}\over66}\;\;,$$
$${\mathbf\kappa}_{\lower}(K)=
{2\over g}\left({{\log 2+\log(3/2)+\log(5/4)+\log(7/6)
+\log(11/10)}}\right),$$
$${\mathbf\kappa}_{\upper}(K)=\BS_{\upper}({\mathcal K})-1+{1\over g}\log{2\pi},$$
the sums being taken over prime $p$.
Numerically
$$(\BS_{\lower}({\mathcal K}),\BS_{\upper}({\mathcal K}))=(1.0482... , 1.0653...),$$
$$({\mathbf\kappa}_{\lower}({\mathcal K}),{\mathbf\kappa}_{\upper}({\mathcal
K}))=(0.1162... , 0.1333...).$$}

{\em Proof.} Along the same lines as the proofs of Theorems 9.4, 9.3 and 9.2,
using Corollary 9.4. ${\Box}$

{\bf Remark 9.3.} Ihara's example of Corollary 9.5
$$K={\mathbb{Q}}(\sqrt{-3\dot 5\dot 7\dot 11\dot 13\dot 17\dot 23\dot 31})$$
with two divisors of 2 splitting in the tower,
has $g\approx 9.5097...$, its deficiency $\delta$ is at most 0.2483... and
$(\BS_{\lower}({\mathcal K}),\BS_{\upper}({\mathcal K}))=(0.9525..., 1.010...)$.

Let us see what can be got without GRH. We consider the same
fields as in GRH Theorems 9.2 and 9.3.

{\bf Theorem 9.7.} {\em Consider the Martinet field
$$K={\mathbb{Q}}(\cos{2\pi\over11},\sqrt{2},\sqrt{-23})$$
of degree $20$ over
${\mathbb{Q}}$.
This field has an infinite unramified $2$-tower ${\mathcal K}$, and $\BS({\mathcal K})\in
(\BS_{\lower}({\mathcal K}),\BS_{\unc,\upper}({\mathcal K}))$, where
$$\BS_{\lower}({\mathcal K})=1-{{10\log(2\pi)}\over g}$$
and
$$\BS_{\unc, \upper}({\mathcal K})=\BS_{\lower}({\mathcal K})+{{10\log({23\over 22})}\over g}
+{{2\log({32\over 31})}\over g}+{{20}\over g}\sum\limits_{p=37}^{97} \log{p\over p-1},$$ i.e.,
approximately in the interval
$$(0.5939..., 0.7108...).$$}

{\em Proof.} We proceed along the same lines as before. Having no GRH at
hand, instead of
GRH Theorem 7.2 we use Theorem 7.3 (the tower being almost normal, as any
2-tower over a normal field), and instead of GRH
Theorem 3.1 we use
either Proposition 3.1 or Proposition 3.2. The latter is easier to
calculate. Knowing the
decomposition law for small primes (cf. the proof of GRH Theorem 9.2), we
see that in $K$
there are 10 infinite complex places, 10 places whose norm is 23, 2 places
of norm 32 and no
other places of norm strictly smaller than 37. Over any other prime
there are at most 20 places.

Then we use the optimization procedure of Section 8, that shows that to get
an upper bound we
can exaggerate the number of places with small norms. Suppose that there
were 20 ideals of each of the norms
from 37 to 97 (in fact, there are much less). Even this would contradict the
inequality of Proposition 3.2,
i.e.,
$${10\over g}(\gamma+\log{2\pi})+{{10}\over g}{\log{23}\over 22}+{2\over
g}{\log{32}\over 31}
+{20\over g}\sum_{p=37}^{97}{\log{p}\over p-1}> 1\;,$$
the sum being taken over primes.
Therefore, by the inequality of Theorem 7.1,
any limit point of the Brauer--Siegel ratio is at most
$$1-{10\over g}\log{2\pi}+{{10}\over g}{\log{{23}\over 22}}+{2\over
g}{\log{{32}\over 31}}
+{20\over g}\sum_{p=37}^{97}{\log{{p}\over p-1}}\approx 0.7108...\;\;{\Box}$$

{\bf Remark 9.4.} Using Proposition 3.1 instead of Proposition 3.2 we can
do better.
We can also use further information on prime decomposition
in $K/{\mathbb{Q}}$.  (In particular, the only possible norms between
37 and 1000 are in fact 121, 353,
439, 463, 593, 967 and 991.) This makes
the constant better. Namely, we can prove that
$\BS({\mathcal K})\le 0.623\dots$ 
  
{\bf Theorem 9.8.} {\em Consider the real Martinet field
$$K={\mathbb{Q}}(\sqrt{2},\sqrt{3\dot 5\dot 7\dot 23\dot 29})$$
of degree $4$ over
${\mathbb{Q}}$. Its genus equals
 
$$
g=g(K)=\log\sqrt{\vert D_K\vert}\approx 13.9293...,
$$
it has an infinite unramified $2$-tower ${\mathcal K}$, and
  $\BS({\mathcal K})\in (\BS_{\lower}({\mathcal K}),\BS_{\unc, \upper}({\mathcal K}))$,
where

$$\BS_{\unc,\upper}({\mathcal K})=\BS_{\lower}({\mathcal
K})+ {\log{2} +2\log{7\over 6} +4\log{11\over 10} +4\log{13\over 12} \over g}+$$
 $$
{{1\over {2gA_{17}}}
\left(g-2\gamma-{\pi}-2\log{8\pi}-A_2-2A_7-4A_{11}-4A_{13}\right)
}\log{17\over 16}\;,$$
where $ A_p={2\log p{\sum_{m=1}^{\infty}(p^m+1)^{-1} }}$,
i.e., approximately in the interval
$$(0.8009\dots , 0.9248\dots).\Box$$
}
{\bf Remark 9.4.} Applying the same technique to the fields of GRH Theorems 9.4 and 9.5 we get the
results presented in the table at the end of Section 1. We do not write out here
the exact
formulae which are rather cumbersome.

\section{Open questions} 
In this section we discuss some open questions concerning the Generalized
Brauer--Siegel Theorem.
First of all, in the proof of the Generalized  Brauer--Siegel Theorem we do not really use
the whole  strength of GRH; moreover, under some mild conditions we have totally dispensed with
GRH (Theorem 7.3).  Therefore, it is but
 natural to ask whether one really needs GRH to prove the result, which leads to 
 
{\bf Problem 10.1.} {\em Prove GRH Theorem} 7.2 {\em unconditionally.}

 Let us now discuss some problems, connected with the Brauer--Siegel ratio introduced
 and studied above. First of all, its very existence (i.e., the existence of the
 corresponding limit) is proved only 
 under GRH or for almost normal asymptotically good infinite global fields, which leads to
 
{\bf Problem 10.2.} {\em Prove unconditionally that for any asymptotically exact family
${\mathcal K}$  of number fields  the Brauer--Siegel ratio $\BS({\mathcal K})$ is well defined. }
 
The following problem is connected with the fact that for an
arbitrary asymptotically exact family unconditionally we have only an  {\em upper}
bound for
$\BS({\mathcal K})$, cf. Theorems 7.1, 7.3 and 8.2.
 
{\bf Problem 10.3.} {\em Give an  unconditional lower bound   for  the Brauer--Siegel
 ratio $\BS({\mathcal K})$ for {\rm any} asymptotically exact family.}
 
Note that for towers of normal number fields this results from Theorems 7.3 and
8.2. One can hope that this problem can be solved if one ameliorates 
  the technique of the usual proof of the Brauer--Siegel theorem, i.e., estimates
in the (adelic) integral representation of the zeta-function (cf. Lemma 3 of Section XVI.2 
of \cite{13}). 

 There also is the question of how good our bounds and examples are.
 
{\bf Problem 10.4.}   {\em Ameliorate on the bounds of GRH Theorem 8.1 and/or of
Theorem 8.2.} 
 
{\bf Problem 10.5.} {\em Construct examples of class field towers  $($or other 
asymptotically exact families$)$  with $\BS({\mathcal K})$ GRH--smaller than those of GRH Theorems 6.2
 and 9.3 or GRH--greater than those of GRH Theorems 9.4 and 9.5. }

Our results in the present paper  are of an  asymptotic nature. However, it is clear that 
a good part of them can be made effective which leads to 
 
{\bf Problem 10.6.} {\em Give effective versions of the above results  with the
remainder terms  as good as possible. }

M.Ts.:

Institut de Math\'ematiques de Luminy, UPR 9016 du CNRS, 

Case 907, 13288, Marseille, FRANCE, 

Independent University of Moscow, and

Dobrushin Math. Lab., Institute for Problems of Information Transmission, Russian
Academy of Sciences. 

E-mail: tsfasman@iml.univ-mrs.fr

\vskip 0.3 cm
S.Vl.:
Institut de Math\'ematiques de Luminy, UPR 9016 du CNRS,
 
Case 907, 13288, Marseille, FRANCE, and 

Dobrushin Math. Lab., Institute for Problems of Information Transmission, Russian
Academy of Sciences. 

E-mail: vladut@iml.univ-mrs.fr.

\end{document}